\documentclass[onefignum,onetabnum]{siamart171218}
\pdfoutput=1

\usepackage{subcaption} %
\usepackage{amssymb} %
\usepackage{mathtools} %
\usepackage{graphicx}
\usepackage{color}
\usepackage{tabularx}
\usepackage{tikz} %
\usepackage{enumitem}
\usepackage[labelfont=bf]{caption}
\usepackage{changepage}
\usepackage{afterpage} %
\usepackage{gensymb} %
\usepackage{placeins} %
\usepackage{enumitem}%
\usepackage{siunitx} %
\usepackage{soul} %
\usepackage[export]{adjustbox} %
\usepackage{tabularx}
\usepackage{wrapfig}
\usepackage[export]{adjustbox}
\usepackage{mwe}
\usepackage{dashrule}

\usepackage{lipsum}
\usepackage{amsfonts}
\usepackage{graphicx}
\usepackage{epstopdf}
\usepackage{algorithmic}
\ifpdf
  \DeclareGraphicsExtensions{.eps,.pdf,.png,.jpg}
\else
  \DeclareGraphicsExtensions{.eps}
\fi

\newsiamremark{remark}{Remark}
\newsiamremark{hypothesis}{Hypothesis}
\crefname{hypothesis}{Hypothesis}{Hypotheses}
\newsiamthm{claim}{Claim}

\headers{Constructing Polynomial Block Methods}{T. Buvoli}

\title{Constructing Polynomial Block Methods \thanks{Submitted to the editors DATE.
\funding{This work was funded by the National Science Foundation, Computational Mathematics Program, under DMS-1115978 and DMS-1216732}}}

\author{Tommaso Buvoli\thanks{Department of Applied Mathematics, University of Washington, Seattle, WA USA. 
  (\email{buvoli@uw.edu}).}
}

\usepackage{amsopn}

\ifpdf
\hypersetup{
  pdftitle={Constructing Polynomial Block Methods},
  pdfauthor={T. Buvoli}
}
\fi

\newcommand{\ODEd}{ODE dataset}
\newcommand{\ODED}{ODE Dataset}

\newcommand{\ODEp}{ODE polynomial}
\newcommand{\ODEP}{ODE Polynomial}

\newcommand{\ODEps}{ODE polynomials}

\newcommand{\ODEdp}{ODE derivative polynomial}

\newcommand{\ODEdps}{ODE derivative polynomials}

\newcommand{\iv}{interpolated value}

\newcommand{\ivs}{interpolated value set}

\definecolor{plot_black}{RGB}{10, 10, 10}
\definecolor{plot_red}{RGB}{192, 57, 43}
\definecolor{plot_orange}{RGB}{230, 126, 34}
\definecolor{plot_yellow}{RGB}{241, 196, 15}
\definecolor{plot_green}{RGB}{39, 174, 96}
\definecolor{plot_blue}{RGB}{41, 128, 185}
\definecolor{plot_violet}{RGB}{155, 89, 182}
\definecolor{plot_grey}{RGB}{149, 165, 166}

\begin{document}

\maketitle

\begin{abstract}

The recently introduced polynomial time integration framework proposes a novel way to construct time integrators for solving systems of first-order ordinary differential equation by using interpolating polynomials in the complex time plane. In this work we continue to develop the framework by introducing several additional types of polynomials and proposing a general class of construction strategies for polynomial block methods with imaginary nodes. The new construction strategies do not involve algebraic order conditions and are instead motivated by geometric arguments similar to those used for constructing traditional spatial finite differences. Moreover, the newly proposed methods address several shortcomings of previously introduced polynomial block methods including the ability to solve dispersive equations and the lack of efficient serial methods when parallelism cannot be used. To validate our new methods, we conduct two numerical experiments that compare the performance of polynomial block methods against backward difference methods and implicit Runge-Kutta schemes.

\end{abstract}

\begin{keywords}
  Time-Integration, Polynomial Interpolation, Complex Time, High-order, Parallelism.
\end{keywords}

\begin{AMS}
 65L04, 65L05, 65L06, 65E99
\end{AMS}

The study of time integration methods for solving the initial value problem
\begin{align}
	\mathbf{y}' = f(t,\mathbf{y}), \quad \mathbf{y}(t_0) = \mathbf{y}_0,
\end{align}
has a rich history \cite{hairer2008solving, hairer1996solving, butcher2016numerical, iserles2009first, gear1971numerical} that has given rise to a variety of different approaches including extrapolation methods, linear multistep methods (LMM), Runge-Kutta (RK) methods, and, more broadly, general linear methods (GLMs) \cite{butcher2006}. The recently introduced polynomial time integration framework \cite{buvoli2018polynomial, buvoli2019constructing} proposes a new approach for constructing time integrators by utilizing interpolating polynomials in the complex time plane. The use of polynomials simplifies high-order method construction while analytic continuation in the complex time plane provides improved stability that cannot be attainted using classical polynomial linear multistep methods. The extension to complex time is not unique to polynomial methods, and has proven beneficial for several other time integration approaches \cite{corliss1980integrating, fornberg2011numerical, HansenOstermann2009, orendt2009geometry}.

Another key component of the polynomial framework is a parameter known as the extrapolation factor that controls the ratio between the stepsize and the width of the interpolation nodes used to construct polynomials. This quantity can be tuned to significantly improve the linear stability properties of an integrator or to ameliorate its susceptibility to numerical roundoff error.

The introductory work \cite{buvoli2018polynomial} laid the foundation for the polynomial framework and introduced two classes of parallel implicit polynomial block methods (PBM) named BBDF and BAM which generalize classical backwards difference methods (BDF) and Adams-Mouton (AM) schemes. Though the PBMs offer significantly improved stability regions compared to BDF and AM methods, they also suffer from two key weaknesses. First, they are inefficient if they cannot be parallelized. The accuracy of BBDF and BAM is nearly identical to that of BDF and AM; however, in the absence of parallelization, the cost per timestep scales linearly with order. Secondly, neither BBDF or BAM methods of order greater than two possess linear stability regions that encompass the imaginary axis. This implies that neither class of method can be used to solve purely dispersive equations. Moreover, the BAM method had a bounded stability region that prohibits the method from being used to solve highly stiff problems.

In this paper we continue to develop the polynomial framework and then demonstrate its flexibility by constructing multiple families of new integrators that address the aforementioned weaknesses of BBDF and BAM. Since the space of polynomial methods is large, we focus exclusively on polynomial block methods with imaginary nodes. When discussing method construction, we emphasize the ability to design integrators by solely considering the underlying interpolating polynomials. We accomplish this by introducing stencil diagrams that are inspired by spatial finite difference stencils and can be used to visualize the polynomials that make up a method.

The primary goals of this work are threefold:
	(1) introduce a broad set of construction strategies for polynomials block methods that include both parallel and serial methods, (2) investigate the possibility of using polynomial methods for dispersive equations, and (3)  briefly introduce the Matlab package PIPack \cite{buvoli2020PIPACK} that provides a simple way to compute linear stability properties and initialization coefficients for all the families of polynomial methods discussed in this work. 
	
This paper is organized as follows. In the first section we provide a brief overview of the polynomial time integrators, describe two new \ODEps{} and the \ivs{}, and introduce several new stencil diagrams.  In Section \ref{sec:constructing_pbms} we introduce a range of strategies for constructing polynomial block methods that posses a variety of different properties including high-order and stability along the imaginary axis. Next, in Sections \ref{sec:linear_stability} and \ref{sec:pipack} we discuss linear stability and briefly mention our Matlab package. Finally in Section \ref{sec:numerical_experiments} we conduct two numerical experiments to validate the efficiency and accuracy for the newly proposed implicit PBMs.

\section{Polynomial time integrators}

The polynomial time integration framework introduced in \cite{buvoli2018polynomial, buvoli2019constructing} encompasses a wide range of methods whose stages and outputs are computed using interpolants named {\em \ODEps{}}. \ODEps{} are a general class of polynomials that approximate the local Taylor series of the ODE solution and encompass a range of well-known approximations including classical Lagrange and Hermite interpolating polynomials. All \ODEps{} are constructed using data from an {\em \ODEd{}} that contains the method's inputs, outputs, stages and their derivatives.

Every polynomial method accepts inputs $y_j^{[n]}$ and produces outputs $y_j^{[n]}$ where ${j = 1, \ldots, q}$. The input and  output nodes of a method scale with respect to the {\em node radius} $r$, and each input or output value approximates the differential equation solution $y(t)$ at a particular time so that
	\begin{align*}
		y_j^{[n]} \approx y\left(t_n + r z_j\right) \quad \text{and} \quad y_j^{[n+1]} \approx y\left(t_n + r z_j + h\right)
	\end{align*}
	where the set $\{z_j\}_{j=1}^q$ is called the node set. A polynomial method may also produce stage values $Y_i$, $i = 1, \ldots, s$, however we will not be discussing methods with stages in this paper.

 Polynomial methods are described in terms of the node radius $r$ and the {\em extrapolation factor} $\alpha$ (rather than $r$ and $h$) since these are natural variables for working with polynomials in local coordinates. The relationship between these three quantities is simply $h = r\alpha$. During time stepping the node radius $r$ serves as the stepsize and $\alpha$ is a constant that parametrizes every polynomial time integrator. A method with a large $\alpha$ has a stepsize that is comparatively larger than the length of the interval containing the interpolating nodes, while the opposite is true for a method with small $\alpha$. By adjusting the extrapolation factor one can fundamentally alter the properties of the underlying method.

In the subsections below, we briefly review the \ODEd{}, the \ODEp{}, and the notation used to describe polynomial time integrators. We also extend the polynomial framework presented in \cite{buvoli2019constructing} by introducing the GBDF \ODEp{}, the \ODEdp{}, and the \ivs{}. We close by briefly reviewing polynomial block methods and introducing several new stencil diagrams for visualizing the methods and their \ODEps{}.

\subsection{The \ODEd{}}

The \ODEd{} is an ordered set containing a method's inputs, outputs and stages at the $n$th timestep. An \ODEd{} of size $w$ can be denoted as
\begin{align*}
	D(r, t_n) =\left\{ \left(\tau_j,~ y_j,~ r f_j \right) \right\}_{j=1}^w
	\quad \text{where} \quad
	y_j \approx y(t(\tau_j)) \quad \text{and} \quad f_j = F(t(\tau_j), y_j).
\end{align*}
The data is represented in the local coordinates $\tau$ where the global time is
	\begin{align}
		t(\tau) = r\tau + t_n,
		\label{eq:global_time}
	\end{align}
	and the scaling factor $r$ is the node radius.
	
\subsection{The \ODEP{}}

The general form for an \ODEp{} of degree $g$ with expansion point $b$ is
			\begin{align}
				p(\tau; b) &= \sum_{j=0}^{g}  \frac{a_{j}(b)(\tau - b)^j}{j!}
				\label{eq:solution_ode_polynomial}
			\end{align}
			where each approximate derivative $a_{j}(b)$ is computed by differentiating interpolating polynomials constructed from the values in an \ODEd{} $D(r, t_n)$. A general formulation for the approximate derivatives $a_j(b)$ is described in \cite{buvoli2018polynomial, buvoli2019constructing}. In this work, we  restrict ourselves to the following three special families of \ODEps{}:

			\vspace{0.5em}
			\begin{enumerate}
				\item {\em Adams \ODEps{}} can be written as the sum of a Lagrange interpolating polynomial $L_y(\tau)$ that approximates $y(t)$ and the integral of a Lagrange interpolating polynomial $L_F(\tau)$ that approximates $ry'(t)$. When expressed in integral form, an Adams \ODEp{} takes the form
				\begin{align}
					p(\tau;b) = L_y(b) + \int_{b}^\tau L_F(\xi) d\xi
				\label{eq:adams_poly_integral_form}
				\end{align}
				where the expansion point $b$ now acts as the left integration endpoint.
				
				\vspace{0.5em}
				\item {\em Generalized BDF (GBDF) \ODEps{}} are polynomials 
				\begin{align}
					p(\tau; b) = H_y(\tau)
					\label{eq:gbdf_poly}
				\end{align}
				where $H_y(\tau)$ is the polynomial of least degree that interpolates one or more solution values and whose derivative $H_y'(\tau)$ interpolates one or more derivative values. We may express these conditions mathematically as
					\begin{align*}
						H_y(\tau_j) &= y_j,  && j \in \mathcal{A}, \\
						H'_y(\tau_k) &= rf_k, && k \in \mathcal{B}
					\end{align*}
					where $\mathcal{A}$ and $\mathcal{B}$ are non-empty sets containing indices ranging from $1$ to the size of the underlying \ODEd{} $w$. All GBDF \ODEps{} do not depend on the parameter $b$. 
				\item {\em BDF \ODEps{}} are GBDF polynomials where $H'_y(\tau)$ interpolates only one derivative value.
			\end{enumerate}
			
			Polynomial integrators can be broadly classified by the families of \ODEps{} used to compute their outputs. An integrator that only makes use of one family of \ODEp{} will inherit the name. For example, the BAM and BBDF methods from \cite{buvoli2019constructing} are respectively Adams and BDF type polynomial methods since they only use these types of polynomials.

	\subsection{\ODEdp{}}
	
	To construct certain explicit polynomial integrators we will also make use of {\em \ODEdps{}} that approximate the solution derivative $y'(t)$ and are denoted using the notation $\dot{p}(\tau;b)$. In the same spirit as \ODEps{}, \ODEdps{} are built from an \ODEd{} $D(r,t_n)$ and are expressible as local Taylor series for $ry'(t)$ where the derivatives have been replaced with approximations obtained by differentiating interpolating polynomials.
	
	 The most general form for the \ODEdp{} is identical to that of the \ODEp{} (\ref{eq:solution_ode_polynomial}), however the approximate derivative $a_{j}(b)$ must be computed differently. For details regarding the general construction of \ODEdps{} we refer the reader to \cite{buvoli2018polynomial}. 
	 
	 In this work we will only use one simple family named the Adams \ODEdp{}. The polynomial $\dot{p}(\tau;b)$ is an Adams \ODEdp{} if
		\begin{align}
			\dot{p}(\tau;b)	= L_F(\tau), \hspace{1 em} \forall b
			\label{eq:der_adams_poly_integral_form}
		\end{align}
		where $L_F(x)$ is a Lagrange interpolating polynomial of degree $g$ that interpolates through at least one derivative value in the \ODEd{}.

	\subsection{Interpolated values}
	
	When constructing methods it is possible that solution or derivative data is not available at a particular temporal node. This can be remedied by using known data to obtain a value at a particular point. We formalize this idea by introducing {\em interpolated values} and the {\em \ivs{}}.
	
	 An \iv{} is an approximate solution or derivative generated by evaluating an \ODEp{} or \ODEdp{} constructed from an \ODEd{} $D(r,s)$. Given a non-empty \ODEd{}, it is possible to generate any number of different interpolated values. Together, these approximations form an  
	\ivs{} $I(r,s)$ constructed from $D(r,s)$.
	
	Every interpolated value can be written as a linear combination of the $y_j$ and $f_j$ values in the generating dataset. Therefore, the operation of computing interpolated values does not introduce new information, but it does provide a convenient way to expand the number of possible polynomial approximations that can be constructed from a fixed amount of data.

\subsection{Polynomial block methods} 

The primary aim of this work is to introduce new construction strategies for polynomial block methods (PBM), which we briefly review here. A PBM is a multivalued method that produces $q$ outputs at each timestep. PBMs have no stages and the output values are each determined by evaluated an \ODEp{} constructed using the methods inputs and outputs. Every polynomial block method can be written as	
	\begin{align}
		y^{[n+1]}_j &= p_j(z_j + \alpha;~ b_j), & j & =1, \ldots,q,
		\label{eq:polynomial_block_method}
	\end{align}
	where each $p_j(\tau;b)$ is an \ODEp{} constructed from the \ODEd{}
	\begin{align*}
		D(r,t_n) = \underbrace{ \left\{ \left(z_j,~ y^{[n]}_j,~ rf^{[n]}_j\right) \right\}_{j=1}^q}_{\text{inputs}} \bigcup \underbrace{ \left\{ \left(z_j + \alpha,~ y^{[n+1]}_j,~ rf^{[n+1]}_j\right)\right\}_{j=1}^q }_{\text{ outputs}},
	\end{align*}	
	and any \ivs{} $I(r,t_n)$ generated from $D(r,t_n)$. All polynomial block methods can be written in coefficient form as 
\begin{align}
	\mathbf{y}^{[n+1]} = \mathbf{A}(\alpha) \mathbf{y}^{[n]} + r \mathbf{B}(\alpha) \mathbf{f}^{[n]} + \mathbf{C}(\alpha) \mathbf{y}^{[n+1]} + r \mathbf{D}(\alpha) \mathbf{f}^{[n+1]},
	\label{eq:block_parametrized}
\end{align}
where $\mathbf{A}(\alpha)$, $\mathbf{B}(\alpha)$, $\mathbf{C}(\alpha)$, $\mathbf{D}(\alpha)$ are $q \times q$ coefficient matrices, and the input and output vector and their derivatives are
\begin{align*}
	\mathbf{y}^{[n]} = \left[ y_1^{[n]},~ \ldots,~ y_q^{[n]}\right]^{\text{T}} \quad
	\mathbf{f}^{[n]} = \left[ f_1^{[n]},~ \ldots,~ f_q^{[n]}\right]^{\text{T}}.
\end{align*}

\subsection{Stencils \& diagrams for visualizing \ODEps{}} 
\label{subsection:stencil_diagrams}

To enable the simple visualization of \ODEps{} we introduce several stencil-based diagrams that are generalizations of spatial finite difference stencils. These illustrations showcase the geometric properties of an \ODEp{} and aid us in the presentation of method construction strategies. For the purposes of method construction, we focus on the following three properties of \ODEps{}:
\begin{enumerate}
	\item {\em The active node set.} The active node set for an \ODEp{} $p(\tau; b)$ is the set of unique temporal nodes $\tau_i$ of the data values used to construct $p(\tau; b)$. Moreover, a temporal node  of an \ODEd{} is said to be {\em active} if the corresponding solution or derivative data is a member of the active node set.	
	\item {\em The expansion point} (denoted by $b$ in (\ref{eq:solution_ode_polynomial})). This property is particularly important for Adams \ODEps{} (\ref{eq:adams_poly_integral_form}) since it represents the left integration endpoint.	Conversely, this quantity is not useful for BDF and GBDF polynomials since they do not depend on the expansion point.
	\item {\em The polynomials for computing approximate derivatives. } Every \ODEp{} has approximate derivatives $a_j(b)$ which are computed by differentiating interpolating polynomials. For a GBDF polynomial (\ref{eq:gbdf_poly}), all approximate derivatives are determined by differentiating $H_y(\tau)$. For an Adams \ODEp{} (\ref{eq:adams_poly_integral_form}) the approximate derivatives are
		\begin{align*}
			 a_0(b) = L_y(b) \quad \text{and} \quad a_{j\ge 1}(b) = \left. \frac{d^{j-1}}{d\tau^{j-1}} L_F(\tau) \right|_{\tau=b}.	
		\end{align*} %
	 Therefore, the stencils for GBDF integrator must show $H_y(\tau)$, while stencils for Adams integrators must show $L_y(\tau)$ and $L_F(\tau)$.
\end{enumerate}

\noindent	For each of these three properties we introduce the corresponding stencil diagrams:
			\begin{enumerate}
				\item A {\em node stencil} illustrates the active and inactive nodes for an \ODEp{}, and a {\em node diagram} displays all the active nodes stencils for the \ODEps{} of a method.  	
				\item An {\em expansion point stencil} shows the expansion points of an \ODEp{}. For Adams \ODEps{} they can also show the integration path connecting the expansion point to the evaluation point. An {\em expansion point diagram} overlays all the expansion point stencils into one figure. 
				\item An {\em ODE polynomial diagram} can be used to visualize the interpolating polynomials for computing approximate derivatives. By combining multiple diagrams we can describe all the \ODEps{} of a method. 
			\end{enumerate}
			
We provide a more detailed explanation of these new diagrams and stencils in the subsections below. In this work we only discuss diagrams for Adams and GBDF \ODEps{}; however in \cite{buvoli2018polynomial}, these diagrams are also extended to more general \ODEps{} (\ref{eq:solution_ode_polynomial}). 

	\subsubsection{The node stencil and diagram}

		The node stencil for an \ODEp{} shows the active and inactive temporal nodes of the corresponding \ODEd{}. These stencils are useful for classifying \ODEps{} that share the same active node set but use different data values at each node.
		
\begin{figure}[h]				
				\centering\includegraphics[width=0.8\linewidth]{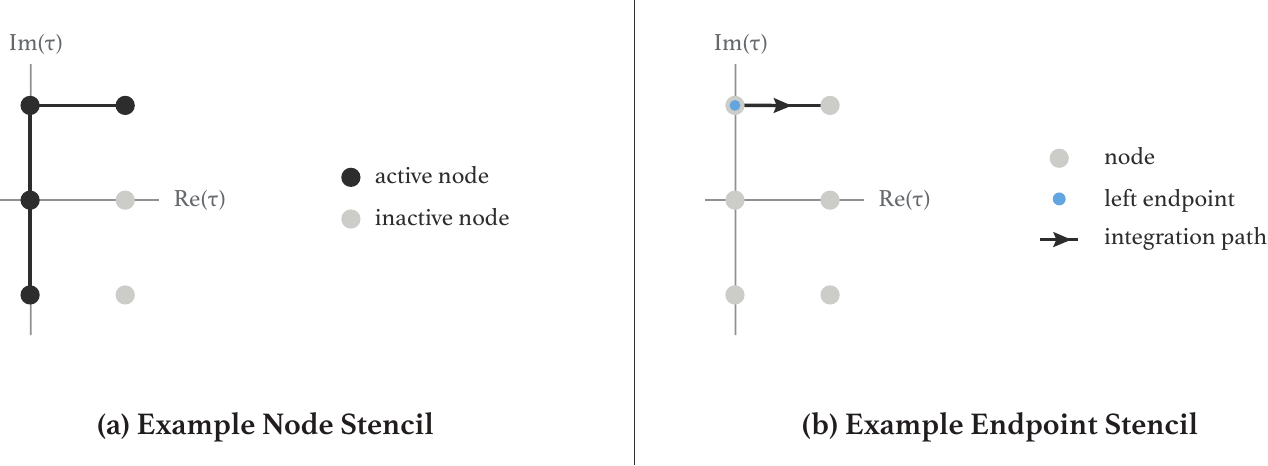}	
			\caption{An example node stencil and endpoint stencil for an \ODEp{} constructed using an \ODEd{} with temporal nodes $\{\tau_j\} = \{-i, 0, i, -i+1, 1, i+1\}$. {\bf (a)} The node stencil for an \ODEp{} that is constructed using the first four solution or derivative values from the \ODEd{}. {\bf (b)} The endpoint stencil for an Adams \ODEp{} that is being evaluated at $\tau = i + 1$ and has a left integration point at $\tau = i$.}
			\label{fig:example_stencils}
		\end{figure}

		Each stencil shows the complex $\tau$-plane where the global time $t(\tau) = r\tau + t_n$. The axes correspond to the real and imaginary components of $\tau$ and the origin marks the point $\tau = 0$. We color the active nodes black and inactive nodes light gray. In Figure \ref{fig:example_stencils} we present an example node stencil for an \ODEp{} with complex nodes.
				
\begin{figure}[h]				
				\centering\includegraphics[width=0.8\linewidth]{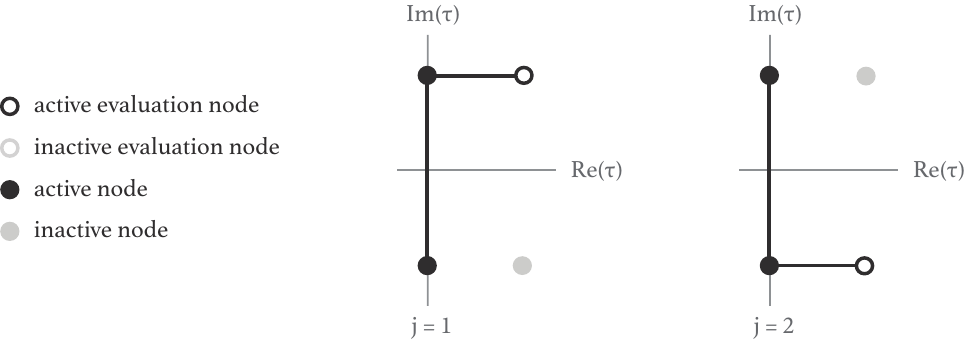}	
			\caption{An example node diagram for a polynomial block method where: $q = 2$, $\{z_j\} = \{i, -i\}$, and the two \ODEp{} are respectively constructed using the data at the nodes $\{i,-i,i + \alpha\}$ and $\{i, -i, -i + \alpha\}$. The BBDF and BAM methods with $q=2$ from \cite{buvoli2019constructing} are two examples methods with this node diagram.}
			\label{fig:example_node_diagram}
		\end{figure}

		To visualize a polynomial block method, we create an active node diagram that combines the $q$ active node stencils for each of the method's \ODEps{}. For additional clarity we distinguish the {\em evaluation point} $\tau = z_j + \alpha$ of the $j$th \ODEp{} from the other temporal nodes using a hollowed circle marker. An example diagram is shown in Figure \ref{fig:example_node_diagram}.

		\subsection{The expansion-point stencil and diagram}	
			An expansion-point stencil for an \ODEp{} shows the temporal nodes of the underlying \ODEd{} along with the expansion point $b$. Expansion-point stencils for the \ODEps{} of an Adams block methods also show the integration path that connects the expansion point (also called the left endpoint for Adams polynomials) to the evaluation node $z_j + \alpha$. The nodes are shown in the complex $\tau$-plane as grey circles, the endpoint is depicted using a blue circle, and the integration paths are drawn with a directed line. 
			Expansion point stencils are useful for presenting a variety of different expansion points $b$ without needing to discuss the underlying \ODEp{}. In Figure \ref{fig:example_stencils} we show an example endpoint stencil for an Adams \ODEp{} with complex nodes. By overlying $q$ expansion point stencils, we create an expansion-point diagram that describes all the endpoint choices for the polynomials comprising an Adams block method.

	\subsubsection{The \ODEp{} diagram}
		
		An \ODEp{} diagram contains multiple stencils that describe each of the interpolating polynomials used to compute approximate derivatives. For an GBDF  \ODEp{} (\ref{eq:gbdf_poly}) we need only a single stencil describing $H_y(\tau)$, while for an Adams \ODEp{} (\ref{eq:adams_poly_integral_form}) we require two stencils to describe the Lagrange interpolating polynomials $L_y(\tau)$ and $L_F(\tau)$. 
			
\begin{figure}[h]
		\begin{center}
			\includegraphics[width=\linewidth]{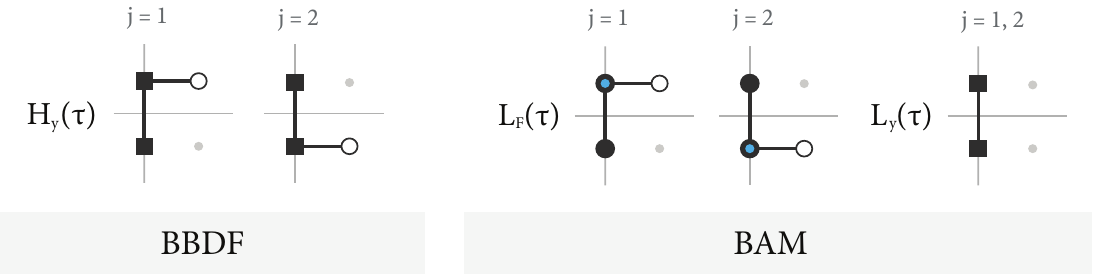}
		\end{center}
		\caption{Example polynomial diagrams for the BBDF and BAM methods with $q=2$. The index $j$ represents the output index of the polynomial $p_j(\tau; b)$.}
		\label{fig:example_poly_diagrams}
	\end{figure}

		The nodes are shown in the complex $\tau$-plane and labeled depending on: (1) the type of data used at the node (solution, derivative, or both), and (2) whether the data is an input or output value from \ODED{} or an interpolated value. The diagrams can also show inactive nodes and expansion points. The following table shows each of the markers:

		\vspace{0.5em}
		\begin{center}
			\includegraphics[width=0.85\linewidth]{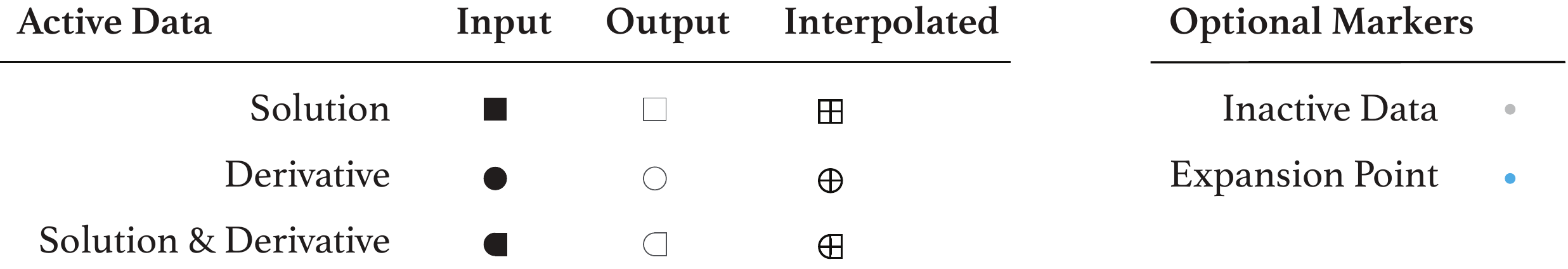}
		\end{center}
		\vspace{0.5em}
	Finally, we can combine $q$ \ODEp{} diagrams to visualize a polynomial block method. The resulting diagram is called a {\bf polynomial method diagram}. In Figure \ref{fig:example_poly_diagrams} we show an example polynomial method diagrams for the BAM and BBDF methods with $q=2$.
			
\section{Constructing polynomial block methods}
\label{sec:constructing_pbms}

	We now present a geometric approach for constructing Adams, BDF, and GBDF polynomial block methods with imaginary nodes. Each of the proposed families of PBMs can  be implemented at any order of accuracy and using multiple node sets. Moreover, the new integrators can also be either implicit or explicit and parallel or serial.
	
	We begin our discussion by restricting ourselves to imaginary nodes that are symmetric about the real axis. We then introduce five different ways to select active nodes and describe how to construct BDF, GBDF, or Adams \ODEps{} from each type of active node set. In short, our proposed procedure for constructing a polynomial block method with imaginary nodes can be summarized in the following three steps:
	\begin{enumerate}
		\item Select an active node set from the possibilities listed in subsection \ref{subsub:active_index_sets}.
		\item Select a method family (Adams, BDF, or GBDF).
		\item Select a set of complex-valued nodes that are symmetric about the real axis.
	\end{enumerate}

	\subsection{Real-symmetric imaginary nodes}
	\label{sec:construction_choosing_nodes}
		
		Before discussing method construction we first select a node family. This choice does not fix the final node set, but rather helps inform future parameter choices and broadly characterizes the type of method we seek to construct. From here on, we restrict ourselves to imaginary node sets $\{z_j\}^q_{j=1}$ that are symmetric with respect to the real-axis so that 
			\begin{align*}
				\text{Re}(z_j) = 0 \quad \text{and} \quad \zeta \in \{ z_j\}^q_{j=1} \iff \zeta^* \in \{ z_j\}^q_{j=1}.
			\end{align*}

\begin{figure}[h]
				\centering 	
					\includegraphics[height=0.4\linewidth]{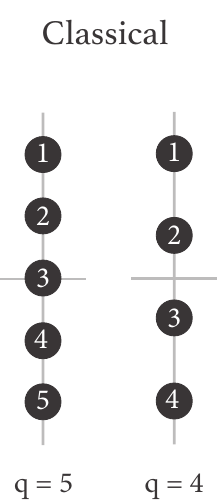}
					\hspace{5em}	
					\includegraphics[height=0.4\linewidth]{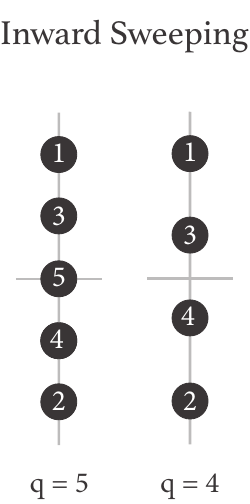}
					\hspace{5em}	
					\includegraphics[height=0.4\linewidth]{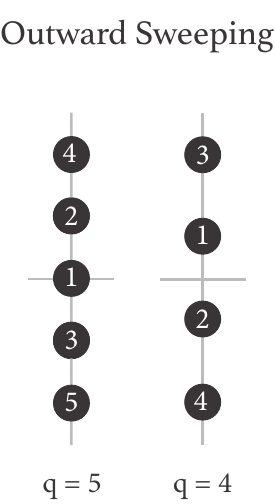}
				\caption{Real-symmetric imaginary node orderings shown for the node sets $\{z_j\}_{j=1}^4 = \{i,~i/3,~i/3,~i\}$ and $\{z_j\}_{j=1}^5 = \{-i,~-i/2,~0,~i/2,~ i\}$.}
				\label{fig:imaginary_node_ordering}
			\end{figure}

		Node ordering is insignificant for parallel block methods since their outputs must be computed independently. However, for serial methods the node ordering dictates which data values can be used to compute an output without making the method fully implicit. We therefore consider three separate orderings for real-symmetric imaginary nodes: 
			\begin{enumerate}
				\item {\bf Classical ordering}: nodes are ordered from top to bottom in the complex plane so that ${iz_1 < iz_2 < \ldots < iz_q}$.
				\item {\bf Outward sweeping}: nodes are ordered so that their magnitude monotonically increases:
				\begin{align*}
					|z_1| \le |z_2| \le \ldots \le |z_q| \quad \text{where} \quad |z_j| = |z_{j+1}| \implies iz_j > iz_{j+1}.
				\end{align*}

				\item {\bf Inward sweeping}: nodes are ordered so that their magnitude monotonically decreases: 
				\begin{align*}
					|z_1| \ge |z_2| \ge \ldots \ge |z_q| \quad \text{where} \quad |z_j| = |z_{j+1}| \implies iz_j > iz_{j+1}.
				\end{align*}
			\end{enumerate}
			For inward sweeping and outward sweeping orderings, the second condition prioritizes points in the upper-half plane over conjugate points in the lower-half plane. In Figure \ref{fig:imaginary_node_ordering} we show each ordering using two node sets of four and five imaginary equispaced nodes.
		
		In general, there are many additional families of nodes that will lead to interesting polynomial methods. For example, polynomial methods with real-valued nodes were discussed in \cite{buvoli2018polynomial}, and a more detailed presentation of these methods will be the topic of future work. Another example of a method with nodes that are neither purely imaginary or real-valued is the explicit integrator based on the roots of unity presented in \cite{buvoli2019constructing}.		
	
	\subsubsection{Conjugate outputs}

For polynomial integrators with complex nodes we can utilize the Swartz reflection principle to reduce the number of linear systems required per timestep when solving real-valued differential equations. This is achieved by creating polynomial block methods with {\em conjugate outputs} $y_j^{[n+1]}$ and $y_k^{[n+1]}$ that satisfy $$y_j^{[n+1]} = \left(y_k^{[n+1]}\right)^*.$$ To construct methods with conjugate outputs it is necessary to treat the upper and lower half planes symmetrically. The outputs $y_j^{[n+1]}$ and $y_k^{[n+1]}$ will always be conjugate if:
\begin{enumerate}
	\item The nodes $z_j$ and $z_k$ are conjugate such that $z_j = z_k^*$.
	\item The \ODEps{} for computing the outputs are conjugate so that $${p_1(z; b_1)^* = p_2(z^*; b_2) \quad  \forall z}.$$ 
\end{enumerate}
The BBDF and BAM methods from \cite{buvoli2019constructing} are two example methods with conjugate outputs. The primary advantage of a method with all conjugate outputs is that the computational cost is cut in half, since only half of the function evaluations and nonlinear solves must be computed at each timestep. This is especially important for serial methods where parallelism cannot be used to offset the additional linear solves required during each timestep. In light of this, all the methods we present in this work will possess conjugate outputs.
		
		\subsection{Choosing an active node index set}
		\label{sec:choosing_active_node_index_sets}
		
		A method's active node set, first defined in subsection \ref{subsection:stencil_diagrams}, determines the degree of implicitness of method and its architecture (parallel or serial). When describing the active nodes for a PBM's \ODEp{}, it is convenient to: (1) split the active node set into input and output nodes to improve readability, and (2) reference the node index rather than the node value. We will call the resulting two sets the {\em active input index set} and the {\em active output index set}. Their precise definitions are below:	
			\begin{definition}[Active input index (AII) and Active output index (AOI) sets] Consider a polynomial method with $q$ outputs. The $j$th AII set $I(j)$ and the $j$th AOI set $O(j)$ are sets containing integers ranging from 1 to $q$ such that:
				\begin{itemize}
					\item $k\in I(j)$ iff $p_j(\tau;b_j)$ is constructed using input data $y_k^{[n]}$ and/or $f_k^{[n]}$.
					\item $k\in O(j)$ iff $p_j(\tau;b_j)$ is constructed using output data $y_k^{[n+1]}$ and/or $f_k^{[n+1]}$.
				\end{itemize}
			\end{definition}
			Note that a polynomial method with $q$ nodes has $q$ total AOI and AOI sets. In Table \ref{tab:block_data_subsets}, we classify different types of PBMs based on the output index set $O(j)$.
			\begin{table}[h]
		\caption{Largest superset of the output index sets for various types of PBMs.}
		\begin{center}
			\renewcommand{\arraystretch}{2}
			\begin{tabular}{l|p{9em}|p{9em}|p{9em}}
				& {\em Explicit} & {\em Diagonally Implicit} & {\em Fully Implicit} \\ \hline
				{\em Parallel} & 
				$O(j) = \emptyset$ &
				$O(j) \subseteq \{j\}$ & 
				no such methods \\	\hline	
				{\em Serial} & 
				$O(j) \subseteq  \{1, \ldots, j-1\}$ &
				$O(j) \subseteq \{1, \ldots, j\}$ &
				$O(j) \subseteq \{1, \ldots, q\}$ \\ 
			\end{tabular}	
		\end{center}	
		\label{tab:block_data_subsets}
		\end{table}

			The first step in our proposed procedure for constructing polynomial block methods is to choose an active node set. This allows us to fix a method's architecture and its degree of implicitness independently of its \ODEps{}. In the following subsections we present five ways to choose the active input index (AII) set and active output index (AOI) set for polynomial block methods with real-symmetric imaginary nodes. Our formulas can be used to construct either explicit or diagonally-implicit integrators.
			
			\subsubsection{A naming convention for active index sets}
			
			We introduce a naming convention for each pair of active index sets that describe the underlying method properties. Since active index sets determine the architecture of a method, each name starts with the word {\em parallel} or {\em serial}. Next, the names reflect the cardinality of the AII and AOI sets. The theoretical maximum order of a polynomial method is proportional to the cardinality of the active index sets. Larger cardinalities mean that  more data values can be used to construct polynomials, leading to increased order of accuracy. Active index sets that achieve a certain property (e.g. parallelism) using the largest amount of data will have the word {\em maximal} in their name. Finally, the cardinalities of the active index sets can remain fixed or may vary for each of the method's outputs. All names will contain the words {\em fixed cardinality} or {\em variable cardinality} to reflect these two possibilities.
			
			\subsubsection{Notation for active index sets}

			The AOI set $O(j)$ must be defined differently for explicit, diagonally-implicit, or fully-implicit integrators. We can write a compact formula for $O(j)$ in terms of another set we call $B(j)$. The AOI set is then defined as
			\begin{align}
				O(j) = 
				\begin{cases}
					B(j) & \text{for explicit methods,} \\
					B(j) \cup \left\{ j \right\}	 & \text{for diagonally-implicit methods,} \\
					B(j) \cup \left\{ 1, \ldots, q \right\} & \text{for fully-implicit methods,}
				\end{cases} 
				\label{eq:AOI_set_b}
			\end{align}
			where the formulas for the set $B(j)$ are be contained in  subsection \ref{subsub:active_index_sets}. Finally, for methods with inwards or outwards node orderings, the AII set and the set $B(j)$ are expressed in terms of the function $\chi_{\text{in}}(j)$ and $\chi_{\text{out}}(j)$ that are defined as:
			\begin{align}
				\begin{aligned}
				&\text{inwards ordering:}&
				&\left\{\begin{aligned}
					\chi_{\text{in}}(j) &= 
					\begin{cases}
						j & j \text{ odd} \\
						j - 1 & j \text{ even}
					\end{cases}	\\[1em]
					\chi_{\text{out}}(j) &= 
					\begin{cases}
						1 & j \text{ odd} \\
						2 & j \text{ even}
					\end{cases}
				\end{aligned}, \right.
				\\[1 em]
				&\text{outwards ordering:}&
				&\left\{ \begin{aligned}
				\chi_{\text{in}}(j) &= 
					\begin{cases}
						\max(1, j-1) & q \equiv j \mod 2 \\
						j & \text{otherwise}
					\end{cases}	\\[1em]
					\chi_{\text{out}}(j) &= 
					\begin{cases}
						1 & q \equiv j \mod 2 \\
						2 & \text{otherwise}
					\end{cases}
				\end{aligned} \right..
				\end{aligned}
				\label{eq:imaginary_shift_functions}
			\end{align}

	\subsubsection{Node diagrams}

\begin{figure}[h]
			\begin{center}
				\hfill
				\begin{minipage}[t]{0.45\textwidth}
					\centering
					\includegraphics[width=0.3\linewidth]{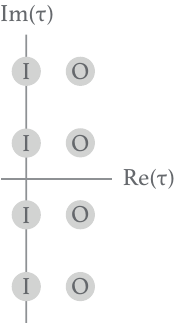} \\[0.5em]
					\begin{small}
						{\bf (a)} Four imaginary equispaced nodes.		
					\end{small}	
				\end{minipage}
				\hfill
				\begin{minipage}[t]{0.45\textwidth}
					\centering
					\includegraphics[width=0.3\linewidth]{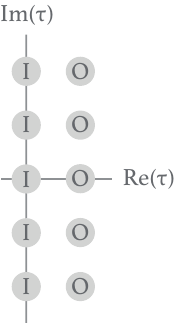} \\[0.5em]
					\begin{small}
					{\bf (b)} Five imaginary equispaced nodes.	
					\end{small}
				\end{minipage}
				\hfill
			\end{center}
	
		    \caption{Empty active node stencils for the \ODEps{} of a block method with: ({\bf a}) nodes $\{z_j\}_{j=1}^4 = \{i,~i/3,~i/3,~i\}$, and ({\bf b}) nodes $\{z_j\}_{j=1}^5 = \{-i,~-i/2,~0,~i/2,~ i\}$. For clarity, the inactive nodes in both diagrams have been enlarged, and inputs are labeled with the letter I while outputs are labeled with the letter O.}
			\label{fig:empty_node_diagram} \label{fig:family_it_prop}
		\end{figure}

	In addition to formula, we also provide active node diagrams to visualize the active index sets. To appreciate the simple geometric construction underlying each active index set, the formula for each AII and AOI set should be read in tandem with their corresponding node diagrams that are shown in Figures \ref{fig:active_node_diagrams_imaginary_implicit_even} and \ref{fig:active_node_diagrams_imaginary_implicit_odd}. All the active node diagrams were created using the imaginary equispaced nodes 
	\begin{align*}
		\{z_j\}_{j=1}^4 = \{i,~i/3,~i/3,~i\}	 \quad \text{or} \quad \{z_j\}_{j=1}^5 = \{-i,~-i/2,~0,~i/2,~ i\}.
	\end{align*}
	 We present visualizations for both sets of nodes since methods with an odd number of nodes have different characteristics due to the real-valued  node at $\tau = 0$.  For clarity we show an empty node stencil for each node sets in Figure \ref{fig:empty_node_diagram}. 
	
	\subsubsection{Five types of active index sets}
	\label{subsub:active_index_sets}
	
\begin{figure}[h]
			\centering
			\includegraphics[width=0.82\linewidth]{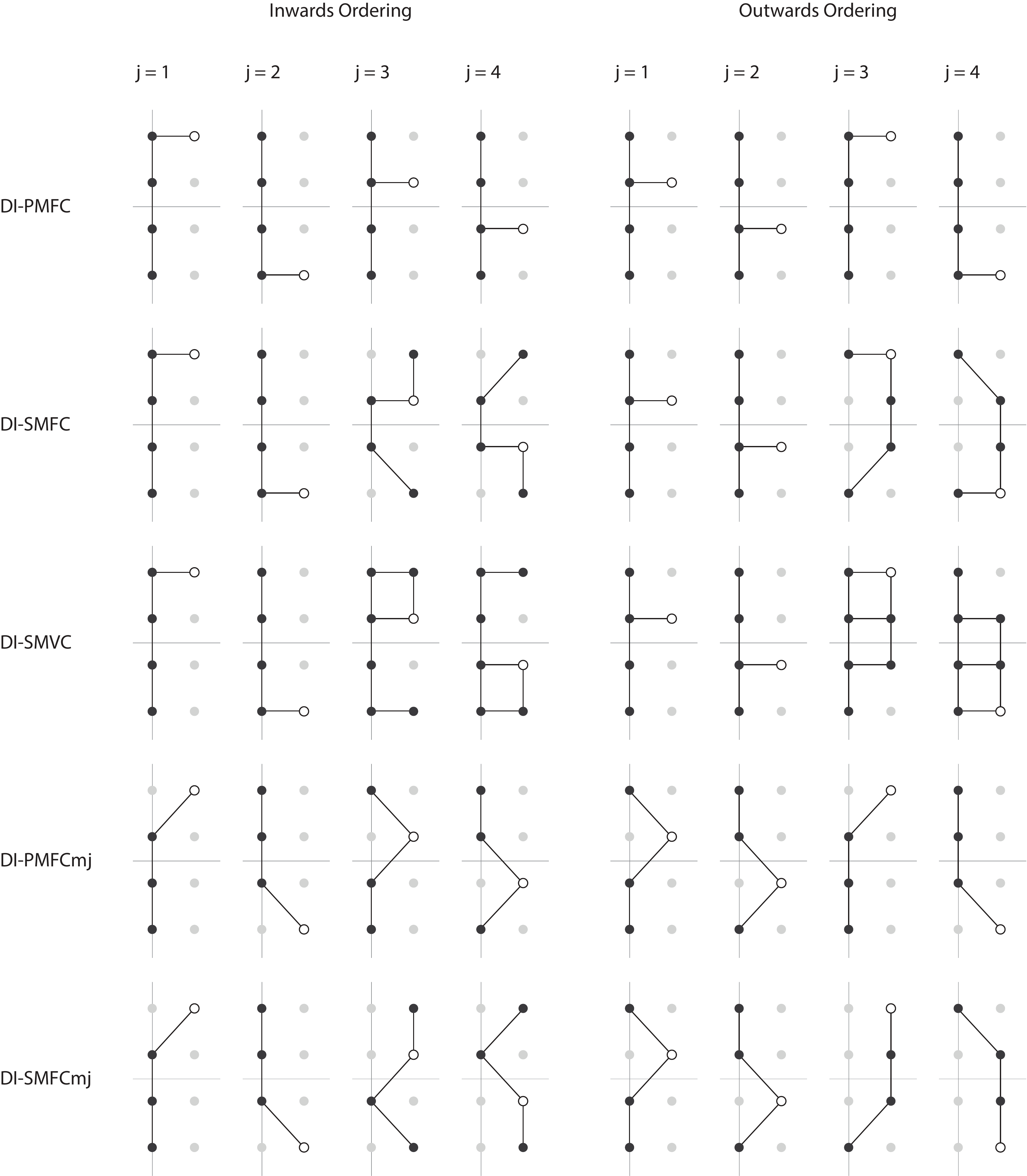}	
			\caption{Node stencils for diagonally implicit block methods with four imaginary equispaced points. Light gray circles denote inactive nodes, black circles denote active nodes, and a white circle with a black border denotes the active output node.}
			\label{fig:active_node_diagrams_imaginary_implicit_even}
		\end{figure}

\begin{figure}[h]
			\centering
			\includegraphics[width=1\linewidth]{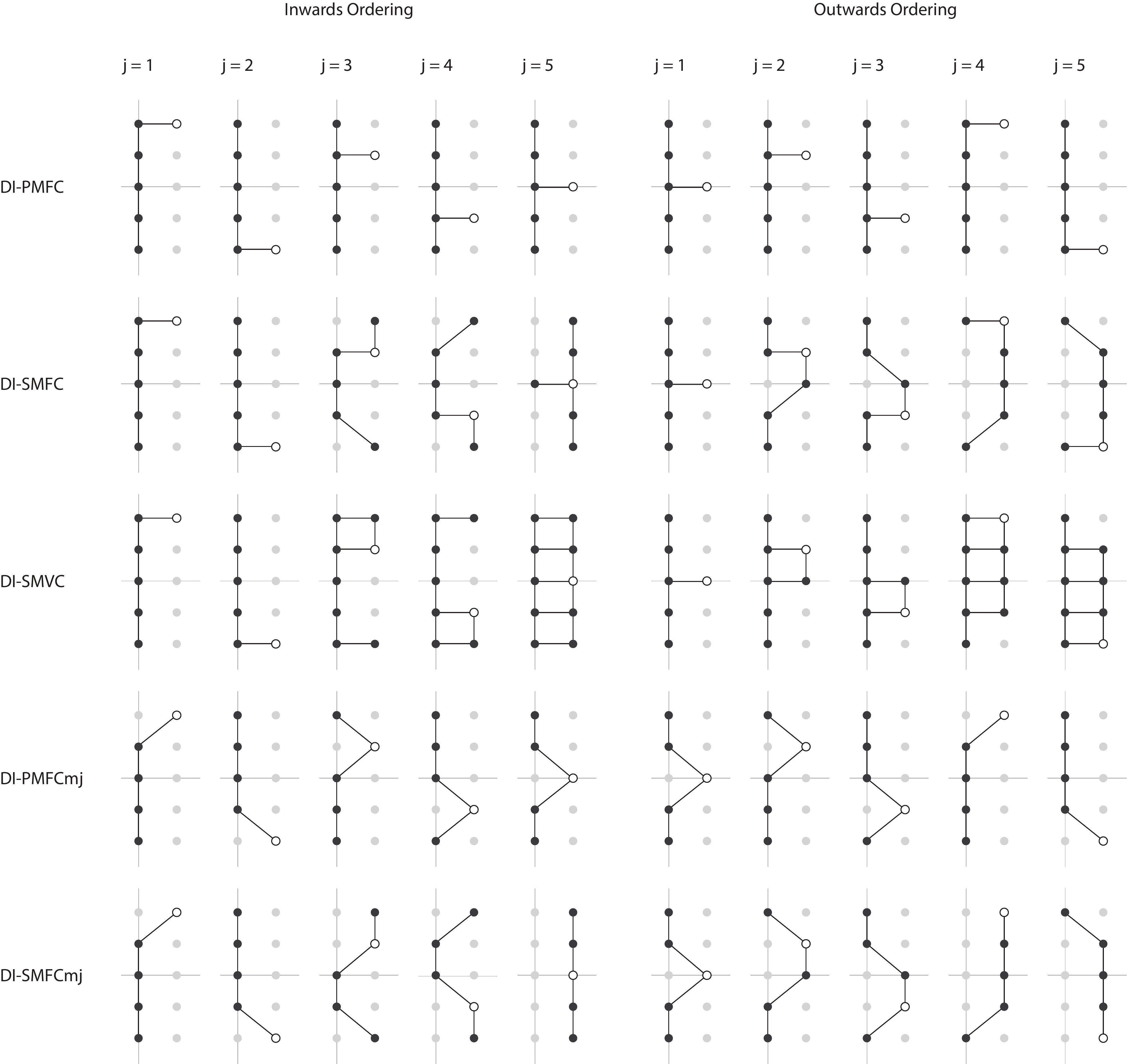}	
			\caption{Node stencils for diagonally implicit block methods with five imaginary equispaced points. Light gray circles denote inactive nodes, black circles denote active nodes, and a white circle with a black border denotes the active output node.}
			\label{fig:active_node_diagrams_imaginary_implicit_odd}
		\end{figure}

	 We now present five choices for selecting the AII set $I(j)$ and the AOI set $B(j)$. The corresponding node diagrams for diagonally implicit integrators are shown in Figures \ref{fig:active_node_diagrams_imaginary_implicit_even} and \ref{fig:active_node_diagrams_imaginary_implicit_odd} and the corresponding formulas are:				
		\begin{enumerate}%
			\item {\em Parallel maximal-fixed-cardinality} ({\bf PMFC}): each output is computed using data at all input nodes. The set $B$ must be empty to retain parallelism. The formulas are:
			\begin{align}
				I(j) = \left\{ 1, \ldots,q\right\} 
				\quad \text{and} \quad
				B(j) = \left\{ \right\}.
			\end{align} 
			
			\item {\em Serial maximal-variable-cardinality} ({\bf SMVC}): each output is computed using data from all input nodes and all previously computed outputs that allow for conjugate polynomials. The cardinality of set $B(j)$ grows as we add new information. The formulas are
			\begin{align}
				I(j) = \left\{ 1, \ldots, q\right\} 
				\quad \text{and} \quad
				B(j) = \left\{ 1, \ldots, \chi_{\text{out}}(j) \right\}.
			\end{align}

			\item {\em Serial Maximal-fixed-cardinality} ({\bf SMFC}): each output is computed using data from all previously computed outputs that allow for conjugate polynomials and some of the inputs. To keep the cardinality of $I(j) \cup O(j)$ fixed for all $j$, we drop inputs in favor of more recently computed outputs. The corresponding formulas are  
			\begin{align}
				I(j) = \left\{ \chi_{\text{in}}(j),~ \ldots~,~ q\right\}
				\quad \text{and} \quad
				B(j) = \left\{ 1, \ldots, \chi_{\text{out}}(j) \right\}.
			\end{align}

			\item {\em Parallel Maximal-fixed-cardinality minus $j$} ({\bf PMFCmj}): similar to PMFC, except now the $j$th output is computed using data from all input nodes, excluding the input with index $j$. The formulas are
			\begin{align}
				I(j) = \left\{ 1, \ldots, q\right\} \setminus \{ j \}
				\quad \text{and} \quad
				B(j) = \left\{ \right\}.
			\end{align}
			
			\item {\em Serial maximal-fixed-cardinality minus j} ({\bf SMFCmj}): similar to SMFC, except the $j$th output cannot use data from the $j$th node. The cardinality of $I(j) \cup B(j)$ remains fixed across all $j$. The formulas are 
			\begin{align}
				I(j) = \left\{ \chi_{\text{in}}(j), \ldots, q\right\} \setminus \{ j \}
				\quad \text{and} \quad
				B(j) = \left\{ 1, \ldots, \chi_{\text{out}}(j) \right\}.
			\end{align}
		\end{enumerate}

		\subsection{Choosing \ODEps{}}
		\label{sec:choosing_odepolynomials}

				The next step in our procedure for constructing polynomial block methods is to choose the \ODEps{} for computing each output. Here we describe how to construct polynomials for Adams, BDF, and GBDF methods using the AII set $I(j)$, the AOI set $O(j)$, and the set $B(j)$ from the previous subsection. Depending on the desired type of method the \ODEp{} can be chosen as follows:

		\begin{enumerate}%
			\item For an {\bf Adams} PBM, select Adams \ODEps{}
				\begin{align*} 
					p_j(\tau; b) = L^{[j]}_y(\tau) + \int_b^\tau L^{[j]}_F(\tau)
				\end{align*}
				 where the Lagrange polynomial $L_F^{[j]}(\tau)$ interpolates the input derivatives $f^{[n]}_k$ for $k\in I(j)$ and the output derivatives $f^{[n+1]}_k$ for $k \in O(j)$. The Lagrange polynomial $L_y^{[j]}(\tau)$ should be constructed differently depending on the particular choice of expansion points $\{ b_j \}_{j=1}^q$. In Subsection \ref{subsection:constructing_adams_endpoint} we present several possibilities.

			\vspace{0.5em}
			\item For a {\bf BDF} PBM, the BDF \ODEps{} $p_j(\tau; b) = H^{[j]}_y(\tau)$ are constructed differently for implicit and explicit methods. For an implicit method, $H^{[j]}(\tau)$ must interpolate the input values $y^{[n]}_k$, for $k\in I(j)$, and the output values $y^{[n+1]}_k$, for $k\in B(j)$, and its derivative $\frac{d}{d\tau}H^{[j]}(\tau)$ must interpolate the output derivative $f^{[n+1]}_j$. 

			For an explicit method, we cannot use the output derivative $f_j^{[n+1]}$. Instead we construct an \ivs{} that contains the interpolated derivatives
				\begin{align*}
					\tilde{f}_j = \dot{p}_j(z_j + \alpha; b) && j = 1, \ldots, q
				\end{align*}
				where $\dot{p}_j(\tau; b)$ is a Lagrange interpolating polynomial that interpolates the input derivatives $f^{[n]}_k$ for $k \in I(j)$ and the output derivatives $f^{[n+1]}_k \text{ for } k \in O(j)$. Then, we construct $H^{[j]}(\tau)$ in the same way as the implicit case, except its derivative now interpolates $\tilde{f}_j$ instead of $f_j^{[n+1]}$.
			
			\vspace{0.5em}
			\item For a {\bf GBDF} PBM, the GBDF \ODEps{} $p_j(\tau; b) = H^{[j]}_y(\tau)$ are also constructed differently for implicit and explicit methods. For implicit methods, $H^{[j]}(\tau)$ interpolates the input values $y^{[n]}_k $ for $k\in I(j)$ and its derivative $\frac{d}{d\tau}H^{[j]}(\tau)$ must interpolate the output derivatives $f_j^{[n+1]}$ for ${k\in O(j)}$.  
		
				For explicit methods, we construct the interpolated derivatives $\tilde{f}_j$, $j = 1, \ldots, q$ in the same way as for an explicit BDF method. Then we construct $H^{[j]}(\tau)$ so that it interpolates the input values $y^{[n]}_k $ for $k\in I(j)$ and its derivative $\frac{d}{d\tau}H^{[j]}(\tau)$ interpolates the output derivatives $f_j^{[n+1]}$ for $k\in B(j)$ and the interpolated derivative $\tilde{f}_j$.

		\end{enumerate}

		\subsubsection{Endpoints for Adams \ODEps{}}
		\label{subsection:constructing_adams_endpoint}

		 Each output of an Adams polynomial block method is computed by evaluating Adams polynomials so that
		\begin{align*}
			y_j^{[n+1]} = p_j(z_j + \alpha; b_j) = L_y^{[j]}(b_j) + \int_{b_j}^{z_j + \alpha} L^{[j]}_F(s) ds.
		\end{align*}
		Here we present three ways to select the endpoints $b_j$ and Lagrange polynomials $L_y^{[j]}(\tau)$. Though both parameters can be chosen independently, it is generally preferable to consider them simultaneously if one wants to obtain the highest order of accuracy using the most compact representation for $L^{[j]}_y(\tau)$. For brevity we only consider endpoint sets where each endpoint $b_j$ is equal to one of the nodes $z_j$. This restriction eliminates any complexity in selecting $L_y^{[j]}(\tau)$, since it is now possible to choose $L_y^{[j]}(\tau)$ to be a constant that is equal to the solution value at the endpoint.

\begin{figure}[h!]
				\centering
					{\small {(a) Endpoints: q = 5}	 \\[0.5em]}
				\includegraphics[width=.75\linewidth]{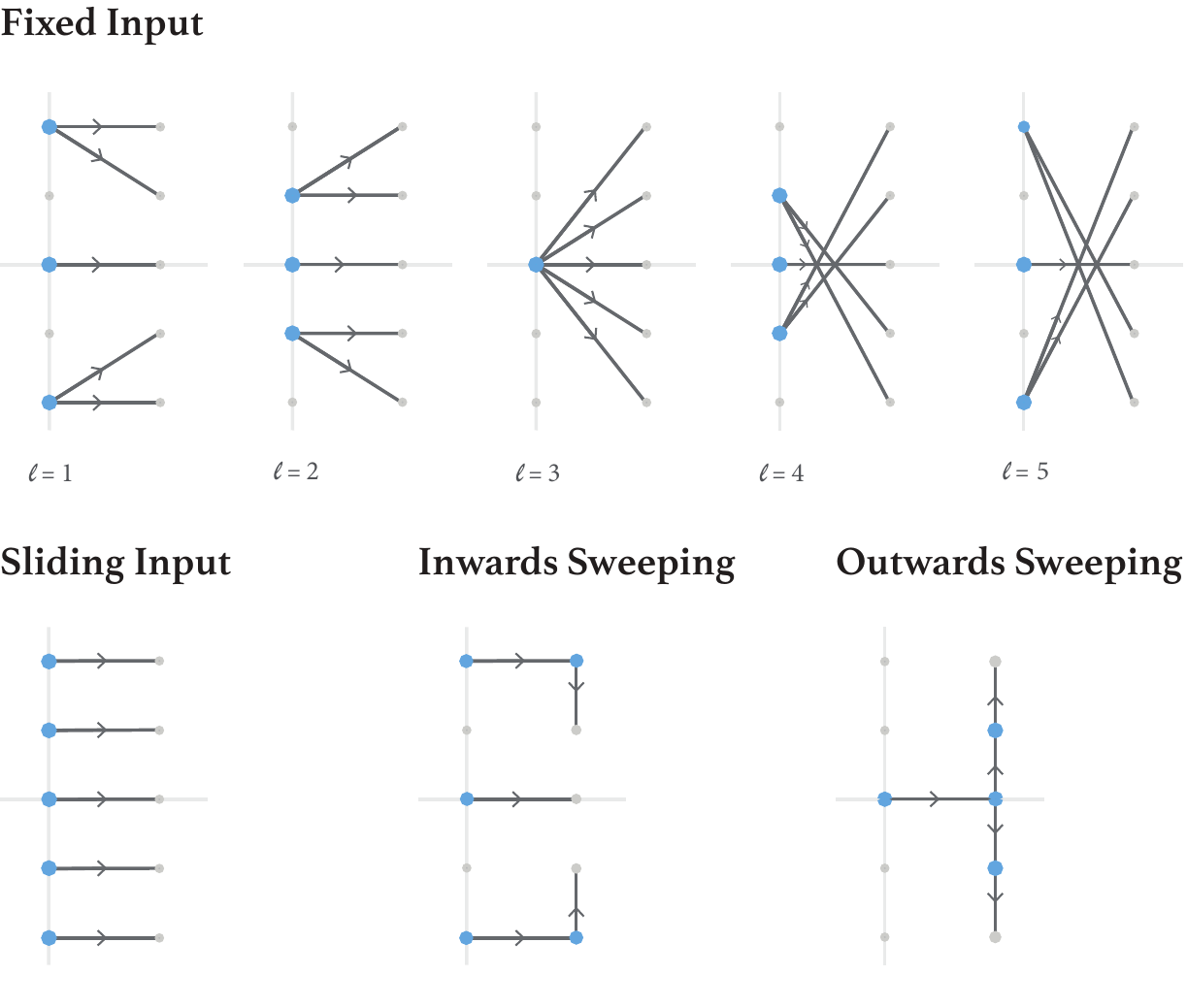}

					{\small \vspace{1em} { (b) Endpoints: q = 4}	 \\[1em]}
				\includegraphics[width=.6\linewidth]{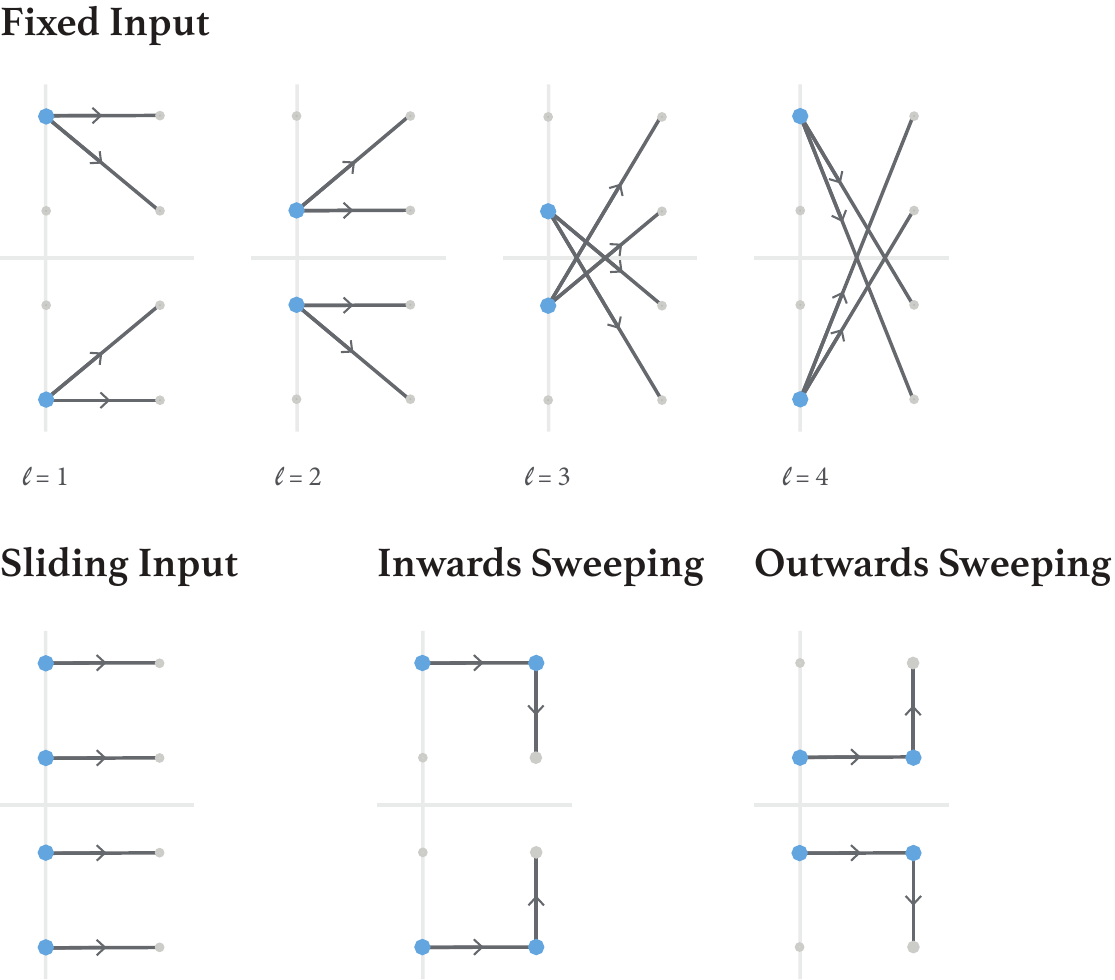}	
				\caption{Endpoint diagrams for Adams PBMs with five and four equispaced Nodes.}
				\label{fig:endpoints}
		\end{figure}

		For improved readability, the formulas for the endpoints are each written in terms of different node orderings. Using the mappings provided in \cite{buvoli2018polynomial, buvoli2020PIPACK} these formulas can be easily re-expressed using any of the three orderings (inwards, outwards or classical). 	In addition to the endpoint formula, we also show endpoint diagrams   in Figure \ref{fig:endpoints} to highlight the simple geometric ideas behind each parameter choice. The three proposed endpoints are:

		\begin{enumerate}
			\item {\em Fixed input} ({\bf FI}): This expansion point can be used for constructing both serial and parallel block methods. For fixed input, the expansion points $b_j$ is either equal to a node $z_\ell$ or its conjugate $z_\ell^*$. For a method with $q$ outputs there are $q$ possibilities for choosing $\ell$. For nodes in classical ordering the formulas are 
			\begin{align*}
				b_j = z_{\text{ind}(j)} \quad \text{and} \quad L_y^{[j]} =  y^{[n]}_{\text{ind}(j)},	\quad j = 1, \ldots, q,
			\end{align*}
			where the function ind($j$) is defined as
			\begin{align*}
				& \text{$q$ even} & \text{ind}(j) &= 
				\begin{cases}
					\ell & j \le \tfrac{q}{2}, \\
					q - \ell + 1 & j > \tfrac{q}{2},
				\end{cases}
				\\[0.5em]
				& \text{$q$ odd} & \text{ind}(j) &= 
				\begin{cases}
					\ell & j < \lceil \tfrac{q}{2} \rceil \\
					\lceil \tfrac{q}{2} \rceil & j = \lceil \tfrac{q}{2} \rceil \\ %
					{q - \ell + 1} & j > \lceil \tfrac{q}{2} \rceil 
				\end{cases}	.
			\end{align*} 

			\item {\em Variable input} ({\bf VI}): This choice of endpoint can be used to construct methods with serial or parallel architectures. The expansion point $b_j$ is equal to the temporal node of the $j$th input, so that
			\begin{align*}
				b_j = z_j	\quad \text{ and } \quad L_y^{[j]} = y_j^{[n]}.
			\end{align*}
 			This formula is valid for all node orderings.

			\item {\em Sweeping}	({\bf S}): This expansion point set leads to a serial method and is defined differently for inwards and outwards node orderings. For inwards sweeping, the expansion points start at the temporal nodes of the inputs that are furthest from the real line. The expansion points then sweep inwards along the temporal nodes of previously computed outputs. If there is a real-valued node, which occurs when $q$ is odd, then the final output is integrated along the real line. The formulas are:
			\begin{align*}
				& \text{q even:}\hspace{2em} & 
				b_j &= 
				\begin{cases}
					z_j & j \le 2, \\
					z_{j-2}+\alpha & j > 2,
				\end{cases} & 
				L^{[j]}_y = & 
				\left\{ \hspace{-.3em} 
				\renewcommand\arraystretch{1.5}
				\begin{array}{ll}
					 y_j^{[n]}  & j \le 2, \\
					 y_{j-2}^{[n+1]}  & j > 2.
				\end{array}
				\right. &  				
				\\[0.5em]
				& \text{q odd:} & 
				b_j &= 
				\begin{cases}
					z_j & j \le 2, \\
					z_{j-2}+\alpha & 2 \le j < q, \\
					z_{q} & j = q,
				\end{cases} & 
				L^{[j]}_y = &				
				\left\{ \hspace{-.3em}
				\renewcommand\arraystretch{1.5}
				\begin{array}{ll}
					 y_j^{[n]} 		& j \le 2, \\
					 y_{j-2}^{[n+1]} 	& 2 \le j < q, \\
					 y_{q}^{[n]} 		& j = q.
				\end{array} 
				\right.
			\end{align*} 
			For outwards sweeping nodes, the expansion points start at the temporal nodes of the input that are nearest to the real line, before sweeping outwards along the temporal nodes of previously computed outputs. The formulae are:
			\begin{align*}
				& \text{q even:}\hspace{2em} & 
				b_j &= 
				\begin{cases}
					z_j & j \le 2, \\
					z_{j-2} + \alpha & j > 2,
				\end{cases} & 
				L^{[j]}_y = & 
				\left\{ \hspace{-.3em}
				\renewcommand\arraystretch{1.5}
				\begin{array}{ll}
					 y_j^{[n]}  & j \le 2, \\
					 y_{j-2}^{[n+1]}  & j > 2,
				\end{array}
				\right. &  				
				\\[0.5em]
				& \text{q odd:} & 
				b_j &= 
				\begin{cases}
					z_1 & j = 1, \\
					z_{1}+\alpha & 1 < j \le 3, \\
					z_{j-2} & j > 3,
				\end{cases} & 
				L^{[j]}_y = &				
				\left\{ \hspace{-.3em}
				\renewcommand\arraystretch{1.5}
				\begin{array}{ll}
					 y_1^{[n]} 		& j = 1, \\
					 y_{1}^{[n+1]} 	& 1 < j \le 3, \\
					 y_{j-2}^{[n]}  	& j > 3.
				\end{array} 
				\right. 
			\end{align*}

	\end{enumerate}

	\subsection{Relation to previously introduced methods}	
		
		Using the new construction strategies, the BAM and BBDF methods from \cite{buvoli2019constructing} can be obtained by pairing imaginary equispaced points with the constructions show in Table \ref{tab:BAM_BBDF_constructions}. Since BAM and BBDF are parallel methods the node ordering is irrelevant.
		\begin{table}[h!]
			\centering
			\renewcommand\arraystretch{1.25}
			\begin{tabular}{r|lll}
			 & AIS & ODEP & EP \\ \hline
			{\bf BAM} & PMFC & Adams & Variable Input \\
			{\bf BBDF} & PMFC & BDF & n.a.	
			\end{tabular}
			\caption{Constructions for BAM and BBDF. The abbreviations AIS, ODEP, and EP are used to abbreviate active index set, \ODEp{}, and expansion points.}
			\label{tab:BAM_BBDF_constructions}
		\end{table}

		\section{Linear stability}
		\label{sec:linear_stability} 
			
		Now that we have introduced a wide range of strategies for constructing PBMs, we investigate linear stability. Like all classical integrators, linear stability for PBMs is determined using the Dahlquist test problem $y' = \lambda y$. When applied to this problem, a PBM reduces to the iteration $\mathbf{y}^{[n+1]} = \mathbf{M}(\zeta,\alpha) \mathbf{y}^{[n]}$ where $\zeta = r \lambda$, $\alpha$ is the extrapolation factor,  and $\mathbf{M}(\zeta,\alpha)$ is a $q \times q$ matrix. The stability region of the method is the subset of the complex $\zeta$-plane
		\begin{align*}
			S(\alpha) &= \left\{ \zeta ~:~ \sup_{n\in \mathbb{N}} \| \mathbf{M}(\zeta/\alpha,\alpha)^n \| < \infty \right\},
		\end{align*}
		where the factor of $\zeta/\alpha$ scales the stability regions relative to the stepsize $h$ rather than the node radius $r$.
		
		Since we are interested in methods for solving stiff initial value problems, we focus solely on diagonally implicit PBMs and their stability angles. Recall that a method is $A(\theta)$ stable if its stability region $S$ satisfies 
		\begin{align*}
		S \supset \left\{z :  |\arg(-z)| < \theta,~ z \ne 0\right\},
		\end{align*}
		where $A(90\degree)$ implies that the entire left-half $z$-plane lies inside the stability region.

		The construction strategies presented in Section \ref{sec:constructing_pbms} can be used with any imaginary node set but this leaves too much generality to present in this initial work. We therefore only consider nodes sets $\left\{ z_j \right\}^{q}_{j=1}$ with $2 \le q \le 8$ and  
		\begin{align*}
				z_j &= -i + 2ij/q  && \text{(imaginary equispaced),} \\
				z_j &= i \cos(\pi(j-1)/(q-1)) && \text{(imaginary Chebyshev).}	
		\end{align*}
		
		One of the primary aims of this work is to look for polynomial methods with stability regions that encompass the imaginary axis. More broadly we are also interested in categorizing $A(\theta)$ stability as a function of the extrapolation parameter $\alpha$. Ideally, we seek methods that are $A(90 \degree)$ stable for large $\alpha$ values. Recall that when the extrapolation factor is small, the ratio between the stepsize and the interpolation interval is also small. Therefore, the interpolation nodes will extend further into the complex plane, which requires either extra analyticity in the solution, or a reduction in the stepsize to prevent instabilities. 
			
		Even after restricting ourselves to diagonally implicit PBMs with imaginary equispaced or Chebyshev nodes, we still cannot present the stability properties for all the types of methods presented in Section \ref{sec:constructing_pbms}. Instead, in Figure \ref{fig:DIPBM_AofTheta_vs_alpha} we show the $A(\theta)$ properties for seven families of methods that possess A($90\degree$) stability for some range of $\alpha$ between zero and one. Amongst the example methods are BDF, BBDF, and Adams type methods with both parallel and serial architectures. For all BDF and BBDF methods, the stability angles increase monotonically as $\alpha$ decreases, while the curves for Adams methods are less regular with certain $\alpha$ ranges allowing for A($90\degree$) stability. We also found many new Adams methods that do not posses A($90\degree$) stability for any range of $\alpha$, but can still be used to solve dissipative problems. In Table \ref{tab:Example_Adams_methods_equi} we present a non-exhaustive list of such methods. In Figure \ref{fig:Adams_stability_angles} we also show the stability angles for two families of Adams methods that we will use later in our numerical experiments.
		
		By searching through the proposed families of PBMs and tuning the extrapolation parameter $\alpha$ we are able construct of a wide range of methods with desirable stability properties. As a general rule, stability improves as the extrapolation factor decreases. This is especially true for high-order PBMs where one must trade analyticity requirements for A($90\degree$) stability. To allow for a more complete analysis of the methods, we postpone further discussion of the stability results until after the numerical experiments.

\newenvironment{GBDFDescTable}[4] %
{
	\begin{scriptsize}
		\begin{tabular}{l}
			{\bf AII/AOI}   \\ 
			#1 \\[0.25em]
			{\bf ODEP}  \\ 
			#2 \\[0.25em]
			{\bf Nodes} \\ 
			#3, #4 \\[2em]
		\end{tabular}
	\end{scriptsize}
}

\newenvironment{AMDescTable}[5] %
{
	\begin{scriptsize}
		\begin{tabular}{l}
			{\bf AII}   \\
			#1 \\[0.25em]
			{\bf ODEP}  \\
			 #2 \\[0.25em]
			{\bf EndPoints}  \\
			 #3 \\[0.25em]
			{\bf Nodes} \\
			 #4, #5
		\end{tabular}
	\end{scriptsize}
}

\newenvironment{MethodNameTable}[1] %
{
	\begin{scriptsize}
		\begin{tabular}{l}
			#1
		\end{tabular}
	\end{scriptsize}
}

\newcommand{\lineLegend}[2]{
	{\tiny \textcolor{#1}{\hdashrule[0.2ex]{1.5em}{2pt}{}} #2 \hspace{-1.5em}}
}

\newcolumntype{f}{>{\hsize=1.2\hsize\linewidth=\hsize}X} %
\newcolumntype{m}{>{\hsize=.60\hsize}X} %

		\begin{figure}[hp!]
		
			\centering

			\hfill
			\begin{tabular}{lllllll}
				\lineLegend{plot_red}{$q=2$} &
				\lineLegend{plot_orange}{$q=3$} &
				\lineLegend{plot_yellow}{$q=4$} &
				\lineLegend{plot_green}{$q=5$} &
				\lineLegend{plot_blue}{$q=6$} &
				\lineLegend{plot_violet}{$q=7$} &
				\lineLegend{plot_grey}{$q=8$} \hspace{1.5em} 
			\end{tabular}

			\begin{tabularx}{\linewidth}{mff}
				\GBDFDescTable{SMVC}{BDF}{Equi}{Inwards}	 & 
					\includegraphics[valign=m,width=1\linewidth]{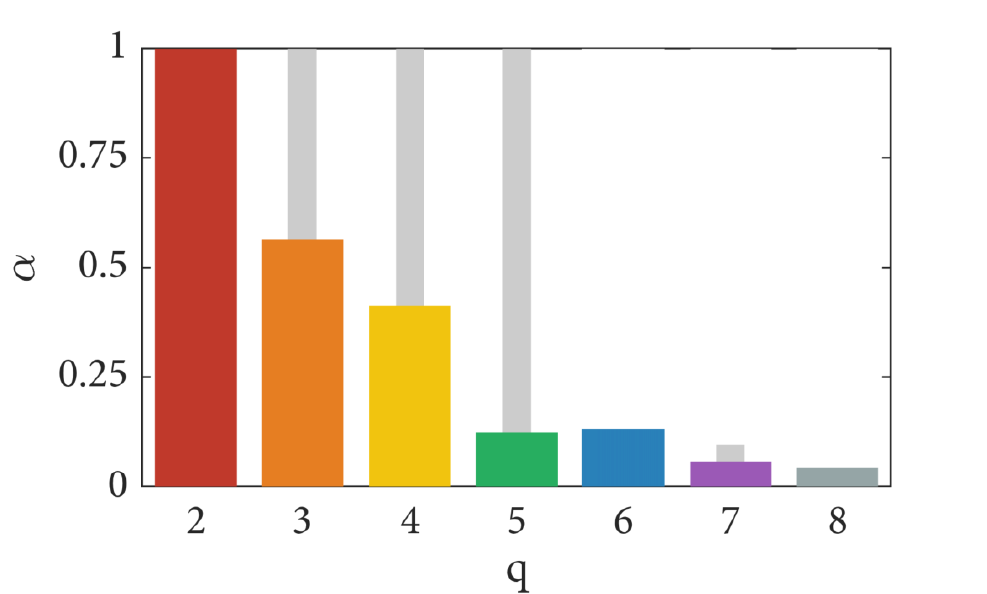}	 &
					\includegraphics[valign=m,width=1\linewidth]{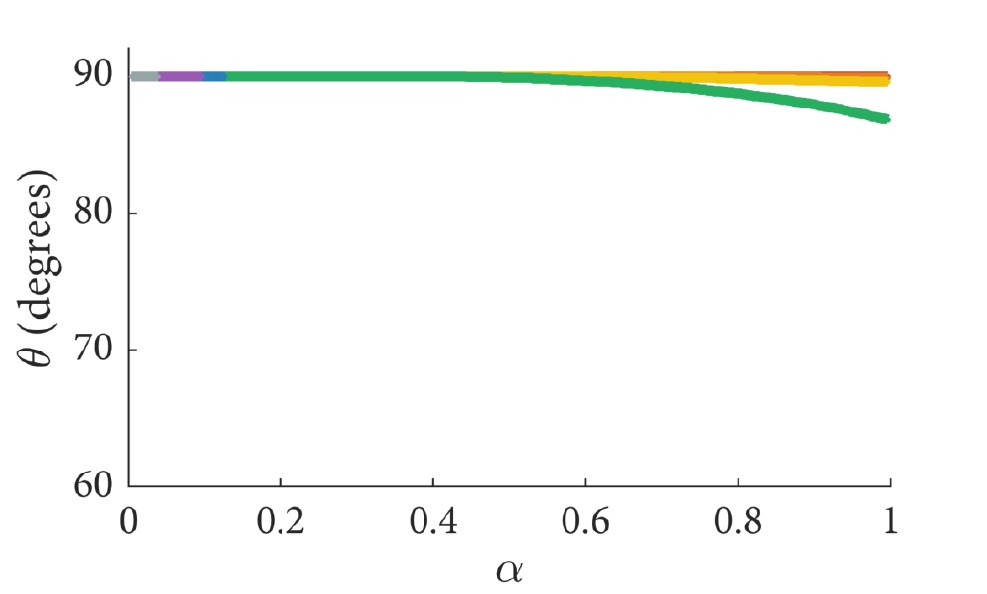} \\[1em]
				\GBDFDescTable{SMFC}{BDF}{Equi}{Inwards} &
					\includegraphics[valign=m,width=1\linewidth]{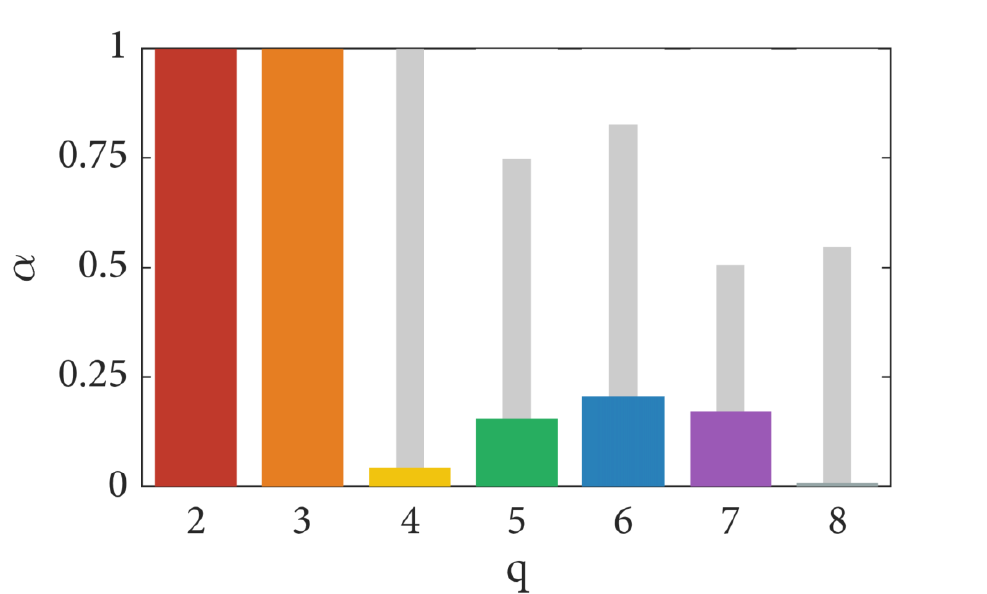}	 &
					\includegraphics[valign=m, width=1\linewidth]{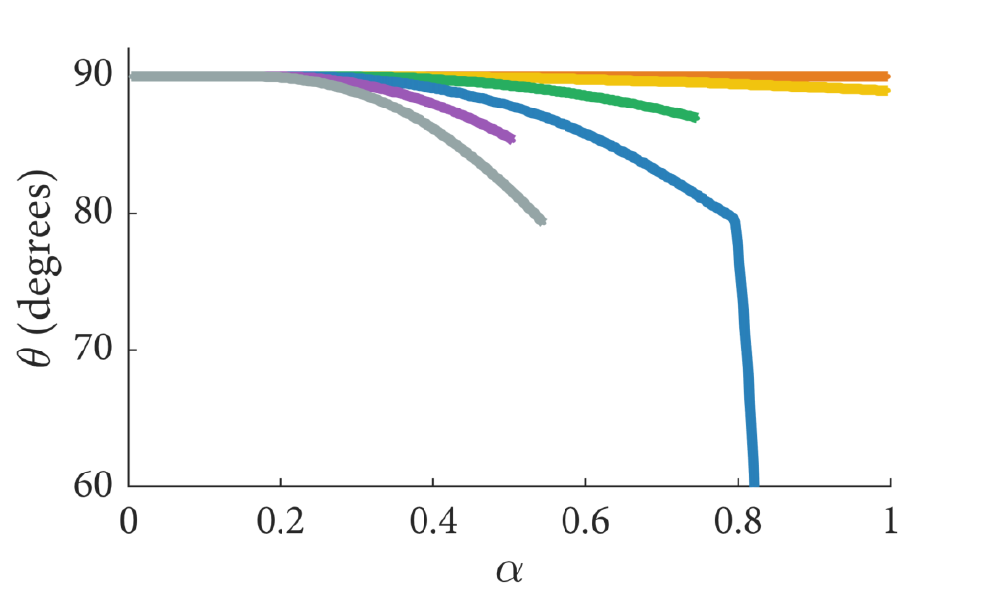} \\[1em]
				\GBDFDescTable{SMVC}{GBDF}{Equi}{Inwards}	 & 
					\includegraphics[valign=m,width=1\linewidth]{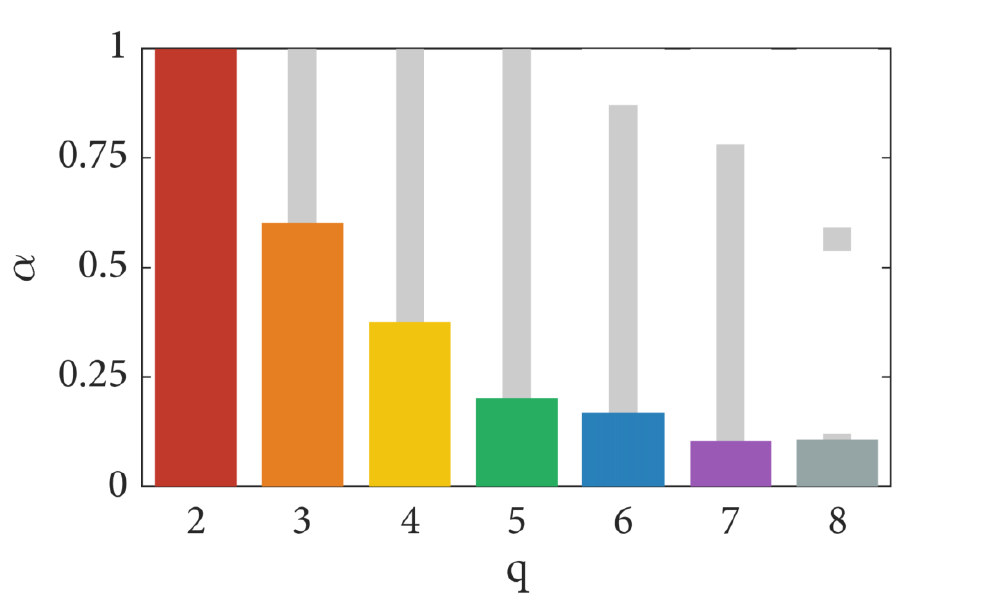}	 &
					\includegraphics[valign=m,width=1\linewidth]{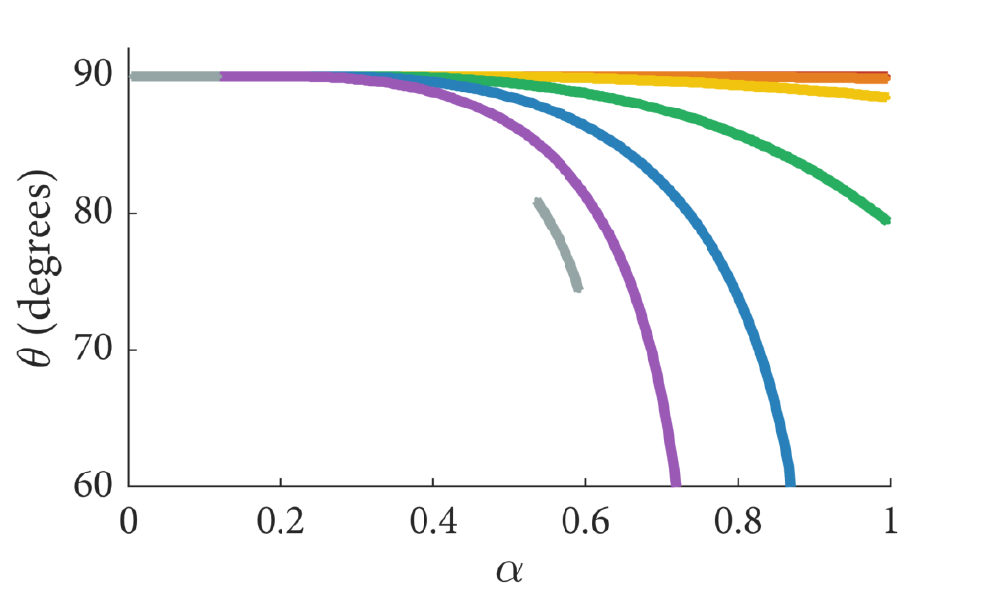} \\[1em]
				\GBDFDescTable{SMFC}{GBDF}{Equi}{Inwards} &
					\includegraphics[valign=m,width=1\linewidth]{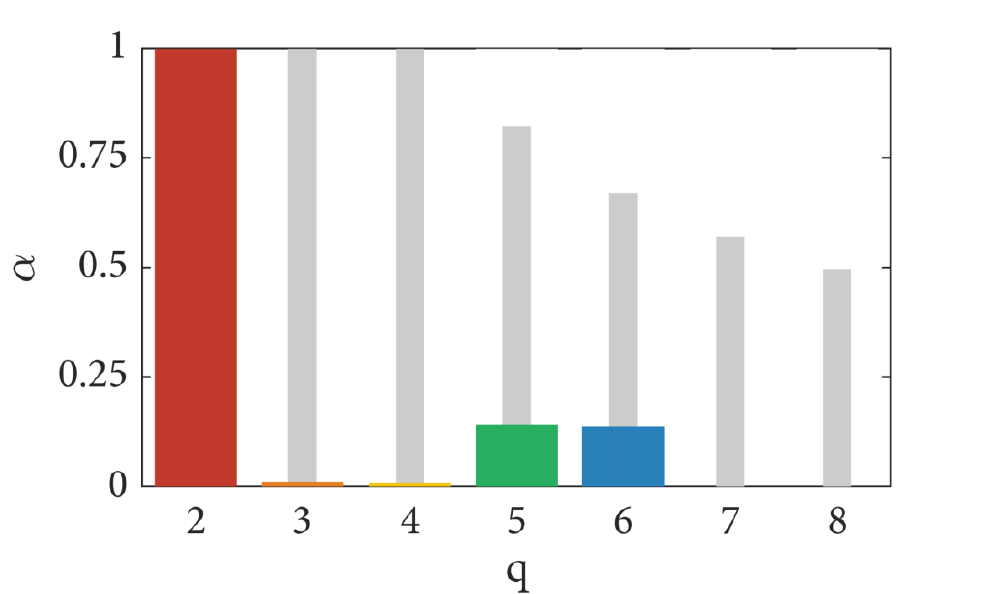}	 &
					\includegraphics[valign=m, width=1\linewidth]{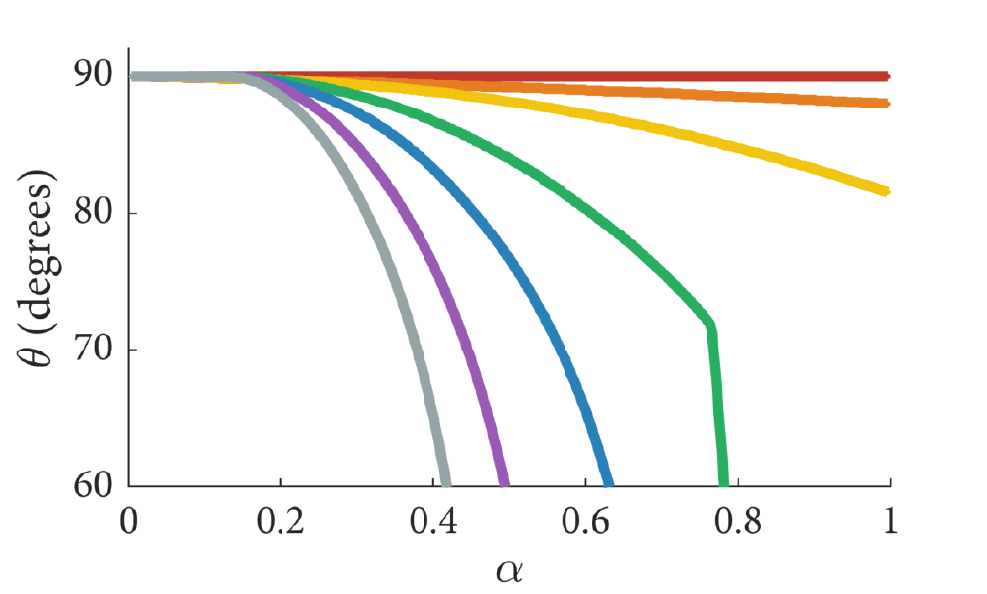} \\[1em]
				\AMDescTable{PMFCmj}{Adams}{FI, $\ell = 3$}{Cheb}{Inwards} &
					\includegraphics[valign=m,width=1\linewidth]{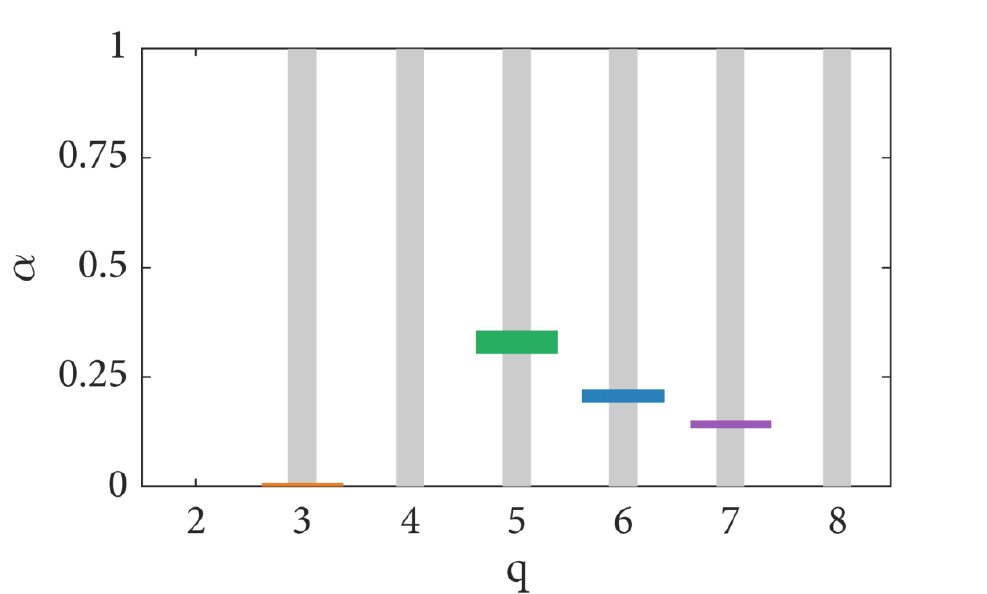}	 &
					\includegraphics[valign=m, width=1\linewidth]{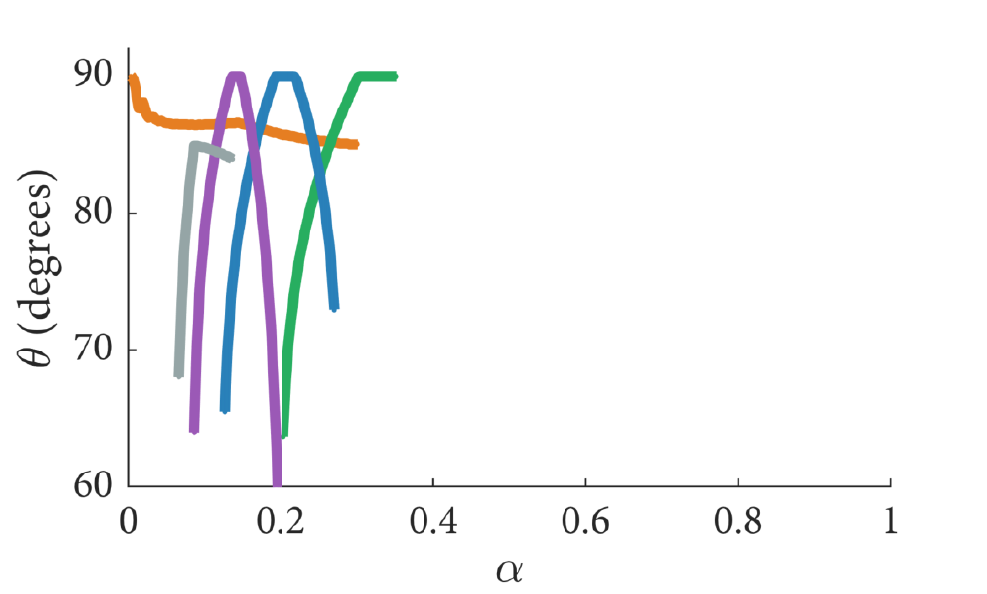} \\[1em]
				\AMDescTable{SMFCmj}{Adams}{FI, $\ell = 3$}{Cheb}{Inwards} &
					\includegraphics[valign=m,width=1\linewidth]{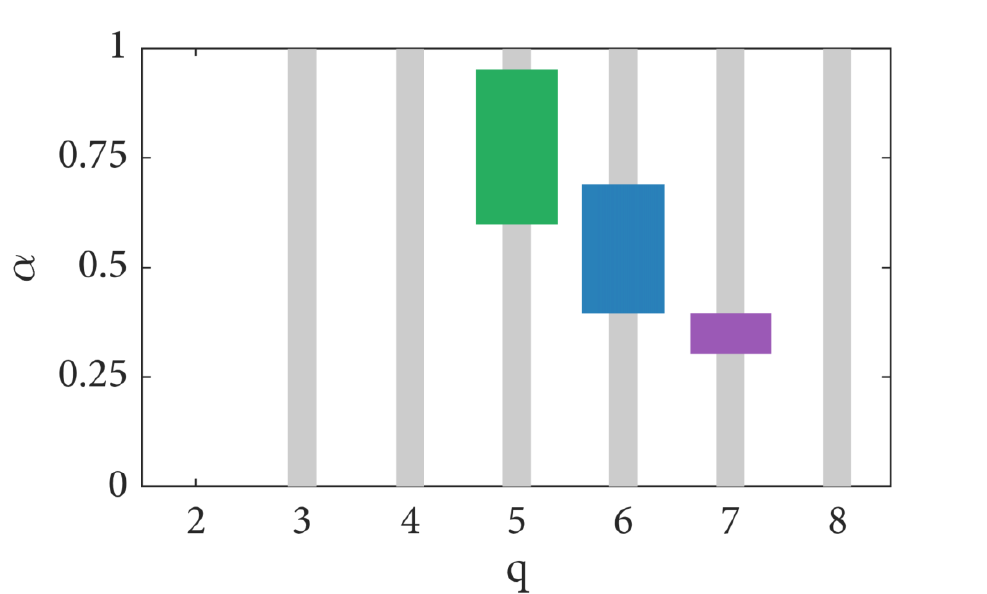}	 &
					\includegraphics[valign=m, width=1\linewidth]{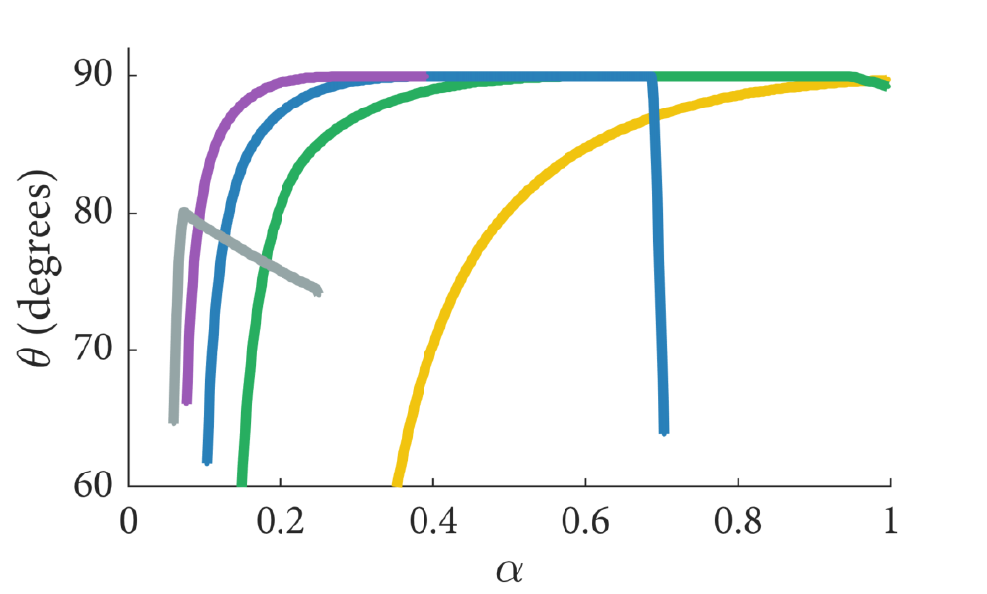} \\[1em]
				\AMDescTable{SMFCmj}{Adams}{VI}{Cheb}{Inwards} &
					\includegraphics[valign=m,width=1\linewidth]{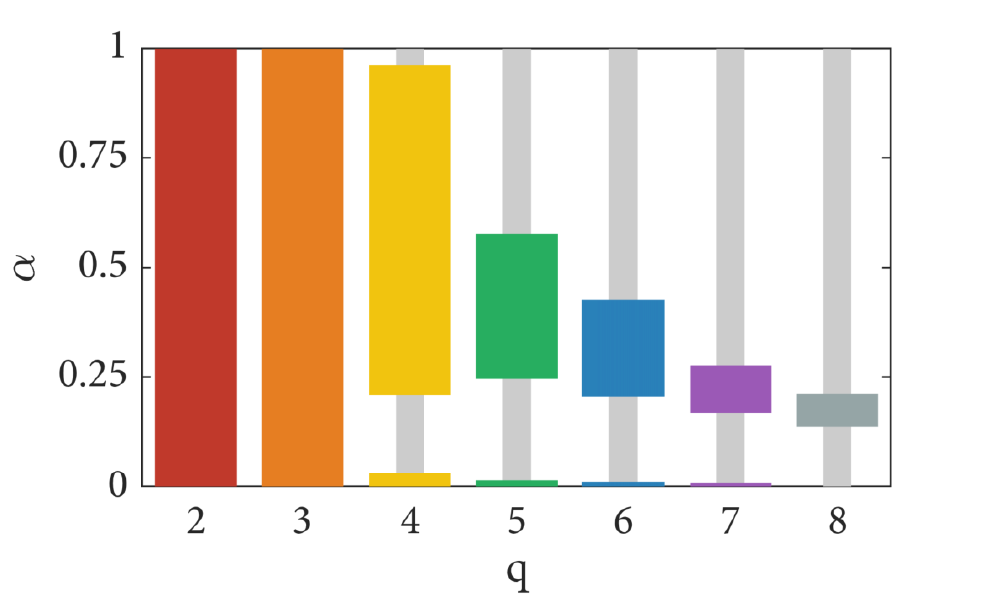}	 &
					\includegraphics[valign=m, width=1\linewidth]{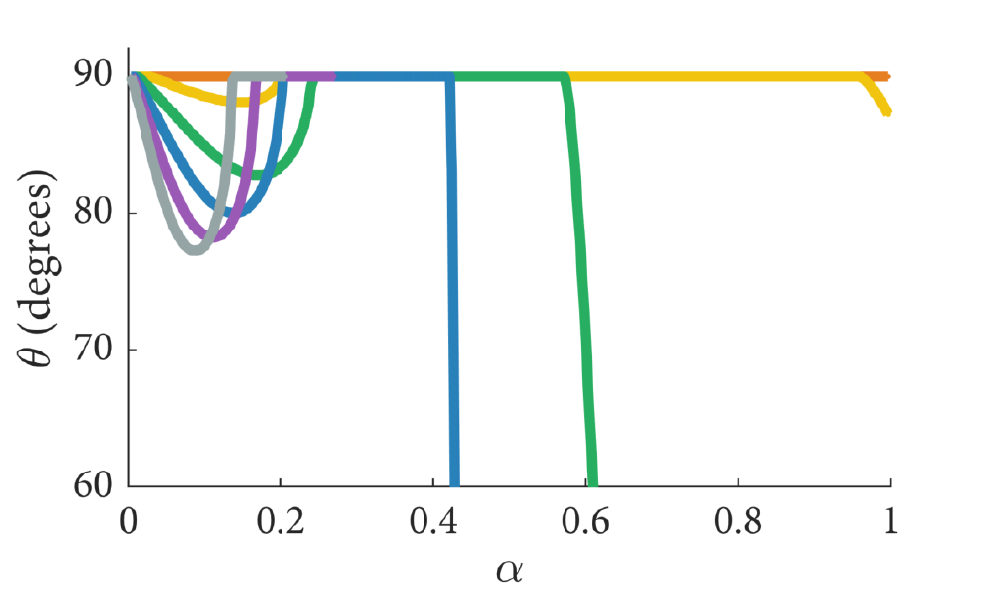} \\[1em]
				
			\end{tabularx}
		
			\caption{Families of PBMs possessing $A(90\degree)$ stability. The bar charts show $A(90\degree)$ stability and root stability. Thick colored bars show the $\alpha$ values for which a method is $A(90\degree)$ stable and light gray bars show the $\alpha$ values for which the method is root stable. The plots in the right column show the stability angle $\theta$ as a function of the extrapolation factor $\alpha$. Color represents the number of nodes $q$.}	
			\label{fig:DIPBM_AofTheta_vs_alpha}
			
		\end{figure}

		\begin{figure}[h!]
		
			\centering

			\begin{tabular}{lllllll}
				\lineLegend{plot_red}{$q=2$} &
				\lineLegend{plot_orange}{$q=3$} &
				\lineLegend{plot_yellow}{$q=4$} &
				\lineLegend{plot_green}{$q=5$} &
				\lineLegend{plot_blue}{$q=6$} &
				\lineLegend{plot_violet}{$q=7$} &
				\lineLegend{plot_grey}{$q=8$} \hspace{1.5em} 
			\end{tabular}

			\begin{tabular}{cc}
				\includegraphics[valign=m,width=0.35\linewidth]{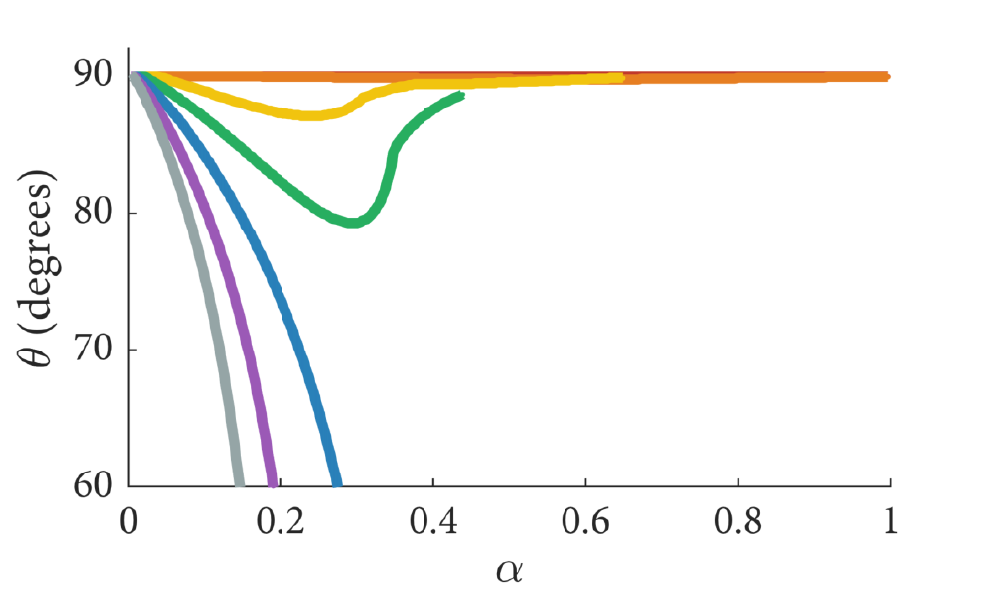}	 &
				\includegraphics[valign=m,width=0.35\linewidth]{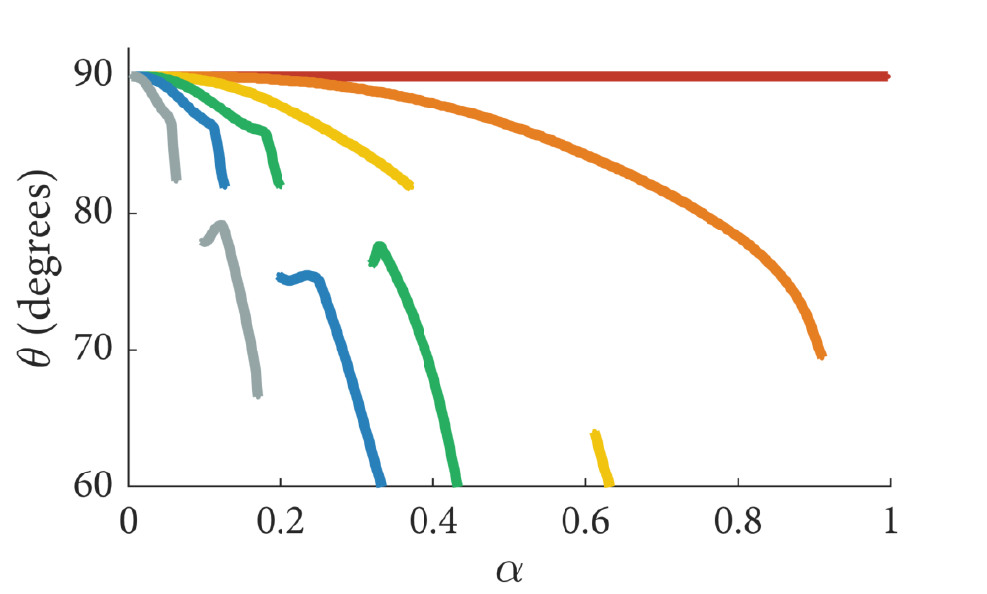} \\
				\\
				\begin{scriptsize}
					\begin{tabular}{ll}
						{\bf AII/AOI:} & PMFCmj \\
						{\bf ODEP:} & Adams \\
						{\bf Nodes:} & Equi, inwards \\
						{\bf Endpoints:} & VI
					\end{tabular}
				\end{scriptsize} &
				\begin{scriptsize}
					\begin{tabular}{ll}
						{\bf AII/AOI:} & SMVC \\
						{\bf ODEP:} & Adams \\
						{\bf Nodes:} & Cheb, inwards \\
						{\bf Endpoints:} & FI, $\ell=2$
					\end{tabular}
				\end{scriptsize} \\[1em]			
			\end{tabular}
		
			\caption{Stability angles for two families of Adams PBMS with $A(\theta > 0)$ stability. Several of these methods will be used in our numerical experiments section.}	
			\label{fig:Adams_stability_angles}
			
		\end{figure}

	\begin{table}[h]
		\centering
		\renewcommand*{\arraystretch}{1.5}
			
		\begin{tabularx}{\linewidth}{XX}
			
			Inwards Sweeping & Outwards Sweeping \\[1em]
			
			\begin{tabular}{llll} \hline
				{\em AIS}    & {\em EP} & {\em q} \\ \hline
				PMFCmj & fixed, $\ell = 1$ & $2-4, *$ \\
					   & sweeping & $2,3,4,5,*$ \\
					   & variable & $2,3,4,5,*$ \\ \hline
				SMFCmj & fixed,  $\ell = 2$ & $2,3,4,5, *$ \\
			 		   & sweeping & $2,3,4,5,*$ \\
			 		   & variable & $2,3,4,5,*$ \\ \hline
				SMFC & fixed, $\ell = 2$ & $2,3$ \\
				     & sweeping & $*$ \\ \hline
				SMVC & fixed, $\ell = 2$ & $2,3$
			\end{tabular}
			&
			\begin{tabular}{lll} \hline
				{\em AIS}    & {\em EP} & {\em q} \\ \hline
				PMFCmj & sweeping & $2,3,4,5,*$ \\	 \\ \\ \hline
				SMFCmj & fixed,  $\ell = 1$ & $2,3,4,5$ \\
			 		   & sweeping & $2,3,4,5, *$  \\
			 		   & variable & $2,3,4,5,*$ \\ \hline
				\\ \\ \\
			\end{tabular}
		\end{tabularx}
						
		\caption{A non-exhaustive list of construction strategies for Adams PBMs using imaginary equispaced nodes with $A(\theta > 0)$ stability for all ${0 < \alpha \le 4/10}$. A star indicates that higher-order methods with $A(\theta > 0)$ exist for a different range of $\alpha$ between zero and one.}
		\label{tab:Example_Adams_methods_equi}
	\end{table}

	\section{A Matlab package for computing coefficients and linear stability}
	\label{sec:pipack} 

	Though we have introduced a wide range of construction strategies, it would be unreasonable to try and discuss them all in this work. We are also faced with the challange of presenting method coefficients, which would fill many pages if we were to include them for all families of PBMs. We therefore created the Matlab library {\em Polynomial Integrator Package} (PIPack) \cite{buvoli2020PIPACK}. This library can be used to easily generate coefficients and analyze the linear stability properties for any PBM that is constructed using any of the strategies proposed in Section \ref{sec:constructing_pbms}. In particular, the code can:
	\begin{enumerate}[leftmargin=*]
		\item Generate PBM coefficients in double precision, variable precision, or symbolically.
		\item Numerically compute stability properties like stability angle and root stability.
		\item Generate linear stability plots and movies that show how linear stability regions change as function of the extrapolation parameter $\alpha$.
		\item Generate \ODEp{} diagrams and expansion point diagrams for any method.
		\item Run the numerical examples presented in the next section of this paper.
	\end{enumerate}
	In the supplemental materials section \ref{sup_sec:example_diagrams} we show two example polynomial diagrams for BDF-SMFC and Adams-PMFCmj that were generated using PIPack. We also include the example code for generating the figures and the corresponding method coefficients. Further documentation and examples are shown on PIPACK's Github repository \cite{buvoli2020PIPACK}.
	
	Finally, we briefly discuss coefficient computation for anyone who wants to write their own code for initializing the matrices $\mathbf{A}(\alpha)$, $\mathbf{B}(\alpha)$, $\mathbf{C}(\alpha)$, and $\mathbf{D}(\alpha)$ in (\ref{eq:block_parametrized}). An algorithm for computing the coefficients of PBMs constructed without an \ivs{} is presented in \cite{buvoli2019constructing} and a more general procedure is described in \cite{buvoli2018polynomial}. %

		\section{Numerical experiments}
		\label{sec:numerical_experiments}
		
		Our numerical experiments are designed to evaluate the performance of the newly introduced PBMs on dispersive equations, and to compare new serial PBMs against the original BAM and BBDF methods from \cite{buvoli2019constructing}. We conduct the numerical experiments by solving two stiff partial differential equations with a variety of diagonally implicit PBMs. For reference, we also compare our methods against classical backwards difference methods (BDF) of orders two through six, and the implicit Runge-Kutta methods from \cite{kennedy2019diagonally} listed in Table \ref{tab:ESDIRK_methods}.
		
		To test methods with A($90\degree$) stability angles, we solve the nonlinear Schr\"{o}dinger (NLS) equation using four of the seven methods described in Figure \ref{fig:DIPBM_AofTheta_vs_alpha}. In Table \ref{tab:NLS_methods} we show the specific $\alpha$ values used in the NLS numerical experiments. To compare the performance of parallel and serial methods we repeat the numerical experiment with Burgers' equation from \cite{buvoli2019constructing} using the PBMS listed in Table \ref{tab:BRG_methods}. Since we are interested in methods for solving stiff PDEs, we only test methods with unbounded stability regions.
		
		We present plots of absolute error versus stepsize and absolute error versus computational time where errors are measured using the $L^{\infty}$ norm. For each equation, we compute the reference solutions and initial values using MATLAB's {\em ode15s} integrator with a tolerance of $1$e-$14$. We solve all nonlinear systems using Newton's method with an exact Jacobian, and MATLAB's backslash to solve the underlying linear systems. All the results presented in this paper have been run on a six-core Intel\textsuperscript{\textregistered} i7-8700 CPU running at 3.70GHz. Serial PBMs were run using a single thread, while parallel PBMs where parallelized using MATLAB's parallel toolkit.  
		
		\subsection{Problems} We briefly describe the partial differential equations used in our numerical experiment including their initial conditions, and spatial discretizations.

\begin{figure}[ph!]
			\centering
			\hfill
			\begin{tabular}{lllllll}
				\lineLegend{plot_red}{$q=2$} &
				\lineLegend{plot_orange}{$q=3$} &
				\lineLegend{plot_yellow}{$q=4$} &
				\lineLegend{plot_green}{$q=5$} &
				\lineLegend{plot_blue}{$q=6$} &
				\lineLegend{plot_violet}{$q=7$} &
				\lineLegend{plot_grey}{$q=8$} \hspace{1.5em} 
			\end{tabular}
			\begin{tabularx}{\linewidth}{mff}
				\AMDescTable{PMFCmj}{Adams}{FI, $\ell = 3$}{Cheb}{Inwards} &
					\includegraphics[valign=m,width=1\linewidth]{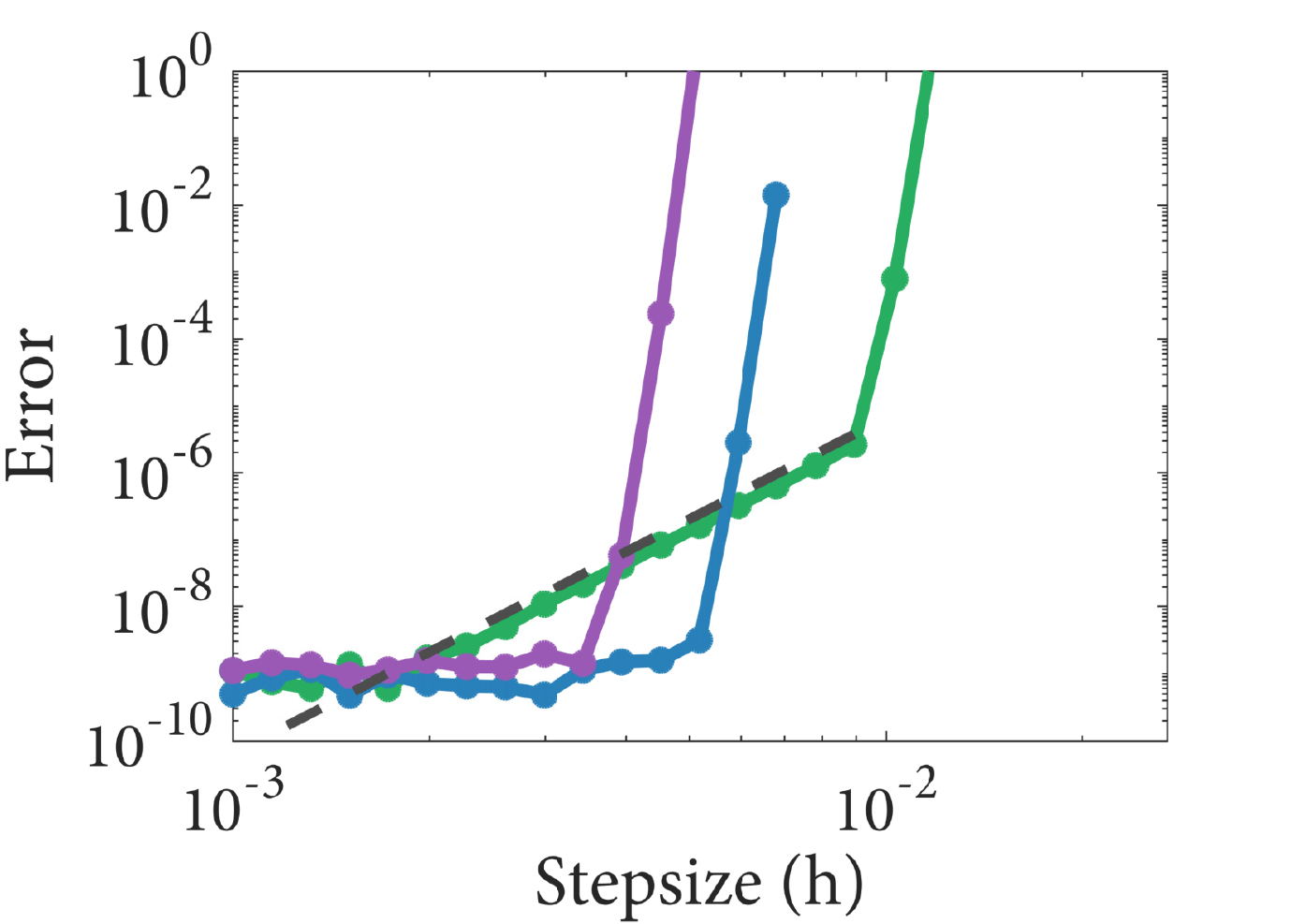}	 &
					\includegraphics[valign=m, width=1\linewidth]{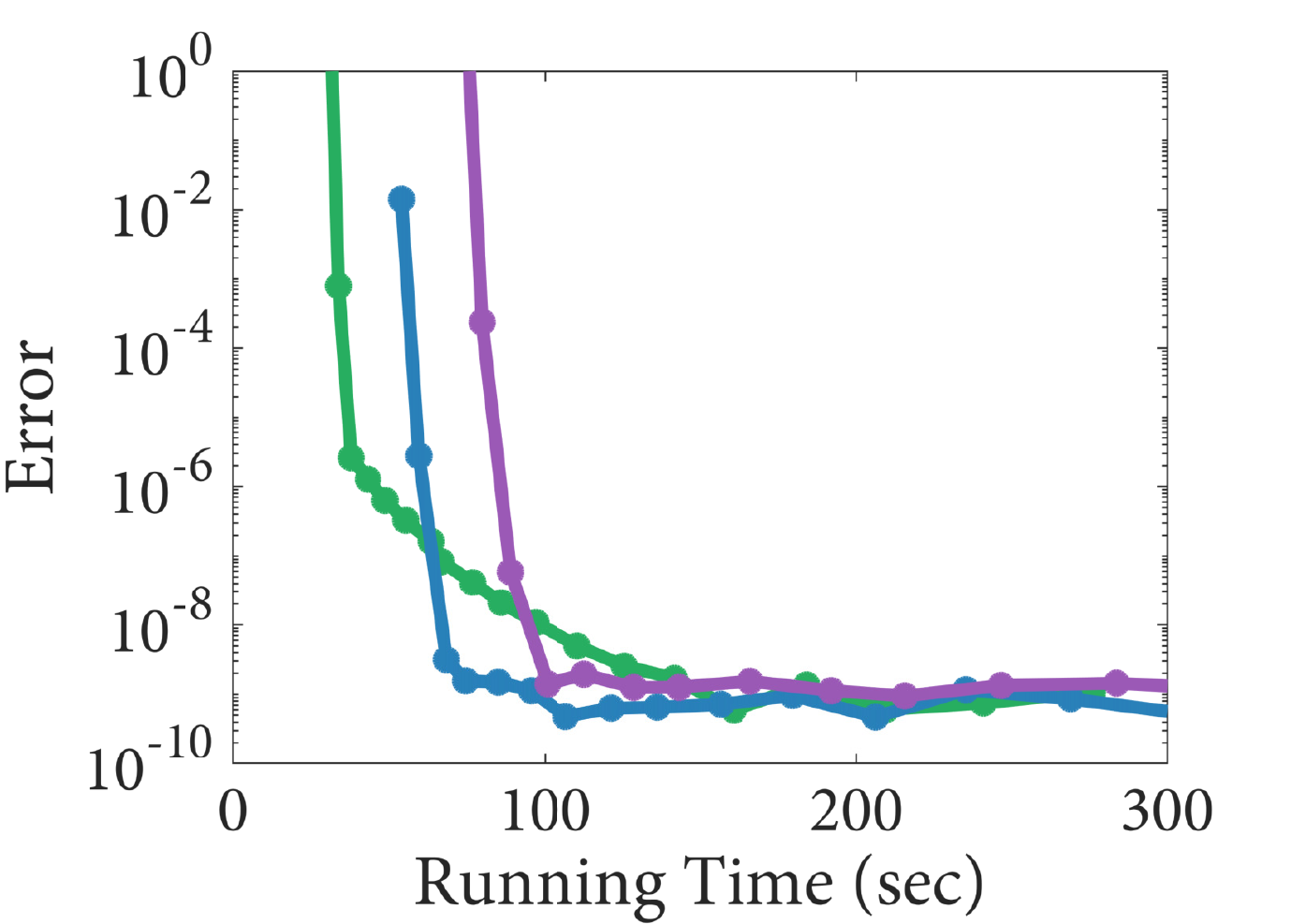} \\[1em]
				\AMDescTable{SMFCmj}{Adams}{FI, $\ell = 3$}{Cheb}{Inwards} &
					\includegraphics[valign=m,width=1\linewidth]{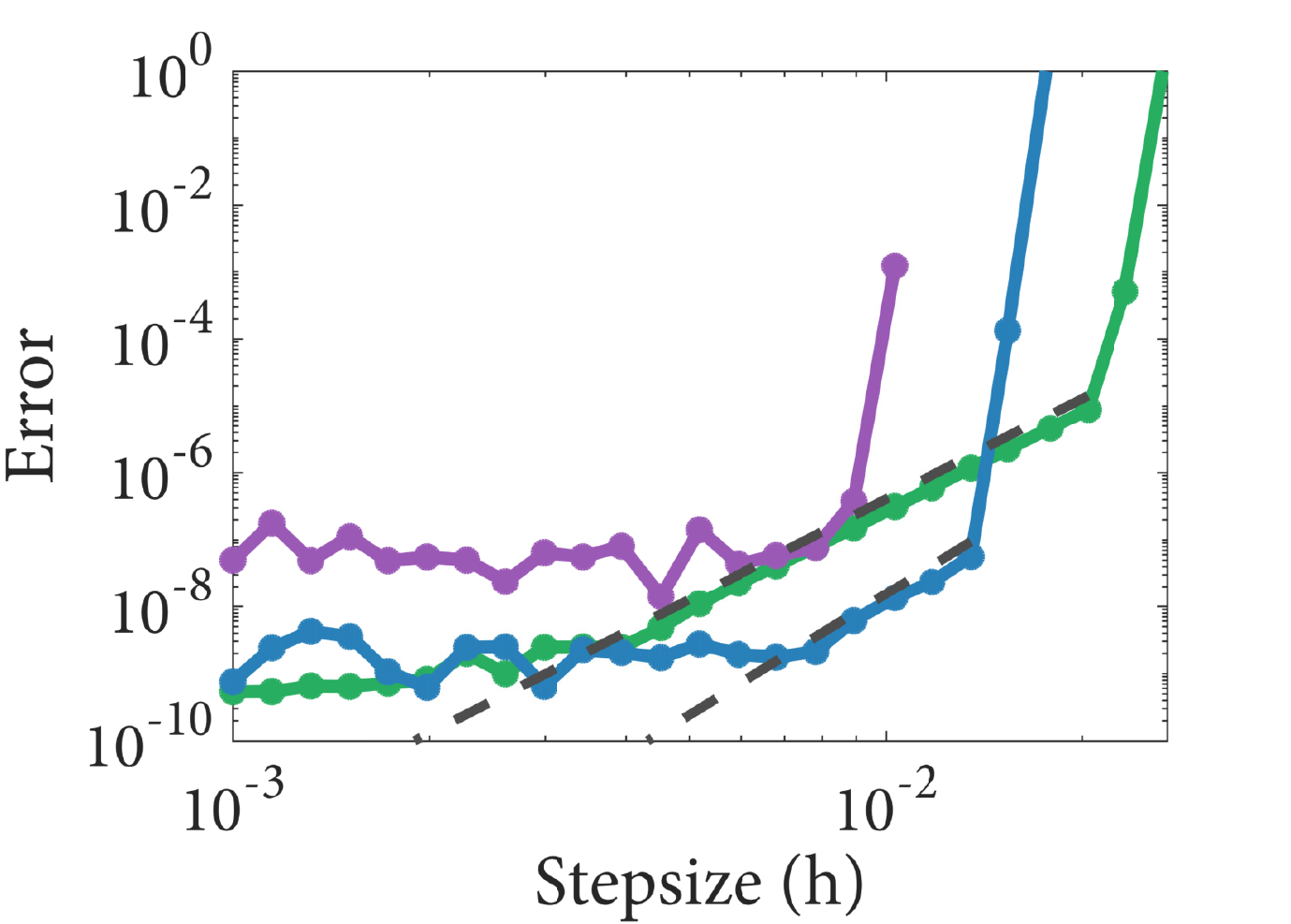}	 &
					\includegraphics[valign=m, width=1\linewidth]{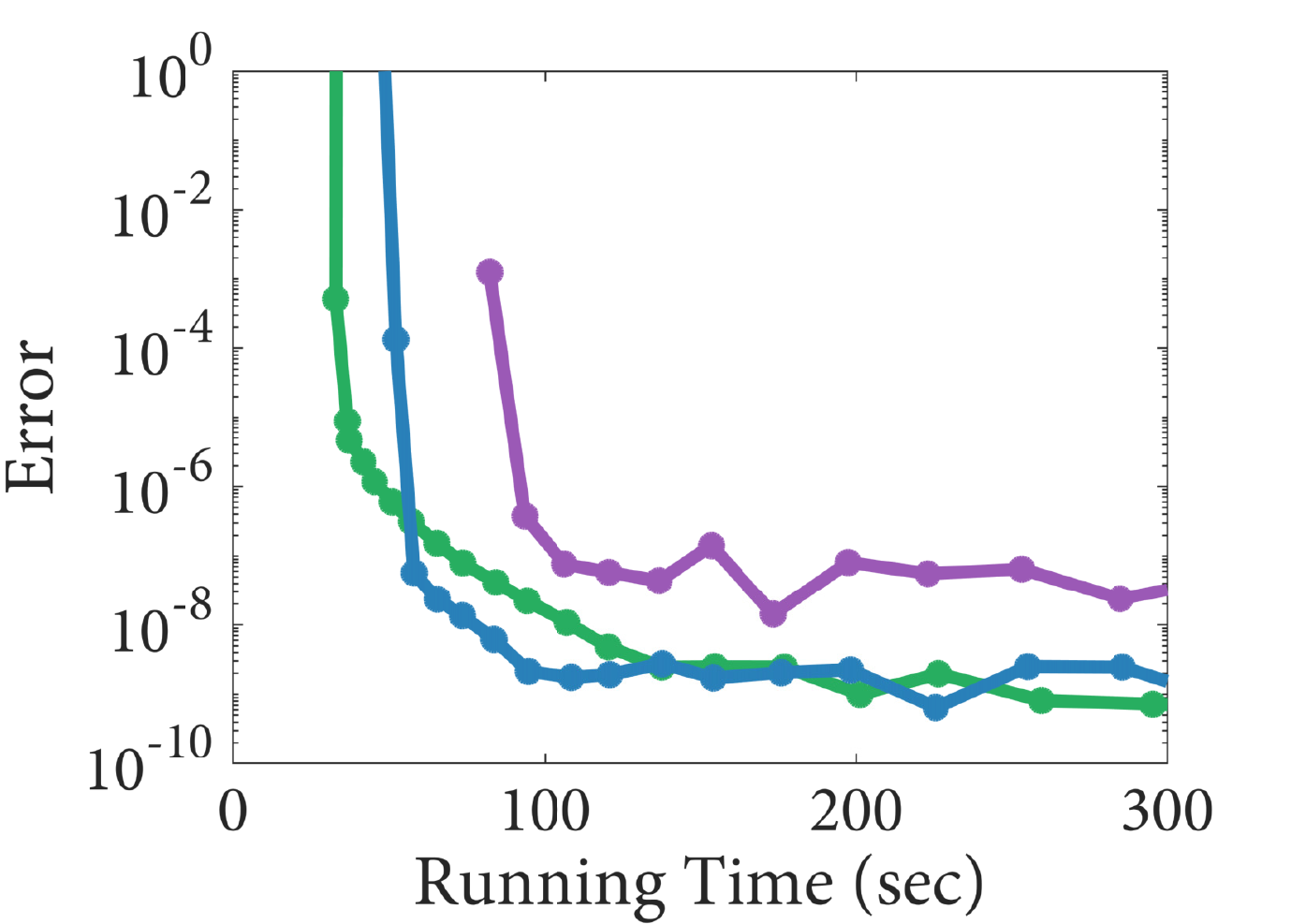} \\[1em]
				\AMDescTable{SMFCmj}{Adams}{VI}{Cheb}{Inwards} &
					\includegraphics[valign=m,width=1\linewidth]{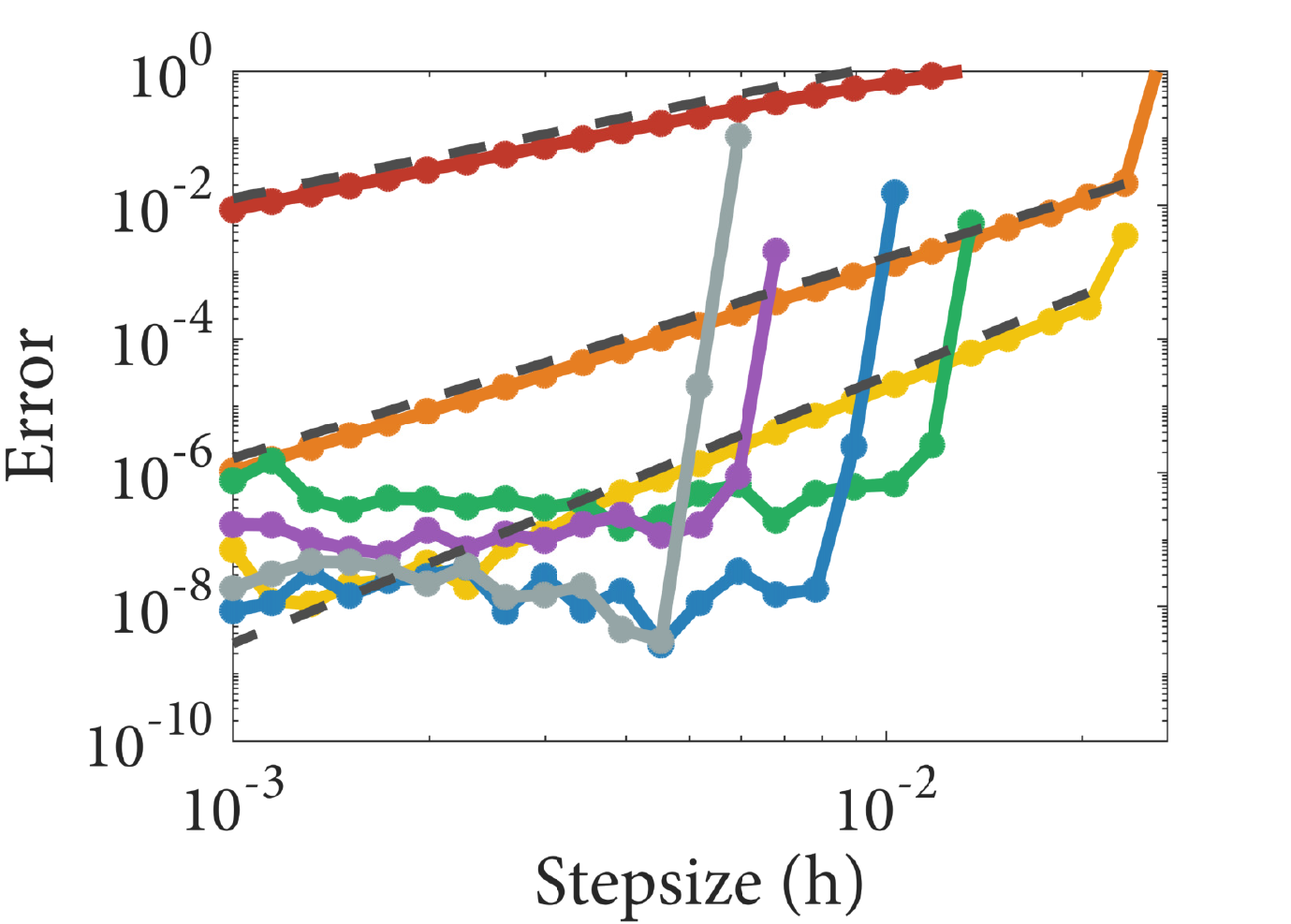}	 &
					\includegraphics[valign=m, width=1\linewidth]{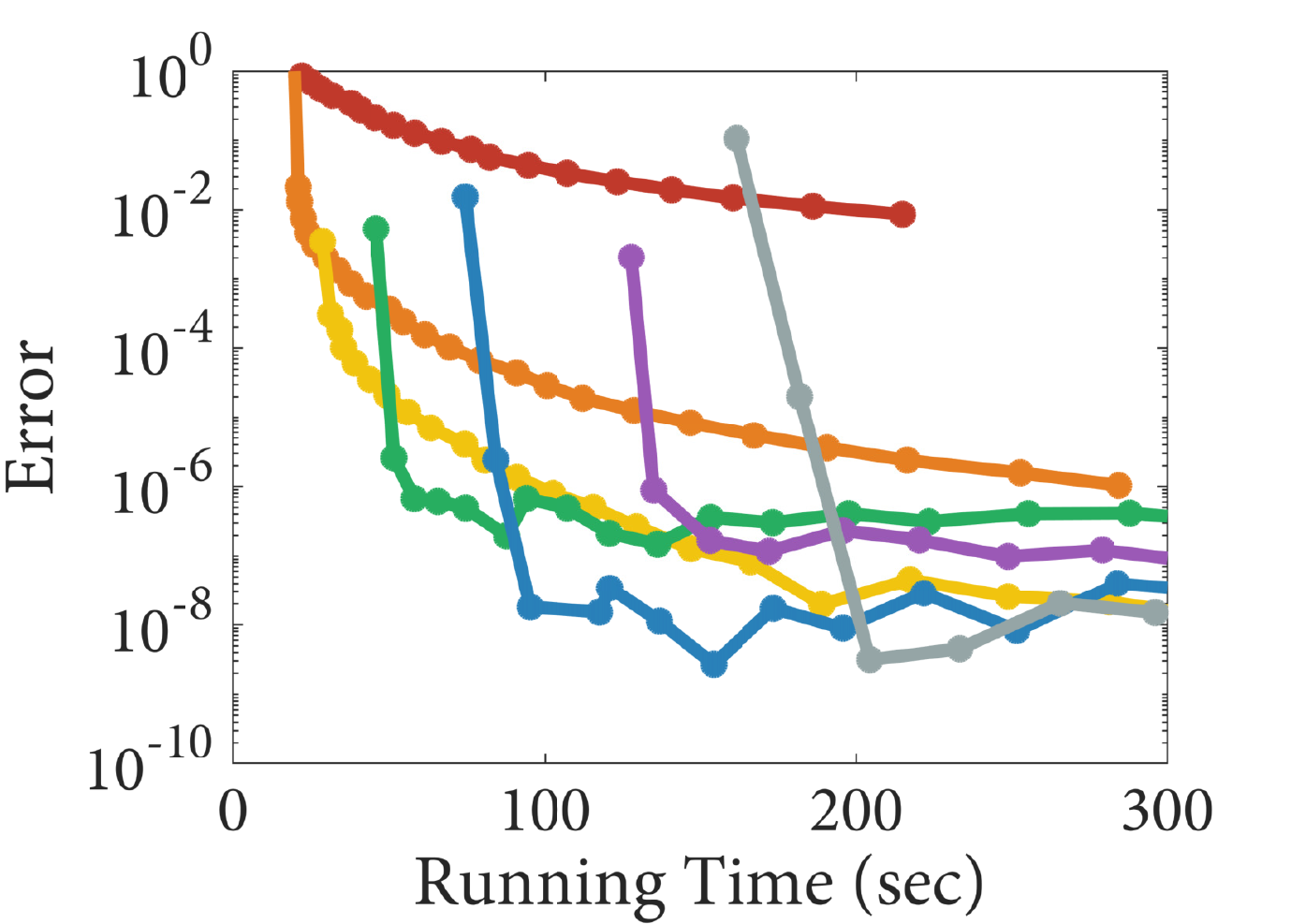} \\[1em]
				\GBDFDescTable{SMFC}{BDF}{Equi}{Inwards} &
					\includegraphics[valign=m,width=1\linewidth]{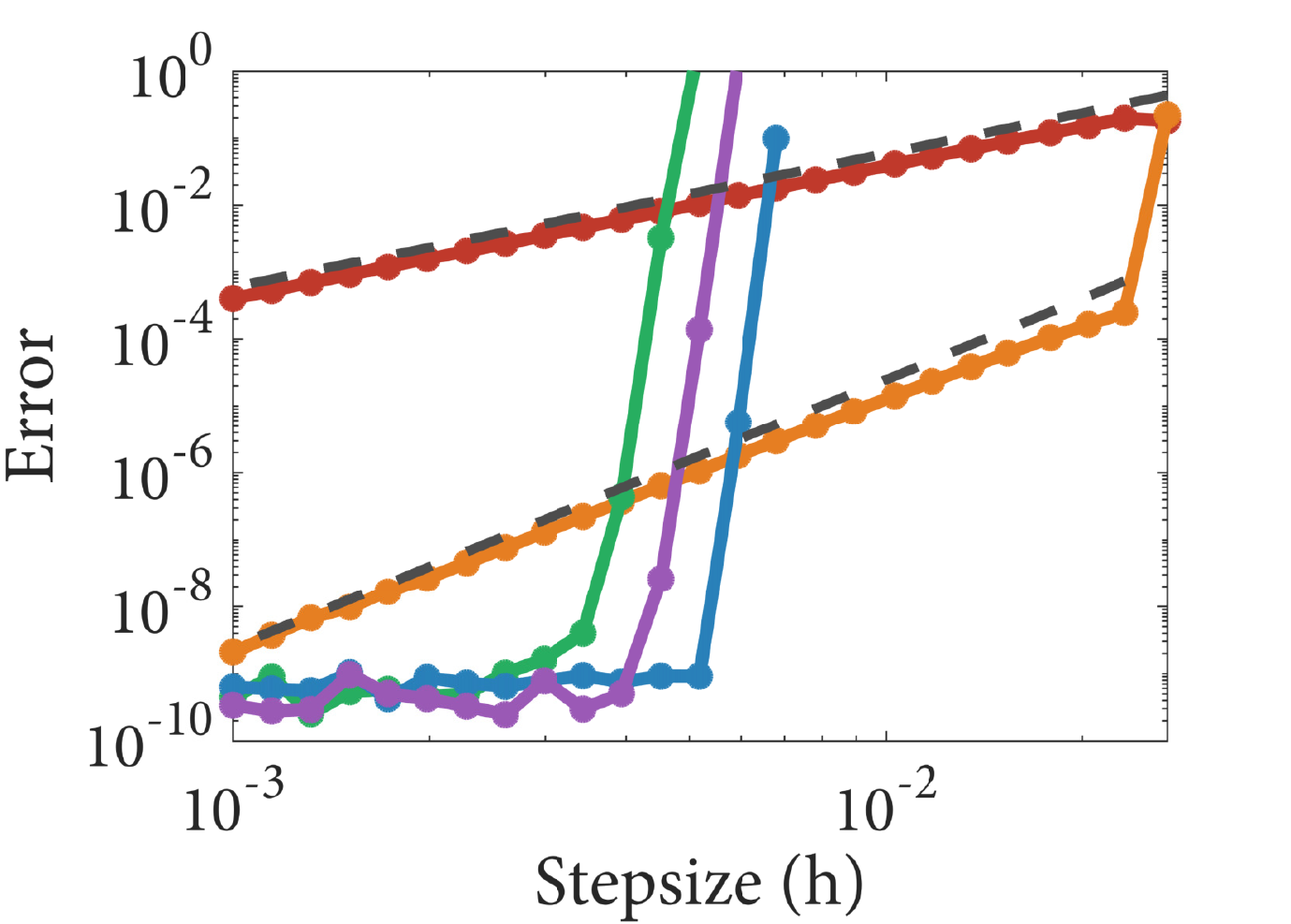}	 &
					\includegraphics[valign=m, width=1\linewidth]{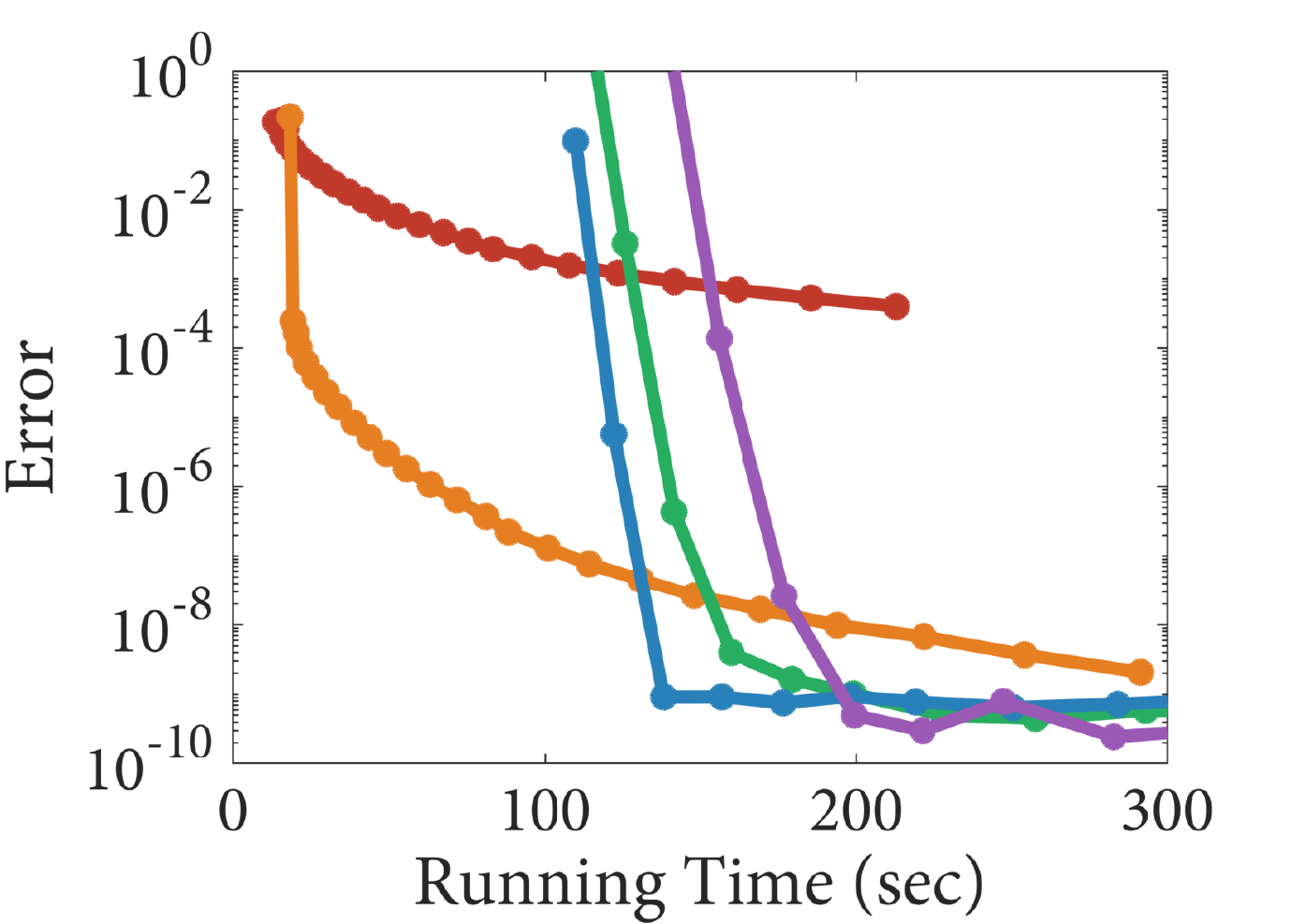} \\[1em]
			\end{tabularx}
			\begin{tabularx}{\linewidth}{llll}
				\hspace{9.5em}
				\lineLegend{plot_orange}{$s=3$} &
				\lineLegend{plot_green}{$s=5$} &
				\lineLegend{plot_violet}{$s=7$} &
				\lineLegend{plot_grey}{$s=8$} 
			\end{tabularx}
			\begin{tabularx}{\linewidth}{mff}
				\MethodNameTable{ESDIRK} & \includegraphics[valign=m,width=0.35\textwidth]{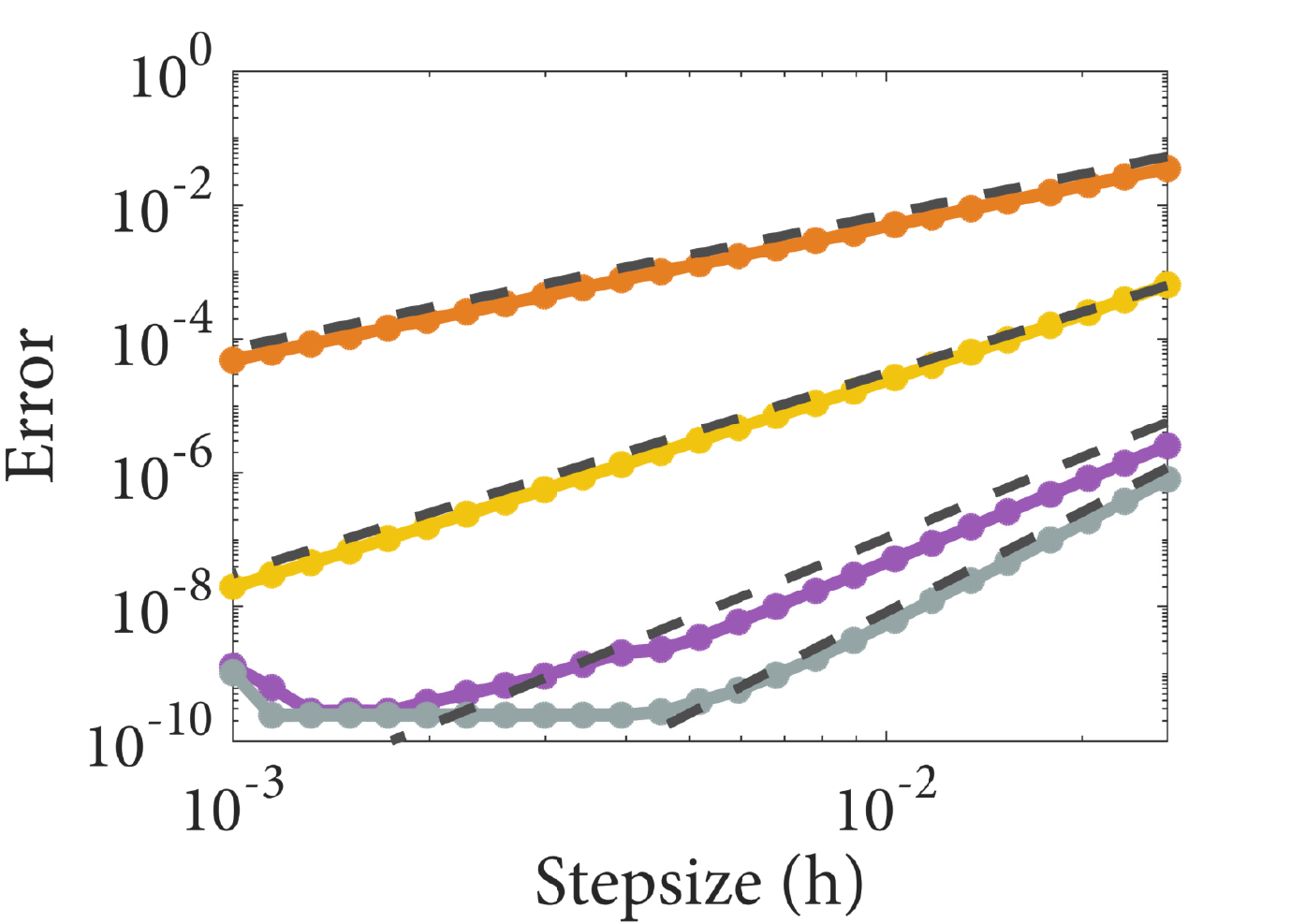}	 &
				\includegraphics[valign=m,width=0.35\textwidth]{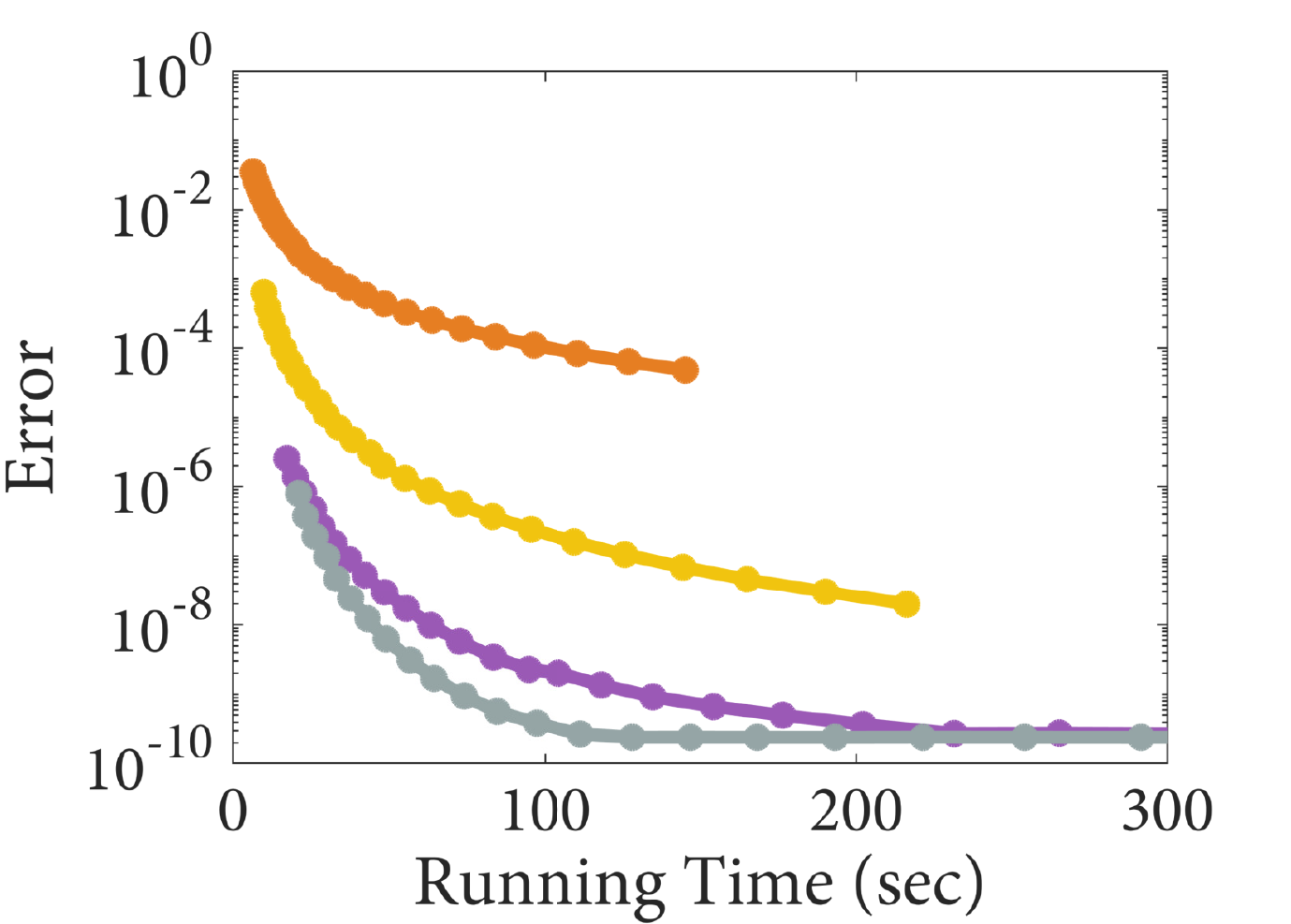} \\[1em]	
				
			\end{tabularx}
			\caption{Accuracy and precision diagrams for the nonlinear Schr\"{o}dinger equation solved using four new PBMs and the ESDIRK methods from \cite{kennedy2019diagonally}. Color is used to represent the number of nodes for polynomial methods and the number of implicit stages $s$ for ESDIRK methods. The BDF-SMFC with $q=3$ converged at fourth order, while all other PBMs with $q$ nodes converged to $qth$ order. ESDIRK methods converged at third, fourth, fifth, and sixth order.}	
			\label{fig:NLS_results}
		\end{figure}

\begin{figure}
			\centering
			\hfill
			\begin{tabular}{lllllll}
				\lineLegend{plot_red}{$q=2$} &
				\lineLegend{plot_orange}{$q=3$} &
				\lineLegend{plot_yellow}{$q=4$} &
				\lineLegend{plot_green}{$q=5$} &
				\lineLegend{plot_blue}{$q=6$} &
				\lineLegend{plot_violet}{$q=7$} &
				\lineLegend{plot_grey}{$q=8$} \hspace{1.5em} 
			\end{tabular}
			\vspace{1em}

			\begin{tabularx}{\linewidth}{mff}
				\GBDFDescTable{SMFC}{BDF}{Equi}{Inwards}	 & 
					\includegraphics[valign=m,width=1\linewidth]{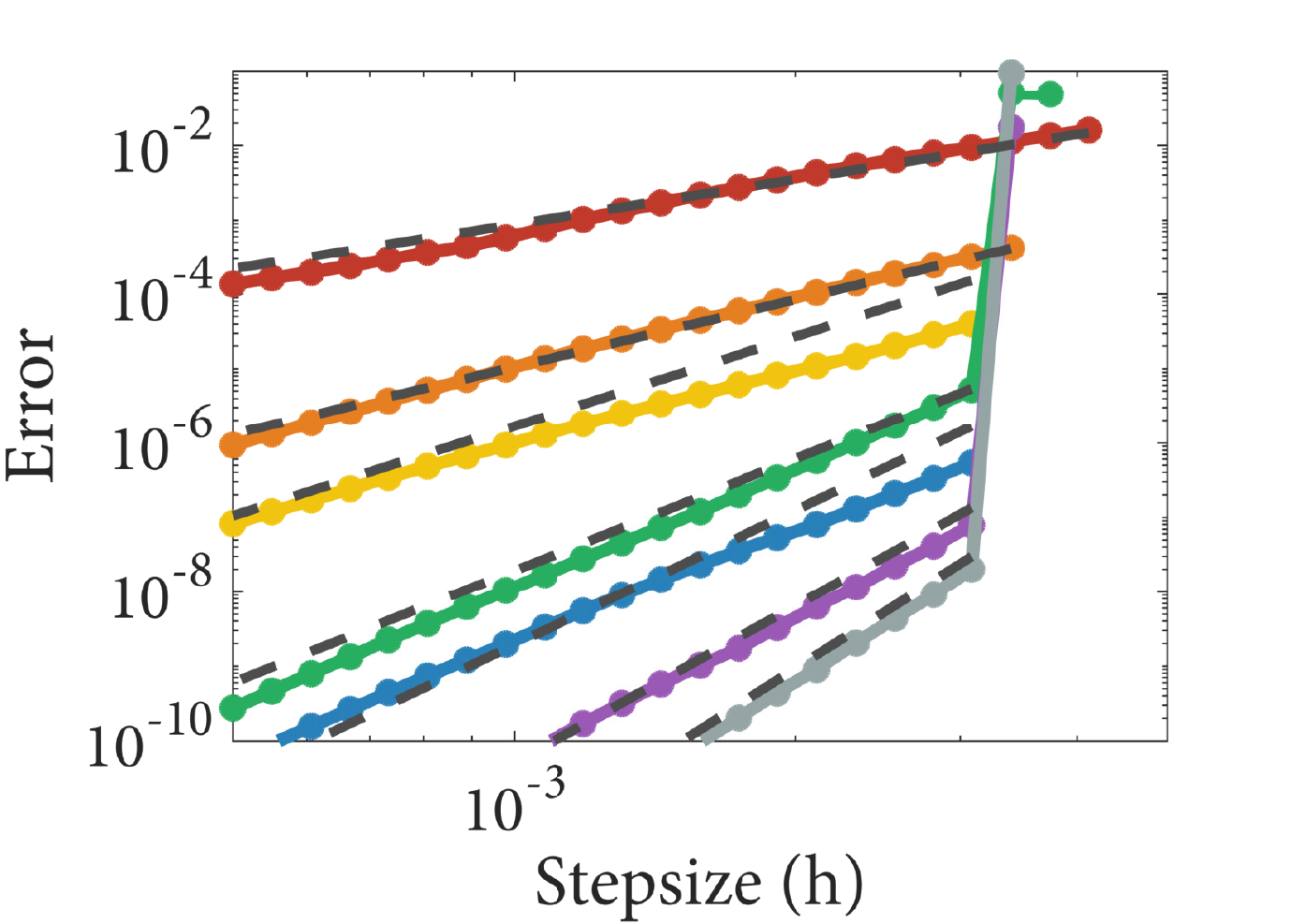}	 &
					\includegraphics[valign=m,width=1\linewidth]{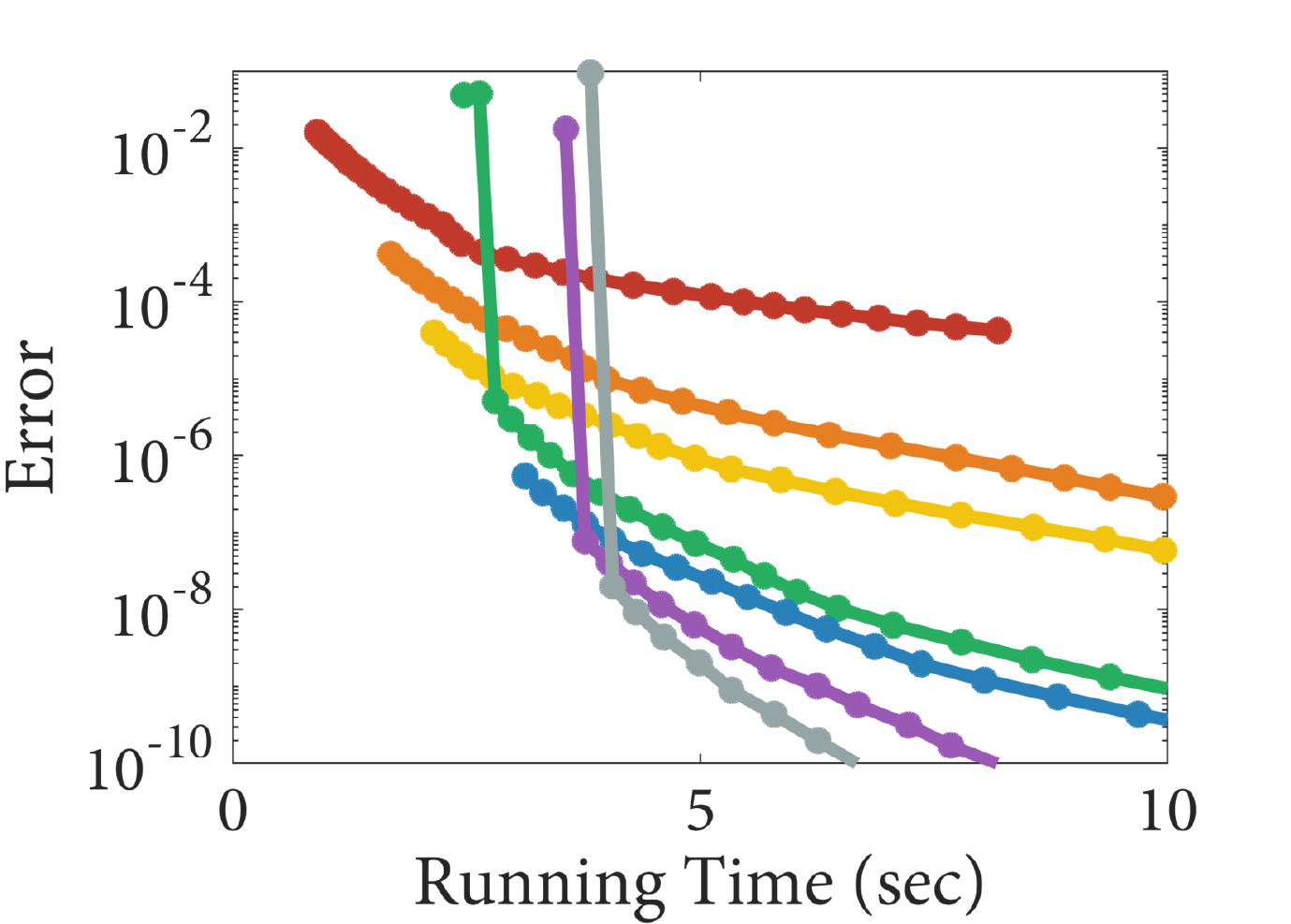} \\[1em]
				\AMDescTable{PMFCmj}{Adams}{Variable Input}{Equi}{Inwards} &
					\includegraphics[valign=m,width=1\linewidth]{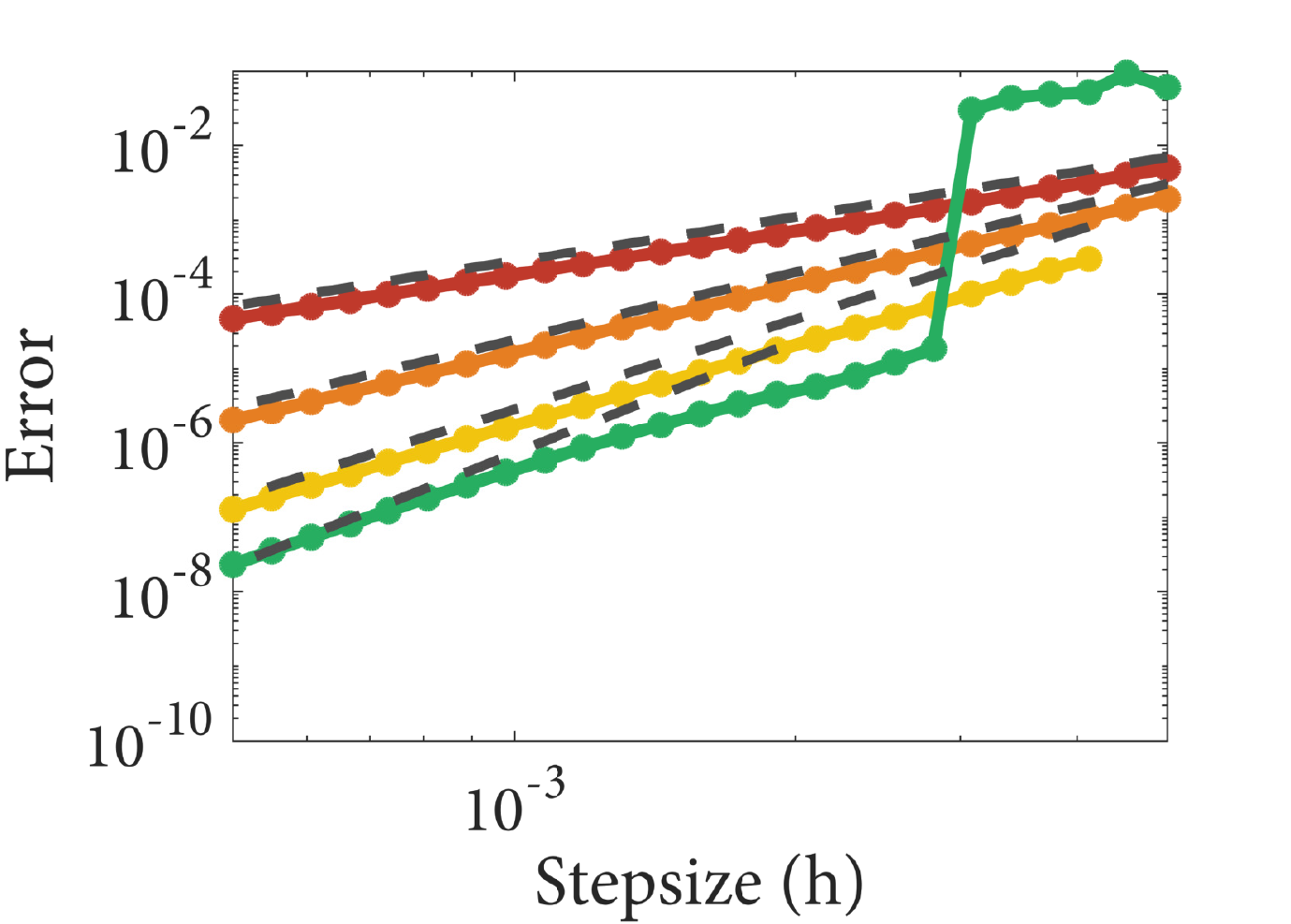}	 &
					\includegraphics[valign=m, width=1\linewidth]{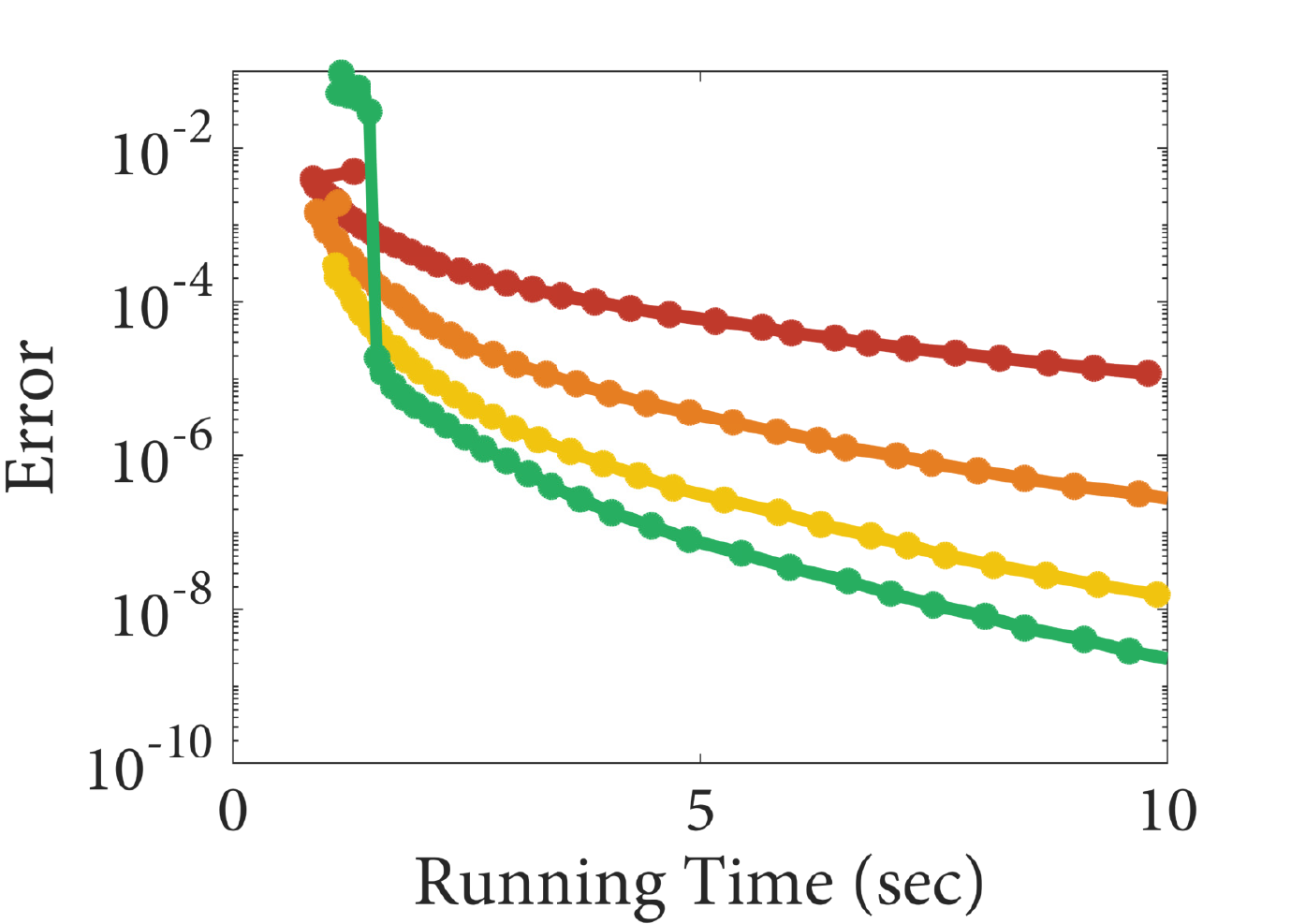} \\[1em]
				\AMDescTable{SMVC}{Adams}{FI, $\ell = 2$}{Cheb}{Inwards} &
					\includegraphics[valign=m,width=1\linewidth]{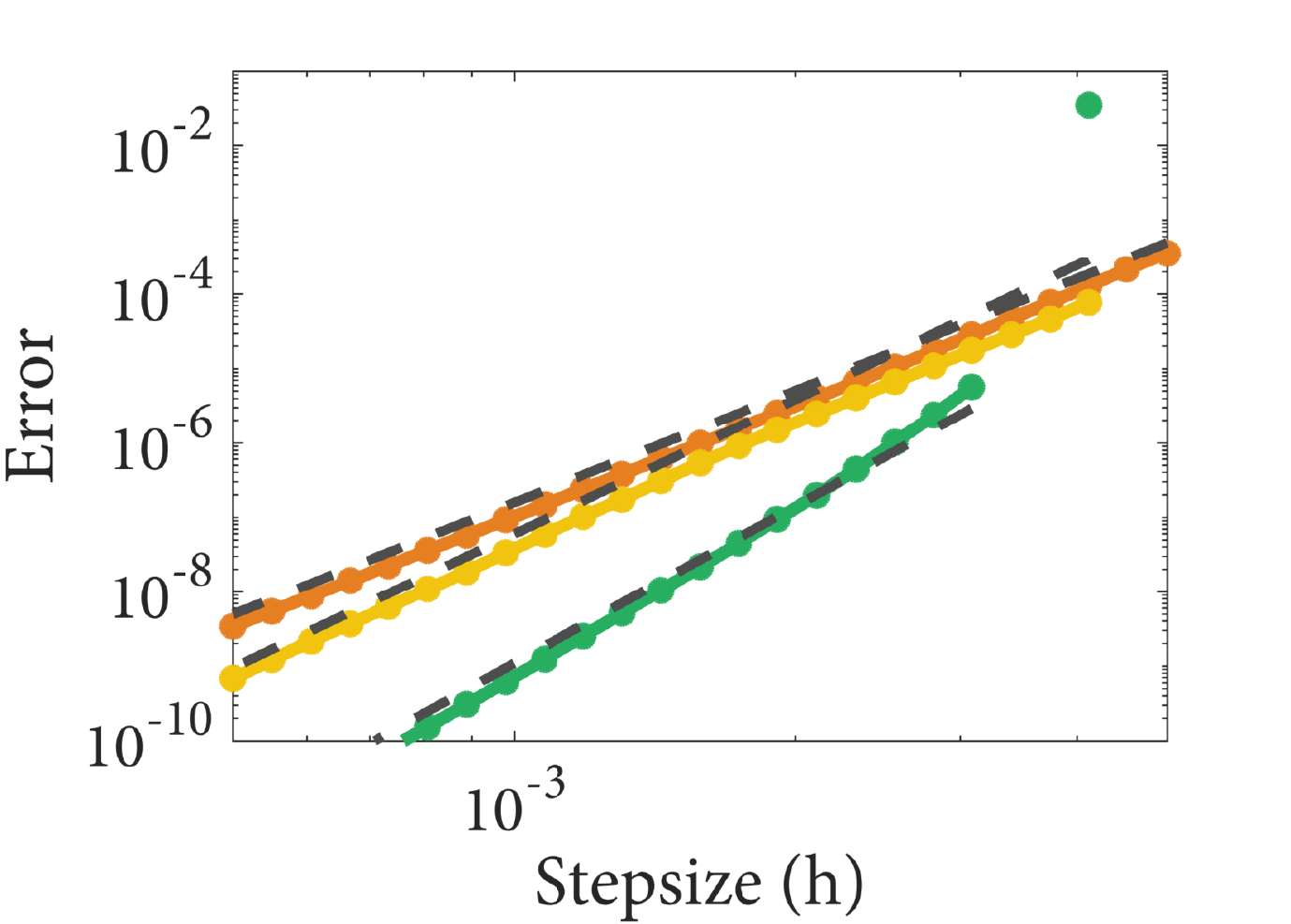}	 &
					\includegraphics[valign=m,width=1\linewidth]{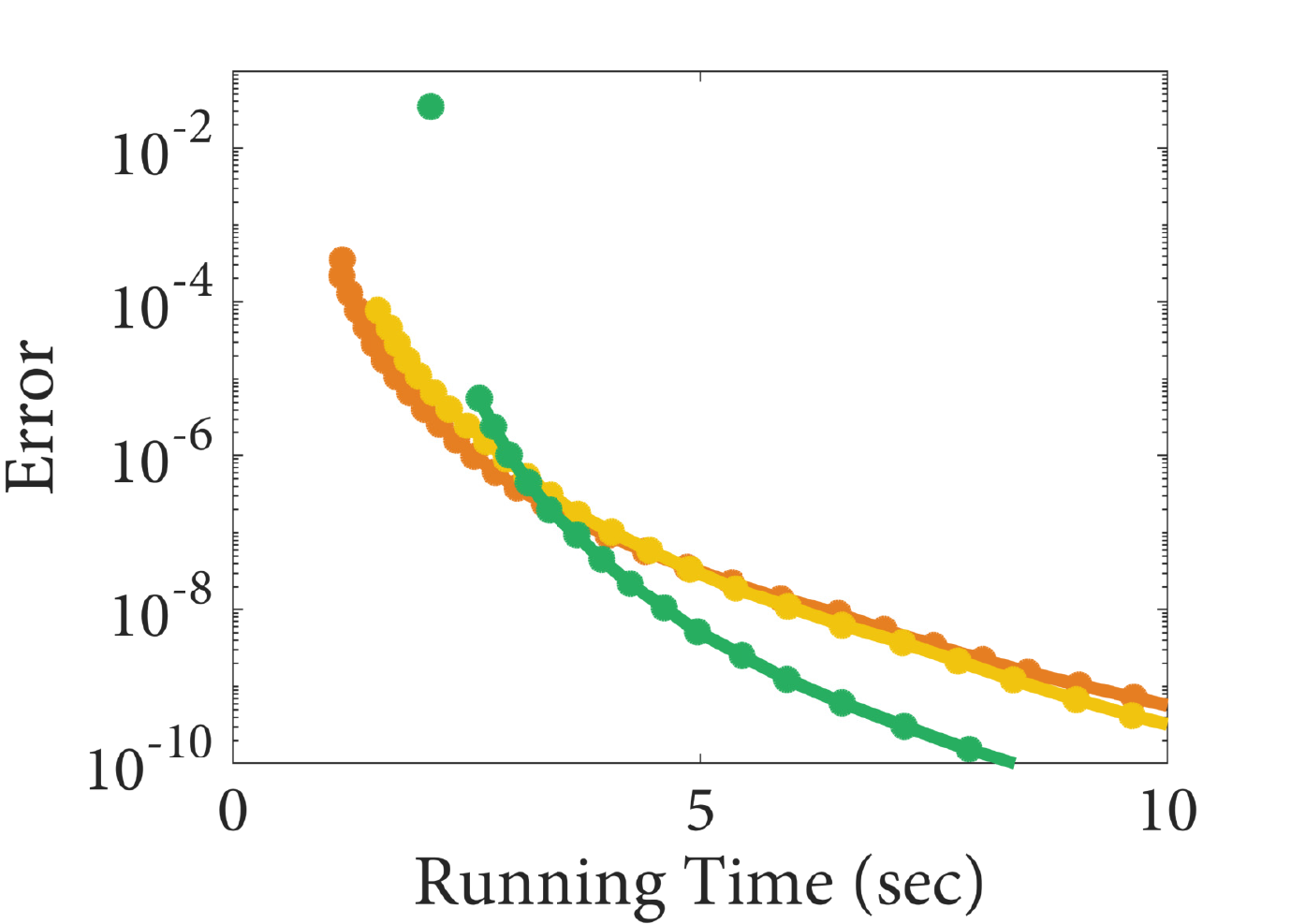} \\[1em]			
			\end{tabularx}
			
			\caption{Accuracy and precision diagrams for Burgers' equation solved using three new PBMs. The color in each plot is used to represent the number of nodes, $q$. BDF-SMFC and Adams-PMFCmj methods with $q$ nodes converged to $q$th order while Adams-SMVC methods with $q$ nodes converged to $(q+2)$th order.}
			\label{fig:Burgers_PBM}
		\end{figure}

\begin{figure}
			\centering
			\hfill
			\begin{tabular}{lllllll}
				\lineLegend{plot_red}{$q=2$} &
				\lineLegend{plot_orange}{$q=3$} &
				\lineLegend{plot_yellow}{$q=4$} &
				\lineLegend{plot_green}{$q=5$} &
				\lineLegend{plot_blue}{$q=6$} &
				\lineLegend{plot_violet}{$q=7$} &
				\lineLegend{plot_grey}{$q=8$} \hspace{1.5em} 
			\end{tabular}
			\vspace{1em}

			\centering
			\begin{tabularx}{\linewidth}{mff}
					\MethodNameTable{BBDF$_{\alpha = \tfrac{1}{2}}$} & \includegraphics[valign=m,width=0.35\textwidth]{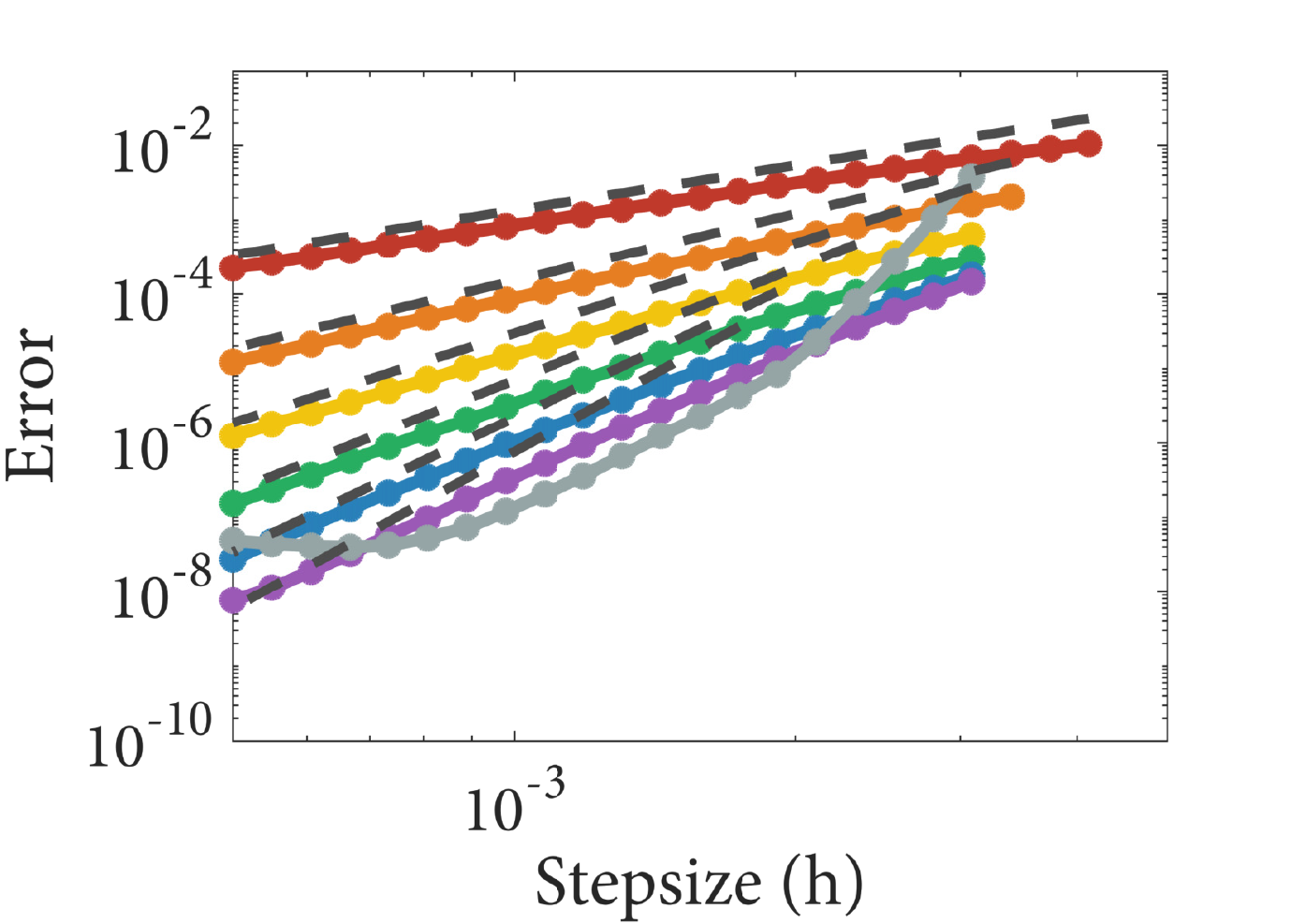} &
					\includegraphics[valign=m,width=0.35\textwidth]{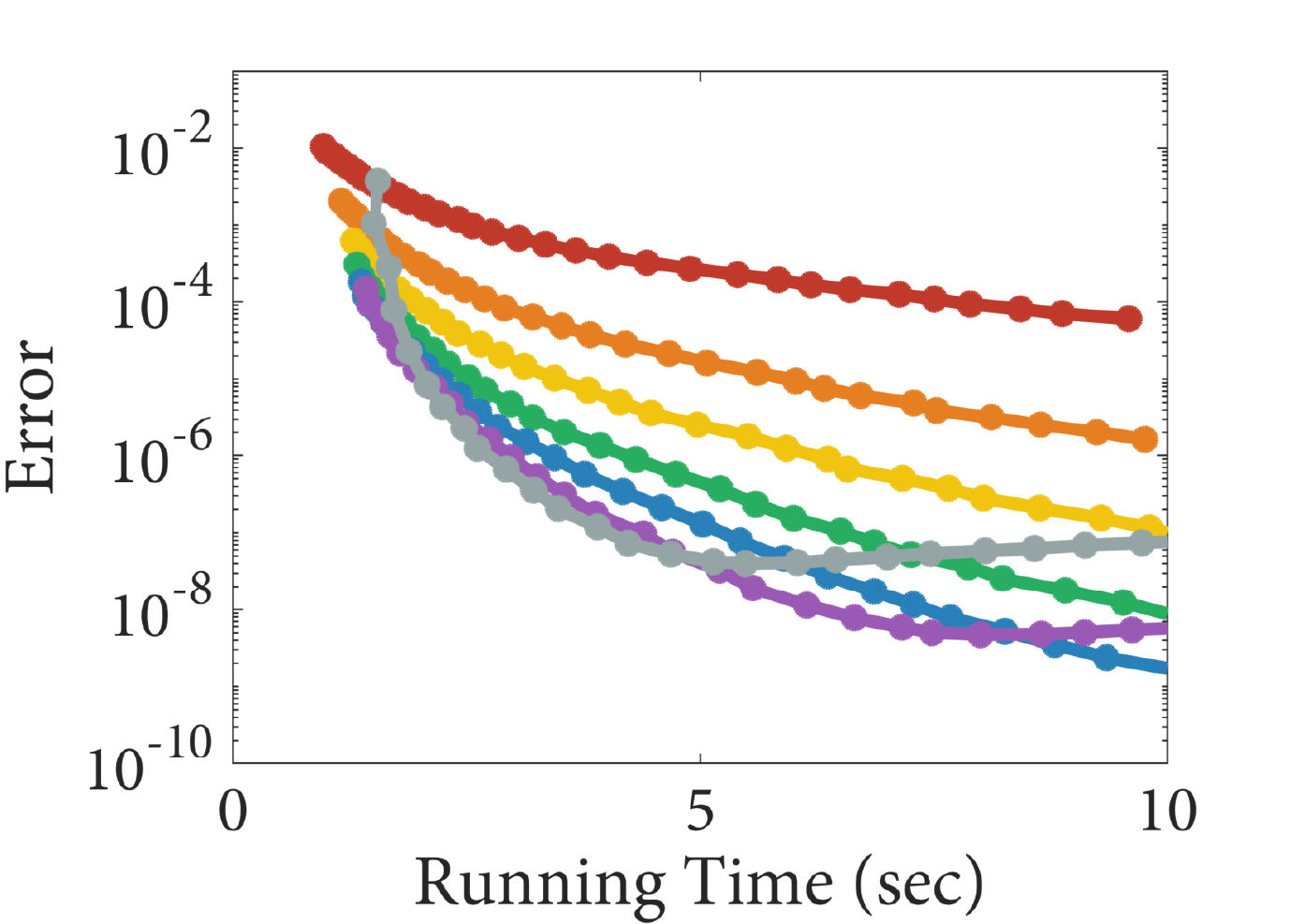} \\[1em]
					\MethodNameTable{BDF} & \includegraphics[valign=m,width=0.35\textwidth]{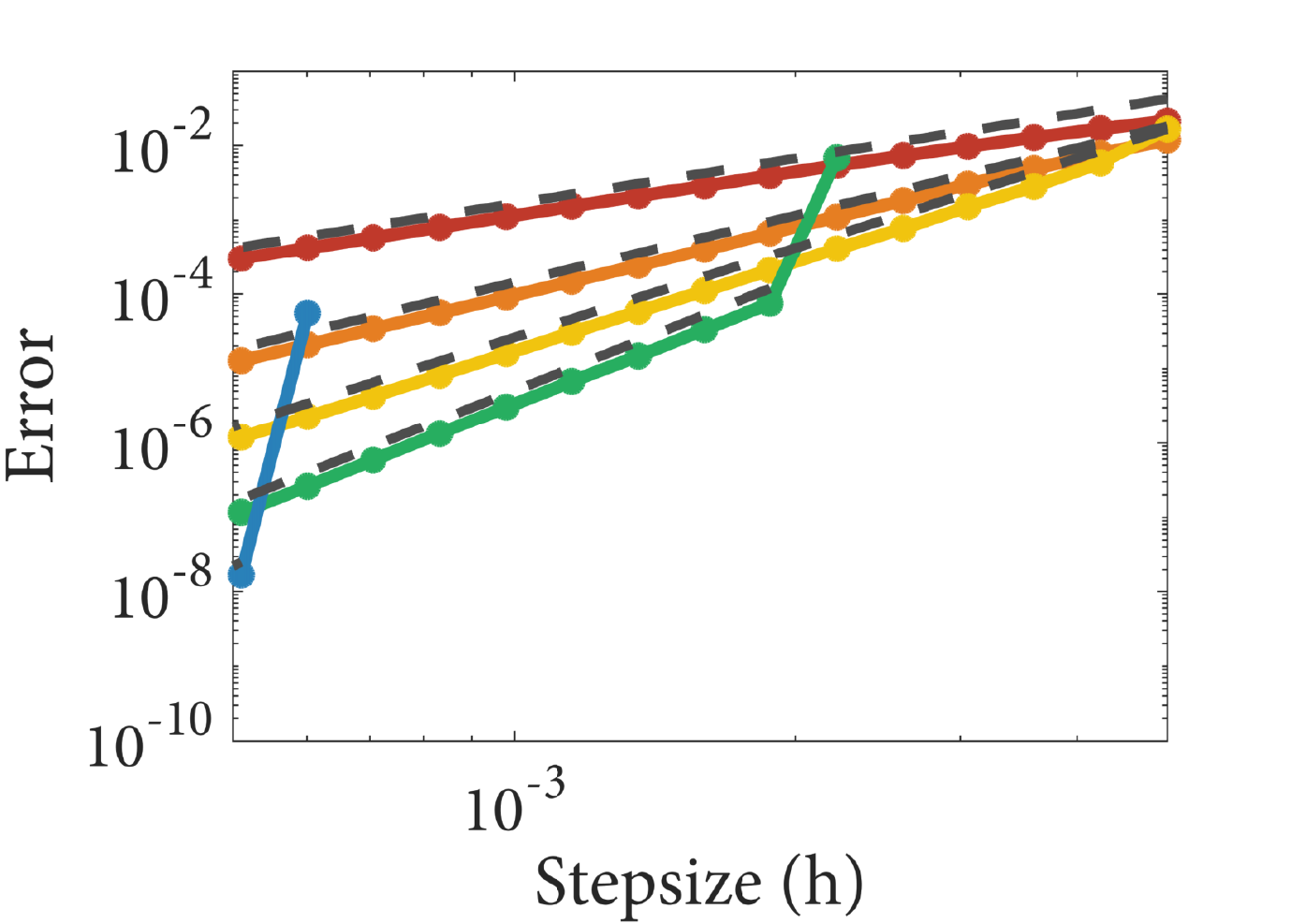}	 &
					\includegraphics[valign=m,width=0.35\textwidth]{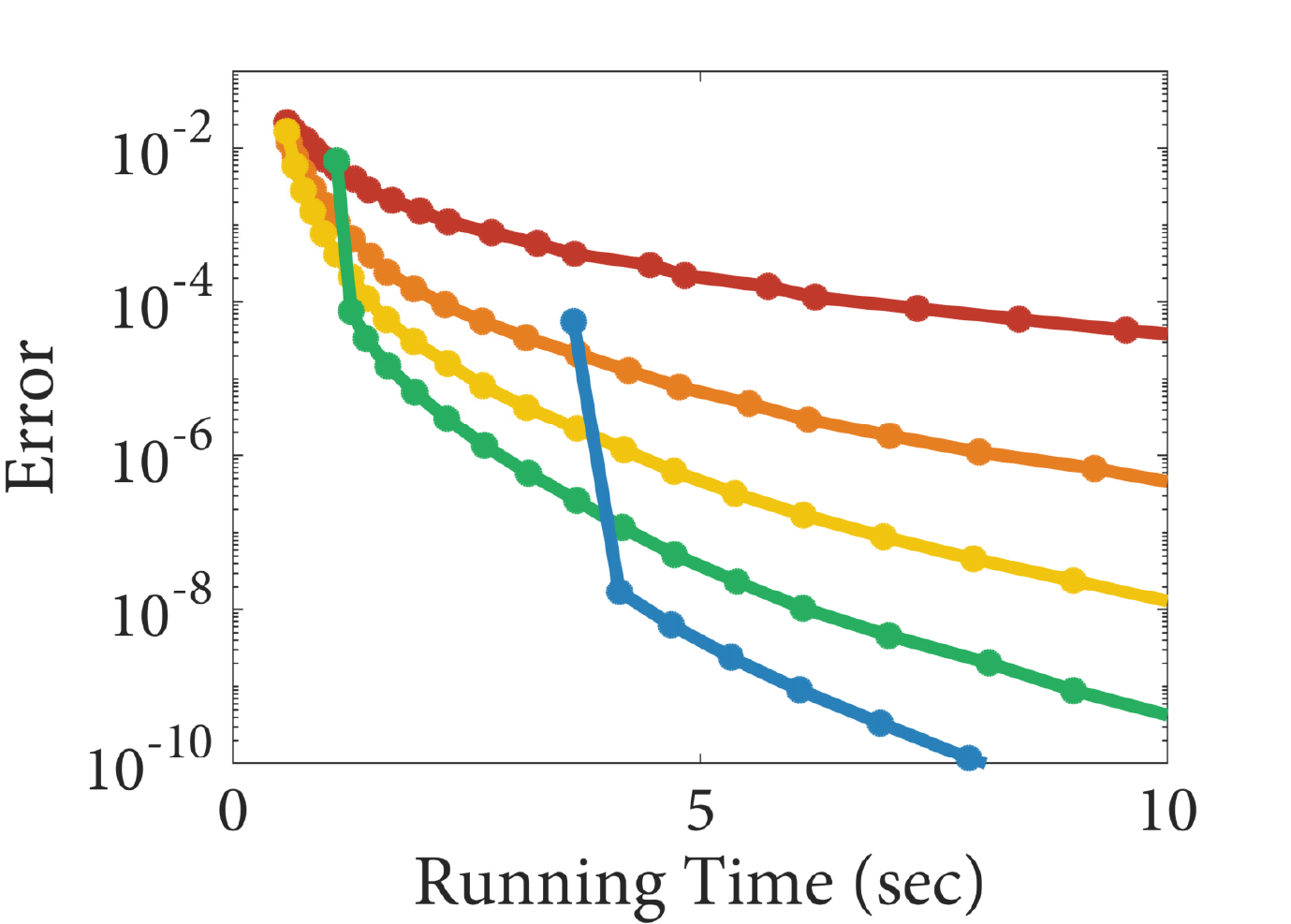} \\[1em]
					\MethodNameTable{ESDIRK} & \includegraphics[valign=m,width=0.35\textwidth]{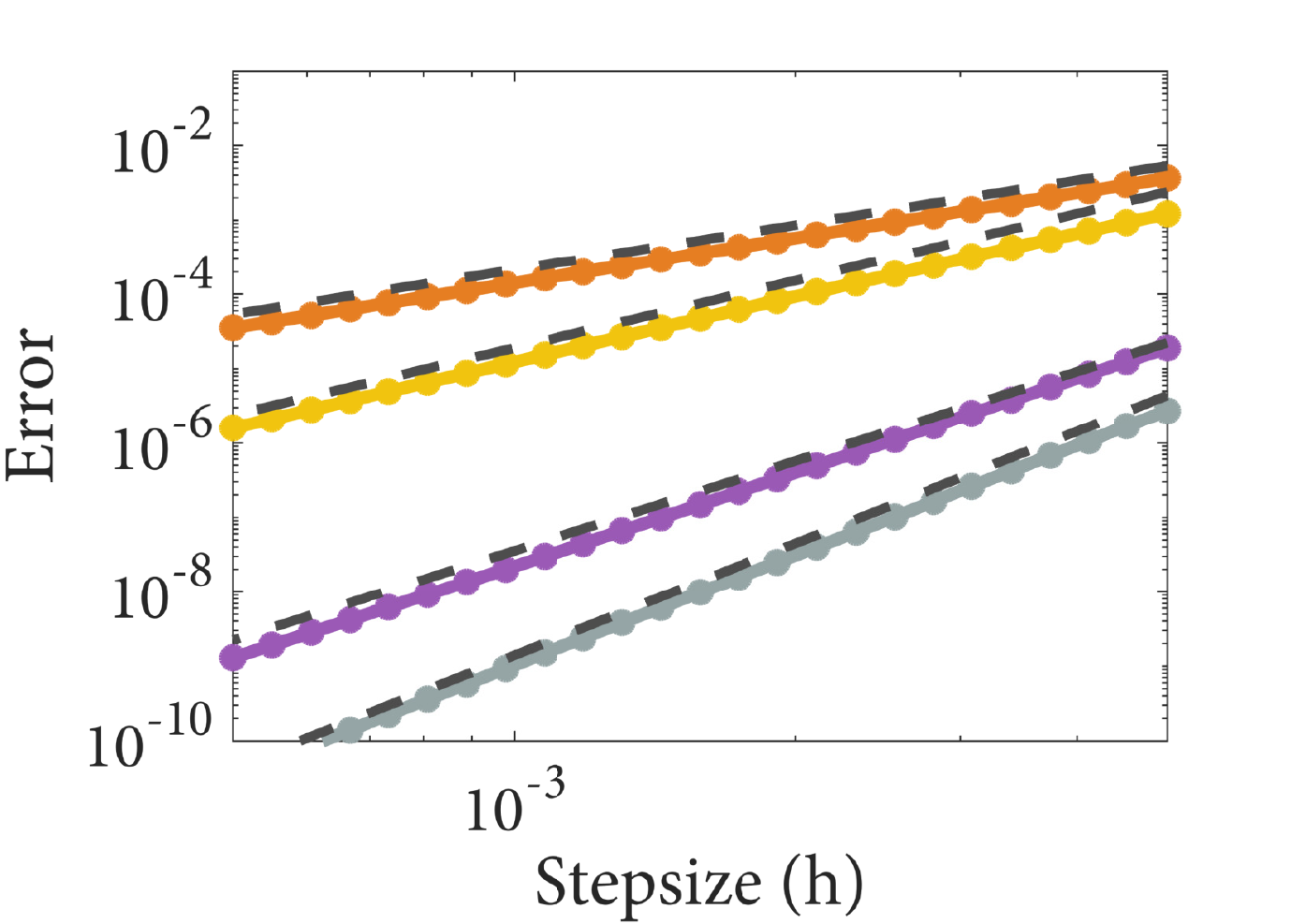}	 &
					\includegraphics[valign=m,width=0.35\textwidth]{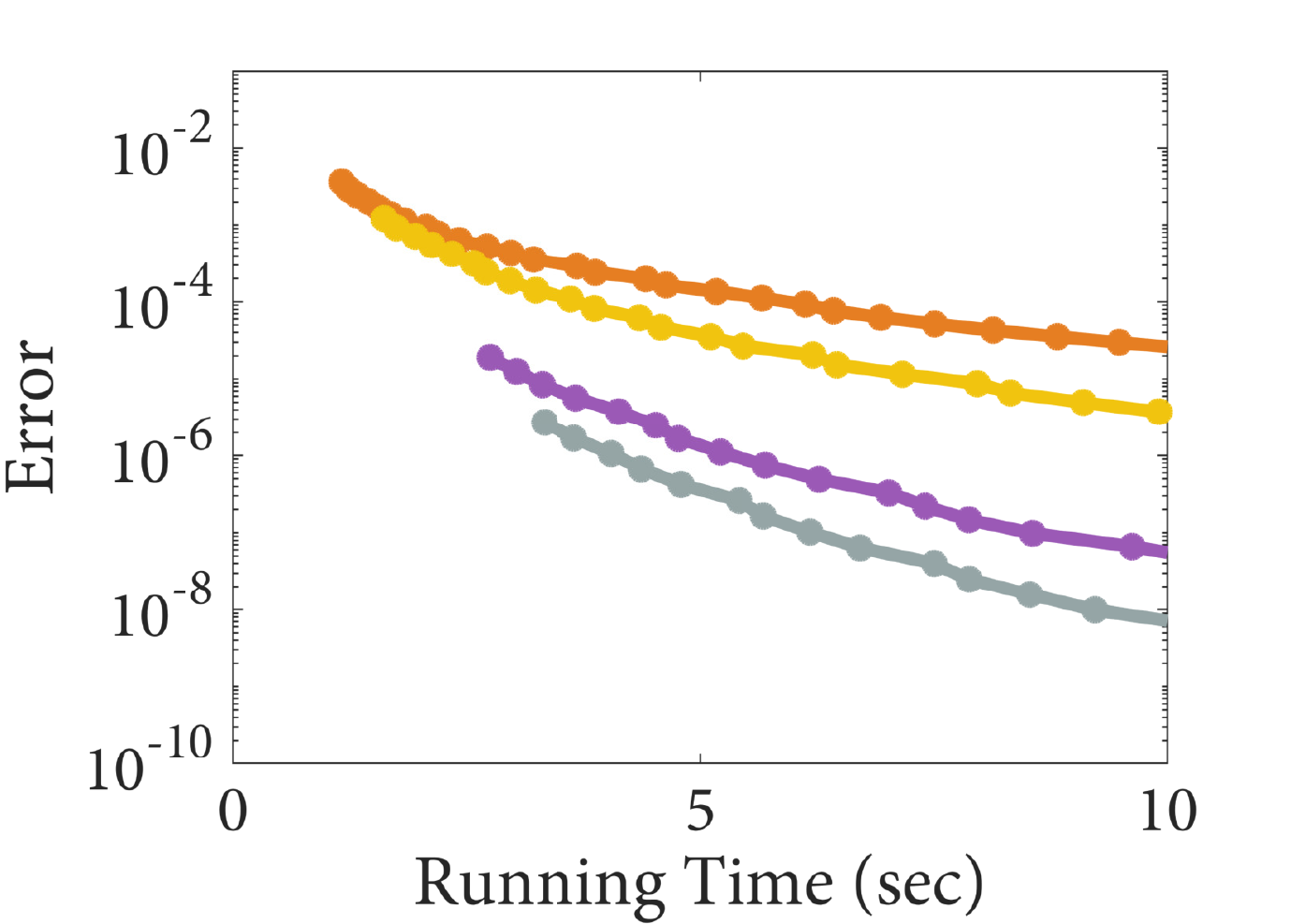} \\[1em]			
			\end{tabularx}

			\caption{Accuracy and precision diagrams for Burgers' equation solved using reference methods. The BBDF method with $\alpha=1/2$ from \cite{buvoli2019constructing} is shown on the top, classical BDF methods are in the middle, and at the bottom are ESDIRK schemes of order three through six listed in Table \ref{tab:ESDIRK_methods}. Color is used to represent the number of nodes for polynomial methods, and the number of implicit stages for ESDIRK methods.}
			\label{fig:Burgers_reference}
		\end{figure}

	\vspace{0.5em}
	\begin{enumerate}[leftmargin=*]
	
		\item  We consider the one dimensional nonlinear Schr\"{o}dinger (NLS) equation with periodic boundary conditions,
\begin{align}
	\begin{aligned}
		& iu_t + u_{xx} + u|u|^2 = 0, \\
		& u(x,t=0) =  1 + \frac{\exp(i x / 4)}{100},
	\end{aligned}
	&&
	\begin{aligned}
		x &\in[-4\pi, 4\pi], \\
		t &\in [0, 11]. \nonumber
	\end{aligned}
\end{align}
In its standard form, NLS cannot be extended into the complex plane due to the lack of analyticity caused by the absolute value in the nonlinear term. To obtain an initial value problem that satisfies the conditions of Cauchy-Kowalevski theorem \cite{ablowitz2003complex} (which allows analytic continuation into the complex time plane), we rewrite the nonlinear Schr\"{o}dinger as the following system of real-valued equations
\begin{align*}
  	a_t &= b_{xx} - (a^2b + b^3) \\
	b_t &= a_{xx} + (a^3 + b^2 a)	
\end{align*}
where $a(x,t)$ and $b(x,t)$ are the real and imaginary components of the solution $u(x,t)$. We then discretize in space using a 256 point Fourier spectral spatial discretization with no antialiasing, and test all time-integrators using twenty-five different stepsizes logarithmically spaced between \num{2.7e-2} and \num{1e-3}. In Figure \ref{fig:NLS_results} we show precision diagrams for the five methods we tested on this problem. 	
	
\item  We consider the one dimensional viscous Burgers' equation with homogenous boundary conditions \cite{buvoli2019constructing, tokman2006efficient},
\begin{align}
	\begin{aligned}
		& u_t = \nu u_{xx} - uu_x \label{eq:vburgers}, \\
		& u(x,t=0) =  \left( \sin(3 \pi x)\right)^2 \left( 1 - x\right)^{3/2},
	\end{aligned}
	&&
	\begin{aligned}
		x &\in[0, 1],\\
		t &\in [0, 1].
	\end{aligned}
\end{align}
where we take $\nu = \num{3e-4}$. We discretize in space using standard second order finite differences with 2000 gridpoints. We test all time-integrators using twenty-five different stepsizes logarithmically spaced between \num{5e-3} and \num{5e-4}. 	In Figure \ref{fig:Burgers_PBM} we show precision diagrams for the three PBMS we tested, and in Figure \ref{fig:Burgers_reference} we show results for BDF, ESDIRK \cite{kennedy2019diagonally} and BBDF.  	

\end{enumerate}

		\begin{table}
			\renewcommand\arraystretch{1.5}
			\centering
			\begin{tabular}{lll}
				Name & Order & NOI \\	\hline
				ESDIRK3(2)4L$[2]$SA 			& 3 & 3 \\
				ESDIRK4(3)7L$[2]$SA 			& 4 & 6 \\
			\end{tabular}
			\hspace{1em}
			\begin{tabular}{lll}
				Name & Order & NOI \\	\hline
				ESDIRK5(4)8$[2]$SA 			& 5 & 7 \\
				ESDIRK6(5)9$[2]$SA 			& 6 & 8 \\
			\end{tabular}
			
			\caption{ESDIRK methods from \cite{kennedy2019diagonally} that are used in our numerical experiments. The acronym NOI abbreviates number of implicit stages. Note that we always use an exact Jacobian so each stage requires a full nonlinear solve.}
			\label{tab:ESDIRK_methods}
		\end{table}

		\begin{table}[h]
			
			\renewcommand\arraystretch{1.5}
			\begin{small}
			\begin{tabular}{lllll}
				{\bf Family} & {\bf AII/AOI} & {\bf Nodes} & {\bf EP} & {\bf Extrapolation Factor }$\alpha$ \\
				BDF & SMVC & iEqui & n.a. & 
				\begin{tabular}{llll}
		        	$q=2,3$ & $q=5$ & $q=6$ & $q=7$ \\ 
		        	0.75 & 0.13 & 0.18 & 0.15
				\end{tabular} \\ \hline
				Adams & PMFCmj & iCheb & FI, $\ell = 3$ & 
				\begin{tabular}{llll}
		        	$q=3$ & $q=4$ & $q=5$ \\ 
		        	0.32 & 0.20 & 0.136
				\end{tabular} \\ \hline			
				Adams & SMFCmj & iCheb & FI, $\ell = 3$ & 
				\begin{tabular}{llll}
		        	$q=3$ & $q=4$ & $q=5$ \\ 
		        	0.75 & 0.50 & 0.33
				\end{tabular} \\ \hline
				Adams & SMFCmj & iCheb & VI & 
				\begin{scriptsize}
				\hspace{-1em} \begin{tabular}{lllll}
		        	$q=2,3,4$ & $q=5$ & $q=6$ & $q=7$ & $q=8$ \\ 
		        	0.75 & 0.45 & 0.32 & 0.22 & 0.18 \\
				\end{tabular}
				\end{scriptsize}				
			\end{tabular}
			\end{small}
				
			\caption{PBMs and extrapolation factors used in the NLS experiment.}
			\label{tab:NLS_methods}

		\end{table}

		\begin{table}[h]
			
			\renewcommand\arraystretch{1.5}
			\begin{small}
			\begin{tabular}{lllll}
				{\bf Family} & {\bf AII/AOI} & {\bf Nodes} & {\bf EP} & {\bf Extrapolation Factor }$\alpha$ \\ 
				BDF & SMVC & iEqui & n.a. & 
				\begin{tabular}{llll}
		        	$q=2,3,\ldots, 8$ \\
		        	0.50
				\end{tabular} \\ \hline
				Adams & PMFCmj & iEqui & VI & 
				\begin{tabular}{llll}
		        	$q=2$ & $q=3$ & $q=4$ & $q=5$ \\ 
		        	3.00 & 1.25 & 0.64 & 0.43
				\end{tabular} \\ \hline
				Adams & SMVC & iCheb & FI, $\ell = 2$ & 
				\begin{tabular}{llll}
		        	$q=3$ & $q=4$ & $q=5$ \\
		        	0.90 & 0.70 & 0.45
				\end{tabular}
			\end{tabular}
			\end{small}
				
			\caption{PBMs and extrapolation factors used in the Burgers' experiment.}
			\label{tab:BRG_methods}

		\end{table}

		\subsection{Discussion}
		
		We split our discussion in three parts. First we discuss PBMs for dispersive equations, then we discuss PBMs for dissipative equations, and finally we conclude with some general comments regarding all the newly introduced PBMs.
		
		\subsubsection{Methods with A($\mathbf{90\degree}$) stability}
		
		The results from the NLS experiment validate linear stability analysis and demonstrate that diagonally implicit PBMS can be used to solve dispersive wave equations. In order to achieve A($90\degree$) stability each method had to be run with a specific $\alpha$ value. Overall we see that parallel methods are the fastest for obtaining the solution to high precision, but limited analyticity in the solution paired with small $\alpha$ requirments prevented convergence at coarse timesteps. On the other hand, serial methods allow for larger $\alpha$ and are able to achieve good efficiency at coarse timesteps, and their improved error constants make them only marginally less efficient than their parallel counterparts at fine timesteps. 
				
		{\em 1. BDF methods:} To achieve $A(90\degree)$ stability with BDF or GBDF methods it is necessary to select a serial method construction. High-order integrators also require a small extrapolation factor $\alpha$, which caused instabilities at coarse time-steps due to limited analyticity in the solution. Despite this limitation, the BDF-SMCF methods with $q=6$ were able to obtain the solution to machine precision with nearly the same efficiency as the highly-optimized sixth order EDIRK method. If only moderate precision is required, then the BDF-SMCF and GBDF-SMVC methods with $q=3$ demonstrated fourth-order convergence and better efficiency than the fourth order EDIRK4 method, while also requiring fewer solution vectors per timestep. Overall, all BDF methods with $2 \le q \le 6$ demonstrated or exceeded their expected orders of accuracy. High-order integrators with $q > 6$ converged rapidly to machine precision and the experiment would have to be repeated using variable precision to see proper convergence curves for these integrators. 
		
		 {\em 2. Adams methods:} $A(90\degree)$ stability for Adams methods can be achieved with both parallel and serial construction strategies. In addition, Adams methods allow for larger extrapolation factors than their GBDF counterparts, enabling them to avoid instabilities at coarse timesteps. Amongst all methods we tested, the parallel Adams PMFCmj scheme with $q=6$ was the fastest at obtaining the solution to machine precision, followed by the sixth order EDIRK method, and the serial Adams SMFCmj schemes with $q=6$. Larger $\alpha$ values allow the serial Adams SMFCmj to converge at coarser timesteps than the Adams PMFCmj scheme; however, for methods with $q=7$ the larger $\alpha$ values cause increased susceptibility to round-off errors as evidenced by the higher noise floor of the SMFCmj method.
			
		Based on our results, we broadly classify the best performing methods with A($90\degree$) stability into the following four cases:
		\begin{center}
			\vspace{1em}
			\renewcommand\arraystretch{1.5}
			\begin{tabular}{l|ll}
			
				& {\em Serial} & {\em Parallel} \\ \hline
				{\em High-Accuracy} & SMFCmj, $q=5,6$ & PMFCmj, $q=6$ \\ 
				{\em Low-Accuracy}	 & BBDF-SMVC $q=3$ & PMFCmj, $q=5$ \\
				 & Adams-SMFCmj $q=2,3,4$ & 
			\end{tabular}
			\vspace{1em}
		\end{center}

		\subsubsection{Methods with A($\theta > 0$) stability}
		
		In our numerical experiment with Burgers' equation we compare a serial BDF PBM, a serial Adams PBM, and a parallel Adams PBM. A key feature of all of the methods is that they all possess  unbounded stability regions. To compare the new PBMs against existing integrators we also integrate in time using the BBDF method \cite{buvoli2019constructing}, BDF, and ESDIRK methods \cite{kennedy2019diagonally}. Before discussing the results, we note that the second fastest method is the sixth order BDF due to its high-accuracy and single nonlinear solve per timestep. However this result is slightly misleading since the method has a poor stability angle. In fact, by simply increasing the spatial resolution of our test problem we can render both BDF6 and BDF5 unstable. However, we use the current spatial grid so that BDF6 can serve as a benchmark for excellent performance.
		
		{\em 1. BDF-SMFC.} This PBM is a serial BDF method with significantly improved error properties and improved susceptibility to round-off errors compared to the parallel BBDF method. BDF-SMFC with $q=8$ is the fastest method to achieve machine precision significantly outperforming all parallel PBMs, ESDIRK methods, and BDF6. Moreover, the stability angle for BDF-SMFC with $q=8$ and $\alpha = 1/2$ is approximately $81.7 \degree$ compared to the $17.8\degree$ for BDF6.
				
		{\em 2. Adams-PMFCmj.} This PBM is a parallel Adams method that is nearly identical to the BAM methods originally presented in \cite{buvoli2019constructing}. However, by switching to the PMFCmj construction (as opposed to PMCF used for BAM) we can obtain a parallel Adams method with an unbounded stability region. In contrast the BAM methods possess a large but bounded stability regions that render them unsuitable for highly stiff equations. Moreover, these new methods have improved error constants compared to the BBDF method so that the Adams-PMFCmj with $q=5$ achieves identical efficiency to the BBDF$_{\alpha=1/2}$ method with $q=7$ or $q=8$. Low order Adams-PMFCmj methods can be run with larger $\alpha$ values and therefore retain stability across a wider range of timesteps. Finally, it is also possible to construct Adams-PMFCmj with $q=6,7,8$, however they are not competitive since they require very small $\alpha$ (e.g. 0.1) and stepsizes to remain stable on our test problem.
		
		{\em 3. Adams-SMVC.} This PBM is a serial Adams method with excellent error properties. Adams-SMVC with $q=3$ and $\alpha=0.9$ is able to outperform all Adams-PMFCmj integrators, all BDF-SMFC with $q \le 6$, and all ESDIRK methods. At coarsest time steps, Adams-PMFCmj methods with $q=2,3,4$ are still slightly faster due to parallelism, but the advantage quickly ends as stepsize decreases. Adams-SMVC with $q=5$ requires a smaller $\alpha$ and is unstable at coarse timesteps, but is even more efficient at fine timesteps, rivaling the performance of BDF-SMFC with $q=7$.

		\subsubsection{General remarks}
		
		Our results demonstrate that we where able to successfully construct PBMs for solving dispersive equations. This opens the possibility of using polynomial based methods in ways that were not not possible with traditional polynomial-based implicit linear multistep methods. One particularly interesting new method with A($90\degree$) is the parallel Adams PMFCmj method that offers imaginary stability despite only requiring one serial nonlinear solve per timestep. On problems where the linear systems can be properly preconditioned this method could offer significant speeds up compared to current approaches.
		
		The second aim of this work was to develop efficient serial PBMs that can be used in replacement of parallel PBMs when additional processing units are unavailable, or when parallel communication costs are no longer negligible compared to the nonlinear solves. Interestingly, many of the serial methods we introduced are often more efficient than their parallel counterparts, especially if high accuracy is desired. For example, provided that one carefully chooses the order depending on their accuracy requirements, the BDF-SMCV scheme will be more efficient than the parallel BBDF scheme.
		
		Looking beyond our two numerical experiments, the efficiency of polynomial PBMs will depend on many factors that are specific to the problem being solved. These include the choice of linear solver and preconditioner, as well as the characteristics of the computer or cluster where the methods are run. Another key importance will be the impact of complex arithmetic which may vary depending on how it is implemented. For details on effects of complex arithmetic for our experiment see the related comments in \cite{buvoli2019constructing} since the numerical experiments in this work are conducted in an identical fashion. 
				
\section{Conclusion}

	In this work we present another key component of the polynomial framework, namely the ability to construct a vast range of integrators with variable architecture and order-of-accuracy. Moreover, these new method constructions are all geometrically motivated by considering interpolating polynomials in the complex plane and do not require discussion of nonlinear order conditions. This presents a different way to approach time-integrator construction that should be more accessible to a wider range of practitioners who are already familiar with spatial finite difference methods. The newly proposed method constructions also create new possibilities to use polynomial block methods including for solving dispersive wave equations, and on machines where parallelism is not feasible. Finally, we also paved a future direction for us to develop optimized PBMs with special nodes and extrapolation factors that are tuned to satisfy certain stability requirements or to minimize analyticity requirements off the real line. 
	
\section*{Acknowledgements} I would like to thank Randy J. LeVeque and Mayya Tokman for their many fruitful discussions during the development of this work. This research was supported in part by funding from the Applied Mathematics Department at the University of Washington and NSF grant DMS-1216732.

		\bibliographystyle{siamplain}
\bibliography{references_nourl/complex_analysis,references_nourl/exponential_integrators,references_nourl/finite_difference_methods,references_nourl/general_linear_methods,references_nourl/implicit_explicit_integrators,references_nourl/order_conditions,references_nourl/rosenbrock_methods,references_nourl/spectral_deferred_correction,references_nourl/spectral_methods,references_nourl/time_integration,references_nourl/buvoli_papers.bib}

\end{document}

\maketitle

\section{Example polynomial diagrams}
\label{sup_sec:example_diagrams}

We show two example polynomial diagrams generated using PIPack along with the Matlab code used to produce the figures and generate method coefficients.

\subsection{Polynomial diagram for BDF SMFC} 
 
The BDF SMFC method is a serial method that can be used to solve both dissipative equations and dispersive equations (if $\alpha$ is chosen to be sufficiently small). Its polynomial diagrams is:

\begin{center}
	\includegraphics[trim={1.5cm 5.5cm 1.5cm 5cm},clip]{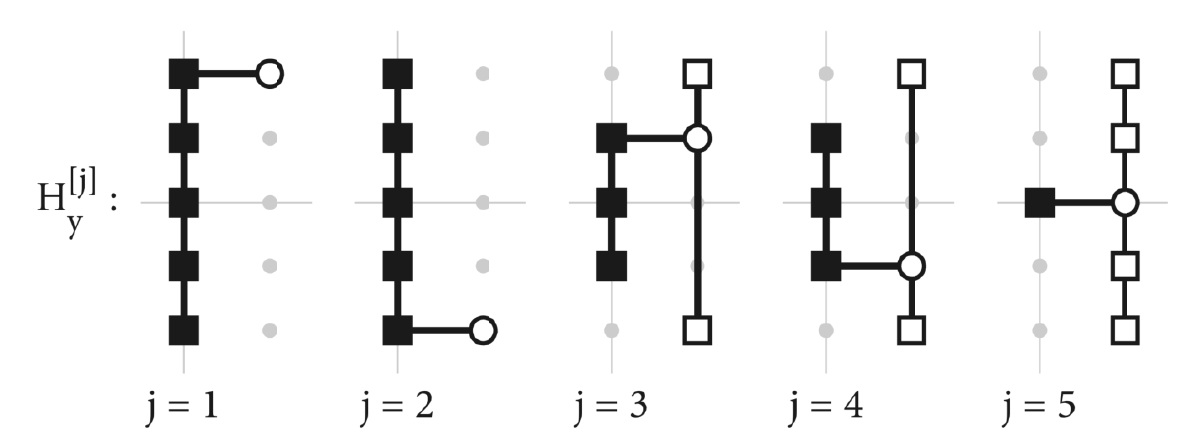}	
\end{center}
The code for generating this figure is shown below. PIPack must be added to Matlab's path for this code to run. We also include code to compute the method coefficients.

\vspace{1em}
\begin{verbatim}
% -- initialize Adams method generator -------------------------------
BMG = PBMGenerator(struct( ...
    'NodeSet_generator', IEquiNSG('inwards', 'double'), ... 
    'ODEPoly_generator', BDF_PG('SMFO', 'diagonally_implicit') ...
));

% -- generate method with 5 nodes ------------------------------------
bdf_method = BMG.generate(5);
alpha = 0.5;

fh = bdf_method.polynomialDiagram( struct('show_numbering', false));
exportFigure(fh, struct('SavePath', 'bdf_smfc_polydiagram'))

% -- generate method coefficients ------------------------------------
% compact_traditional produces a method of the form:
%   y^{[n+1]} = A y^{[n]} + h B f^{[n]} + h C f^{[n+1]}
[A, B, C] = bdf_method.blockMatrices(alpha, 'compact_traditional');


\end{verbatim}

\subsection{Polynomial Diagrams for Adams PMFCmj with VI endpoints}

The Adams PMFCmj method is a parallel method that can be used to solve dissipative equations. The polynomial diagrams for Adams methods show the polynomials $L_y(\tau)$ and $L_F(\tau)$:  
\begin{center}
	\includegraphics[trim={1.5cm 5.5cm 1.5cm 5cm},clip]{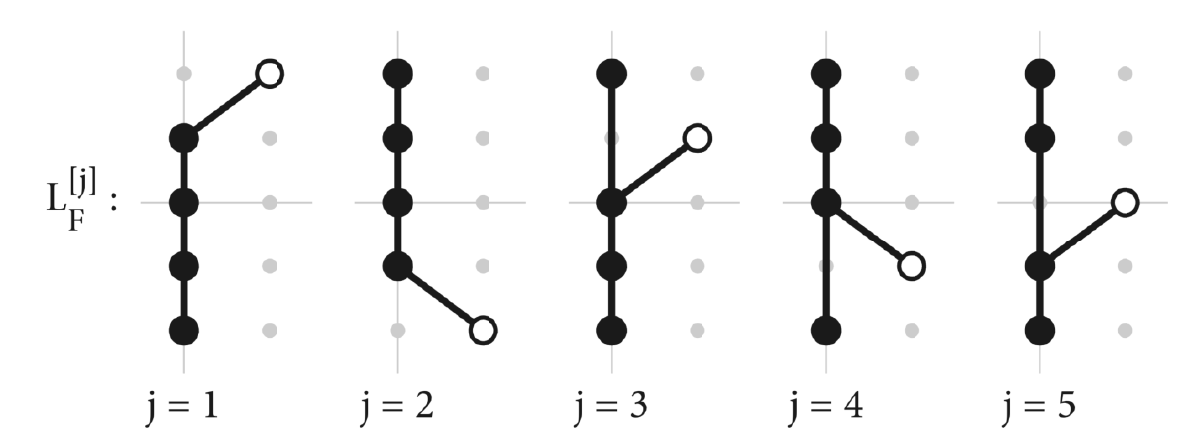}
	\includegraphics[trim={1.5cm 5.5cm 1.5cm 5cm},clip]{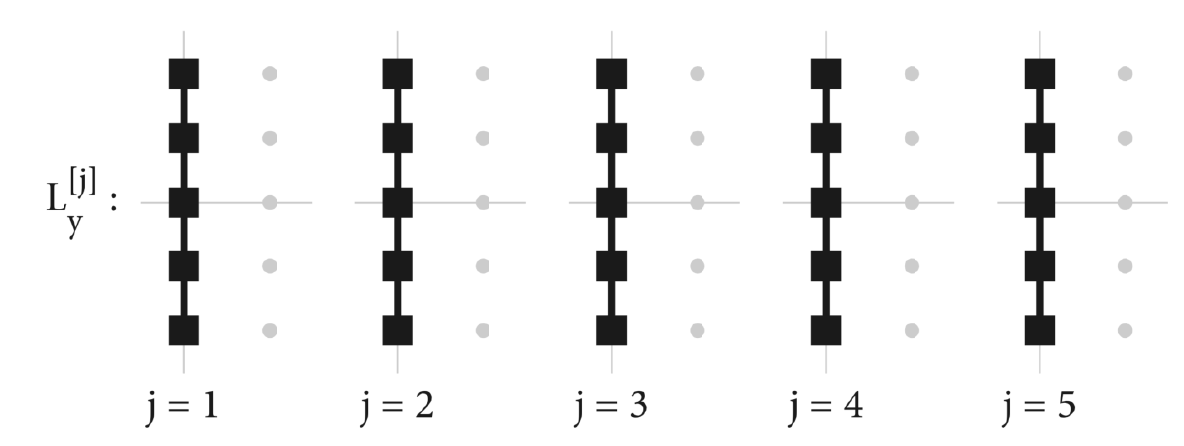}
\end{center}
Since this is an Adams methods, it also has an endpoint diagram:
\begin{center}
	\includegraphics[trim={6cm 2cm 5cm 1cm},clip,width=0.12\textwidth,valign=m]{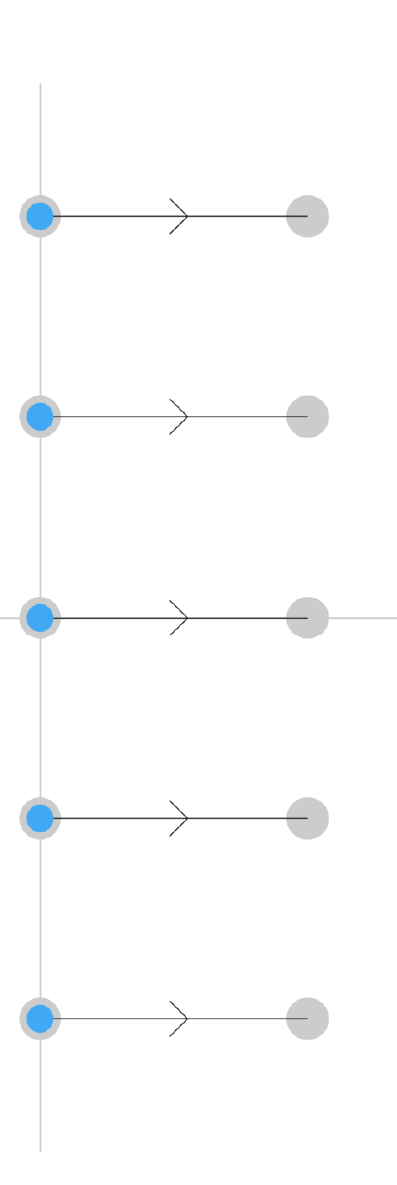}
\end{center}

\vspace{1em}
\noindent The code for generating this figure is shown on the following page. We also include code for computing the method coefficients.

\newpage
\begin{verbatim}
% -- initialize Adams method generator -------------------------------
AMG = PBMGenerator(struct( ...
    'NodeSet_generator', IEquiNSG('inwards', 'double'), ...
    'ODEPoly_generator', Adams_PG( ...
         'PMFO', 'PMFOmj', ...
         VariableInputEPG(), ...
         diagonally_implicit') ...
)); 

% -- generate method with 5 nodes ------------------------------------
adams_method = AMG.generate(5);
alpha = 0.5;

[Ly_fh, ~, LF_fh] = adams_method.polynomialDiagram( ...
     struct('show_numbering', false));
exportFigure(Ly_fh,  struct('SavePath', 'adams_PMFCmj_Ly_polydiagram'))
exportFigure(LF_fh,  struct('SavePath', 'adams_PMFCmj_LF_polydiagram'))

fh = adams_method.expansionPointDiagram();
exportFigure(fh,  struct('SavePath', 'adams_PMFCmj_endpointdiagram'))	


% -- generate method coefficients ------------------------------------
% compact_traditional produces a method of the form:
%   y^{[n+1]} = A y^{[n]} + h B f^{[n]} + h C f^{[n+1]}
[A, B, C] = bdf_method.blockMatrices(alpha, 'compact_traditional');


\end{verbatim}